\documentclass[11pt, leqno]{amsart}
\usepackage{amscd, amsmath, amssymb, amsfonts}
\usepackage{latexsym}
\usepackage{mathrsfs}

\newtheorem{Theorem}{Theorem}[section]
\newtheorem{Proposition}[Theorem]{Proposition}
\newtheorem{Lemma}[Theorem]{Lemma}
\newtheorem{Corollary}[Theorem]{Corollary}

\newtheorem{TheoremNoNum}{Theorem}[Theorem]
\newtheorem{PropositionNoNum}{Proposition}[Theorem]
\newtheorem{CorollaryNoNum}{Corollary}[Theorem]

\theoremstyle{definition}
\newtheorem{Definition}[Theorem]{Definition}
\newtheorem{Remark}[Theorem]{Remark}
\newtheorem{Example}[Theorem]{Example}

\newtheorem{Question}[Theorem]{Question}

\renewcommand{\theTheorem}{\arabic{section}.\arabic{Theorem}}
\renewcommand{\theClaim}{\arabic{section}.\arabic{Theorem}.\arabic{Claim}}
\renewcommand{\theequation}{\arabic{section}.\arabic{equation}}

\newcommand{\KK}{{\mathbb{K}}}

\newcommand{\ZZ}{{\mathbb{Z}}}
\newcommand{\QQ}{{\mathbb{Q}}}
\newcommand{\RR}{{\mathbb{R}}}

\newcommand{\PP}{{\mathbb{P}}}
\newcommand{\OO}{{\mathcal{O}}}

\newcommand{\Xscr}{{\mathscr{X}}}
\newcommand{\Dscr}{{\mathscr{D}}}
\newcommand{\Escr}{{\mathscr{E}}}

\newcommand{\Ucal}{{\mathcal{U}}}
\newcommand{\Vcal}{{\mathcal{V}}}

\newcommand{\Div}{\operatorname{Div}}

\newcommand{\Image}{\operatorname{Image}}

\newcommand{\ord}{\operatorname{ord}}

\newcommand{\rest}[2]{\left.{#1}\right\vert_{{#2}}}  
\renewcommand{\setminus}{\smallsetminus}

\newcommand{\Spec}{{\operatorname{Spec}}}

\newcommand{\outdeg}{\operatorname{outdeg}}

\newcommand{\id}{\operatorname{id}}

\newcommand{\Supp}{\operatorname{Supp}}
\newcommand{\alg}{\operatorname{alg}}

\newcommand{\Ob}{\operatorname{Ob}}

\newcommand{\Proof}{{\sl Proof.}\quad}
\newcommand{\QED}{{\unskip\nobreak\hfil\penalty50\quad\null\nobreak\hfil
{$\Box$}\parfillskip0pt\finalhyphendemerits0\par\medskip}}

\newcommand{\Def}{\operatorname{Def}}
\newcommand{\Spf}{\operatorname{Spf}}

\newcommand{\val}{\operatorname{val}}
\newcommand{\Isom}{\operatorname{Isom}}

\newcommand{\Hom}{\operatorname{Hom}}

\begin{document}

\title[Rank of divisors under specialization]{Rank of divisors on 
hyperelliptic curves and graphs under specialization}
\author{Shu Kawaguchi}
\address{Department of Mathematics, Graduate School of Science,
Kyoto University, Kyoto 606-8502, Japan}
\email{kawaguch@math.kyoto-u.ac.jp}
\author{Kazuhiko Yamaki}
\address{Institute for Liberal Arts and Sciences, 
Kyoto University, Kyoto, 606-8501, Japan}
\email{yamaki@math.kyoto-u.ac.jp}
\thanks{The first named author partially supported by KAKENHI 24740015, and the second named author partially supported by KAKENHI 21740012.}

\newcommand{\Proj}{\operatorname{\operatorname{Proj}}}
\newcommand{\Prin}{\operatorname{Prin}}
\newcommand{\Rat}{\operatorname{Rat}}
\newcommand{\zero}{\operatorname{div}}

\newenvironment{notation}[0]{%
  \begin{list}%
    {}%
    {\setlength{\itemindent}{0pt}
     \setlength{\labelwidth}{4\parindent}
     \setlength{\labelsep}{\parindent}
     \setlength{\leftmargin}{5\parindent}
     \setlength{\itemsep}{0pt}
     }%
   }%
  {\end{list}}

\newenvironment{parts}[0]{%
  \begin{list}{}%
    {\setlength{\itemindent}{0pt}
     \setlength{\labelwidth}{1.5\parindent}
     \setlength{\labelsep}{.5\parindent}
     \setlength{\leftmargin}{2\parindent}
     \setlength{\itemsep}{0pt}
     }%
   }%
  {\end{list}}
\newcommand{\Part}[1]{\item[\upshape#1]}

\begin{abstract}
Let $(G, \omega)$ be a hyperelliptic vertex-weighted graph of genus $g \geq 2$. 
We give a characterization of $(G, \omega)$ for which there exists a smooth projective curve $X$ of genus $g$ over a complete discrete valuation field with reduction graph $(G, \omega)$ such that the ranks of any divisors are preserved under specialization. 
We explain, for a given vertex-weighted graph $(G, \omega)$ in general,  
how the existence of such $X$ relates the Riemann--Roch 
formulae for $X$ and $(G, \omega)$, and also how the existence of such $X$ 
is related to a conjecture of Caporaso. 
\end{abstract}

\maketitle

\section{Introduction and statements of the main results}
\subsection{Introduction}
The theory of divisors on smooth projective curves has been actively and deeply studied since the nineteenth century (cf. \cite{ACGH, ACGH2}).  It has been found that, also on graphs, there exists a good theory of divisors (including such notions as linear systems, linear equivalences, canonical divisors, degrees, and ranks).  A Riemann--Roch formula, one of the most important formulae in the theory of divisors, was established by Baker and Norine on finite loopless graphs in their foundational 
paper~\cite{BN}. A Riemann--Roch formula on tropical curves was independently proved by Gathmann and Kerber~\cite{GK} and Mikhalkin and Zharkov~\cite{MZ}. Further, a Riemann--Roch formula on vertex-weighted graphs was proved by Amini and Caporaso \cite{AC}, and on metrized complexes by Amini and Baker~\cite{AB}. 

As Baker \cite{B} revealed, the above parallelism between the theory of divisors  on curves and that on graphs is not just an analogy. Let  $\KK$ be a complete discrete valuation field with ring of integers $R$ and algebraically closed residue field $k$. Let $X$ be a geometrically irreducible smooth projective curve over $\KK$. 
An $R$-curve means an integral scheme of dimension $2$ that is projective and 
flat over $\Spec(R)$. A {\em semi-stable model} of $X$ is an $R$-curve $\Xscr$  
whose generic fiber is isomorphic to~$X$ and whose 
special fiber is a reduced scheme with at most nodes  
as singularities. For simplicity, suppose that there exists a semi-stable model $\Xscr$ of $X$ over $\Spec(R)$. 
Let $(G, \omega)$ be the (vertex-weighted) reduction graph 
of $\Xscr$, where $G$ is the dual graph of the special fiber of $\Xscr$ with natural vertex-weight function $\omega$ on $G$ (see  \S\ref{sec:prelim} for details). 
Let $\Gamma$ be the metric graph associated to  $G$, where each edge of $G$ is assigned length $1$. 
To a point $P \in X(\KK)$, one can naturally assign  
a vertex $v$ of $G$. This assignment is called the specialization map, 
and extends to  $\tau: X(\overline{\KK}) \to \Gamma_\QQ$, where 
$\overline{\KK}$ is a fixed algebraic closure of $\KK$ and  $\Gamma_\QQ$ is the set of points on~$\Gamma$ whose distance from every vertex of $G$ is rational. 
Let $\tau_*: \Div(X_{\overline\KK}) \to \Div(\Gamma_\QQ)$ be 
the induced map on divisors, and let $r_{X}$ (resp. $r_\Gamma$, 
$r_{(\Gamma, \omega)}$) denotes 
the rank of divisors on $X$ (resp. $\Gamma$, $(\Gamma, \omega)$) (see \S\ref{sec:prelim} for details). 
In \cite{B}, Baker showed that 
$r_{\Gamma} (\tau_*(\widetilde{D})) \geq r_{X}(\widetilde{D})$ for any 
$\widetilde{D} \in \Div(X_{\overline{\KK}})$, 
a result now called Baker's Specialization Lemma (see \cite{AB, AC} for generalizations of the specialization lemma). 
This interplay between curves and graphs has yielded several applications to the classical algebraic geometry such as a tropical proof of the famous Brill--Noether theorem~\cite{CDPR} 
(see also~\cite{Ca, LPP}). 

In the specialization lemma, it is often that $r_{\Gamma} (\tau_*(\widetilde{D}))$ is larger than $r_{X}(\widetilde{D})$ (see e.g. Example~\ref{eg:2}). 
In this paper, we study when the ranks of divisors are preserved under the specialization map (see Proposition~\ref{prop:relation} for our original motivation).  
By a finite graph, we mean an unweighted, finite connected multigraph, where loops are allowed. A vertex-weighted graph $(G, \omega)$ is the pair of 
a finite graph $G$ and a function $\omega: V(G) \to \ZZ_{\geq 0}$, where 
$V(G)$ denotes the set of vertices of $G$. 

\begin{Question}
\label{q}
Let $(G, \omega)$ be a vertex-weighted graph, and 
let $\Gamma$ be the metric graph associated to $G$. 
Under what condition on $(G, \omega)$, does there exist a regular, generically smooth, semi-stable $R$-curve $\Xscr$ with reduction graph $(G, \omega)$ satisfying the following condition? 
\begin{enumerate}
\item[(C)]
Let $X$ be the generic fiber of $\Xscr$, and $\tau: X({\overline\KK}) \to \Gamma_\QQ$ the specialization map. 
Then, for any $D \in \Div(\Gamma_\QQ)$, 
there exists a divisor $\widetilde{D} \in \Div(X_{\overline\KK})$ 
such that $D = \tau_*(\widetilde{D})$ and 
$r_{(\Gamma, \omega)}(D) = r_{X}(\widetilde{D})$. 
\end{enumerate}
\end{Question}

The purpose of this paper is to answer Question~\ref{q} for {\em hyperelliptic} graphs. Here, a vertex-weighted graph $(G, \omega)$  is hyperelliptic if the genus of $(G, \omega)$ is at least $2$ and there exists a 
divisor $D$ on $\Gamma$ such that $\deg(D) = 2$ and $r_{(\Gamma, \omega)}(D) =1$ (see Definition~\ref{def:hyp:vw}). 
An edge $e$ of $G$ is called a {\em bridge} if the deletion of $e$ makes $G$ 
disconnected. Let $G_1$ and $G_2$ denote the connected components of $G\setminus\{e\}$, which are respectively equipped with the vertex-weight functions $\omega_1$ and $\omega_2$ given by the restriction of $\omega$. 
A bridge is called a {\em positive-type} bridge if 
each of $(G_1, \omega_1)$ and $(G_2, \omega_2)$ has genus at least $1$. 

With the notation in Question~\ref{q}, we also consider the following 
condition \textup{(C')}, which implies (C)  (see Lemma~\ref{lemma:prime}). 
\begin{enumerate}
\item[(C')]
For any $D \in \Div(\Gamma_\QQ)$, 
there exist a divisor $E = \sum_{i=1}^k n_{i} [v_i] 
\in \Div(\Gamma_\QQ)$ that is linearly equivalent to $D$ 
and a divisor $\widetilde{E} = \sum_{i=1}^k n_{i} P_i \in \Div(X_{\overline\KK})$
such that $\tau(P_i) = v_i$ for any $1 \leq i \leq k$ and 
$r_{(\Gamma, \omega)}(E) = r_{X}(\widetilde{E})$. 
\end{enumerate}

Our main result is as follows. 

\begin{Theorem}
\label{thm:main}
Let $\KK$ be a complete discrete valuation field with ring of integers $R$ 
and algebraically closed residue field $k$. Assume that ${\rm char}(k) \neq 2$. 
Let $(G, \omega)$ be a hyperelliptic vertex-weighted graph. Then the following are equivalent. 
\begin{enumerate}
\item[(i)]
For every vertex $v$ of $G$, there are at most $(2\, \omega(v) + 2)$ positive-type bridges emanating from~$v$.  
\item[(ii)]
There exists a regular, generically smooth, semi-stable $R$-curve $\Xscr$ with reduction graph $(G, \omega)$ 
which satisfies the condition \textup{(C)}.  
\item[(iii)]
There exists a regular, generically smooth, semi-stable $R$-curve $\Xscr$ with reduction graph $(G, \omega)$ 
which satisfies the condition \textup{(C')}.
\end{enumerate}
\end{Theorem}

In fact, we will see that the condition (i) is equivalent to the existence of a regular, generically smooth, semi-stable $R$-curve $\Xscr$ with reduction graph $(G, \omega)$ such that 
$\Xscr_\KK$ is {\em hyperelliptic} (see Theorem~\ref{thm:lifting}), and that any such $R$-curve $\Xscr$ satisfies the conditions (C) and (C'). 

As a corollary, we have the following vertex-weightless version. 
A semi-stable $R$-curve $\Xscr$ is said to be 
{\em strongly semi-stable} if every component of the special fiber is smooth, and 
{\em totally degenerate} if every component of the special fiber is 
a rational curve.
Let $(G , \omega)$ be the vertex-weighted reduction graph of 
an $R$-curve $\Xscr$.
Note that,
if $\Xscr$ is strongly semi-stable,
then $G$ is loopless,
and if $\Xscr$ is totally degenerate, then $\omega = \mathbf{0}$.

\begin{Corollary}
\label{cor:main}
Let $\KK, R$ and $k$ be as in Theorem~\ref{thm:main}. Let $G = (G, \mathbf{0})$ be a loopless hyperelliptic graph. Then the following are equivalent. 
\begin{enumerate}
\item[(i)]
For every vertex of $G$, there are at most $2$ positive-type bridges emanating from it. 
\item[(ii)]
There exists a regular, generically smooth, strongly semi-stable, totally degenerate $R$-curve $\Xscr$ with reduction graph $G$ which satisfies the condition \textup{(C)} \textup{(}with $r_\Gamma$ in place of $r_{(\Gamma, \omega)}$\textup{)}.  
\item[(iii)]
There exists a regular, generically smooth, strongly semi-stable, totally degenerate $R$-curve $\Xscr$ with reduction graph $G$ which satisfies the condition \textup{(C')} \textup{(}with $r_\Gamma$ in place of $r_{(\Gamma, \omega)}$\textup{)}.  
\end{enumerate}
\end{Corollary}

We have come to consider Question~\ref{q}  in our desire to understand relationship between the Riemann--Roch formula on graphs and that on curves. Indeed, we have the following Proposition~\ref{prop:relation}. (Since the Riemann--Roch formula on vertex-weighted graphs  
is a corollary of that on vertex-weightless graphs, we give the  vertex-weightless version.) Recall that the Riemann--Roch formula on a metric graph asserts that 
\begin{equation}
\label{eqn:1:1}
 r_\Gamma(D) - r_\Gamma(K_\Gamma -D) = \deg(D) + 1 - g(\Gamma) 
\end{equation}
for any $D \in \Div(\Gamma)$ (cf. \cite{BN, GK, MZ}), where 
the canonical divisor of a compact connected metric graph 
$\Gamma$ is defined to be 
$
  K_\Gamma := \sum_{v \in \Gamma} (\val(v)-2) [v]
$
(cf. \cite{Zh}). 

\begin{Proposition}
\label{prop:relation}
Let $G$ be a finite graph 
and $\Gamma$ the metric graph associated to $G$. 
Assume that there exist a complete discrete valuation field $\KK$ 
with ring of integers $R$, and a regular, generically smooth, strongly semi-stable, totally degenerate  $R$-curve $\Xscr$ with reduction graph $G$ which satisfies 
the condition~\textup{(C)}. Then the Riemann--Roch formula on 
$\Gamma$ is deduced from the Riemann--Roch formula on $X_{\overline\KK}$, where $X$ is the generic fiber of $ \Xscr$. 
\end{Proposition}

Let $G$ be a loopless hyperelliptic graph. Let $\overline{G}$ be the 
hyperelliptic graph that is obtained by contracting all the bridges of $G$.  
Then Corollary~\ref{cor:main}, Proposition~\ref{prop:relation} and comparison of divisors on $G$ and $\overline{G}$ gives a proof of the Riemann--Roch formula on a loopless hyperelliptic graph $G$ (see Remark~\ref{rmk:relation}).  
It should be noted, however, that, as the original proof by Baker--Norine, 
this proof uses the theory of reduced divisors (in the proof of 
Theorem~\ref{thm:main}). 

Let $(G, \omega)$ be a vertex-weighted graph, and $\Gamma$ the metric graph associated to $G$. Question~\ref{q} is also of interest from 
the viewpoint of the Brill-Noether theory: 
For fixed integers $d, r \geq 0$, we put $W_{d}^r(\Gamma_\QQ, \omega) := \{D \in \Div(\Gamma_\QQ) \mid \deg(D) = d, r_{(\Gamma, \omega)}(D) \ \geq r \}$; If the condition (C) is satisfied with an $R$-curve $\Xscr$ with generic fiber $X$, 
then we will have $\tau_*(W_{d}^r(X_{\overline{\KK}})) = W_{d}^r(\Gamma_\QQ, \omega)$. 

\bigskip
Caporaso has kindly informed us that the condition (C) is related to 
her conjecture \cite[Conjecture~1]{Ca2}. Let $(G, \omega)$ be a vertex-weighted graph, and let $D \in \Div(G)$. 
The {\em algebro-geometric rank} $r_{(G, \omega)}^{\alg, k}(D)$ of $D$ is 
defined by  
\begin{align*}
  r_{(G, \omega)}^{\alg, k}(D)  & := \max_{X_0}\, r(X_0, D), \\
  r(X_0, D)  & := \min_{E}\, r^{\max}(X_0, E), \\
  r^{\max}(X_0, E)  & := \max_{\Escr_0}\, \left( h^0(X_0, \Escr_0) -1 \right), 
\end{align*}
where $X_0$ runs over all connected reduced projective nodal curves defined over $k$ with dual graph $(G, \omega)$,  
$E$ runs over all divisors on $G$ that are linearly equivalent to $D$ in $\Div(G)$,  and $\Escr_0$ runs over all Cartier divisors on $X_0$ such that $\deg\left(\Escr_0\vert_{C_v}\right) = E(v)$ for any $v \in  V(G)$. (Here $C_v$ denotes the irreducible component of $X_0$ corresponding to $v$.)  
In \cite[Conjecture~1]{Ca2}, Caporaso has conjectured that 
\begin{equation}
\label{eqn:Caporaso:conj}
  r_{(G, \omega)}^{\alg, k}(D) = r_{(\Gamma, \omega)}(D)
\end{equation}
and showed that \eqref{eqn:Caporaso:conj} holds
in the following four cases: (1) $g(\Gamma, \omega) \leq 1$;  (2) $\deg(D) \leq 0$ or $\deg(D) \geq 2 g(\Gamma, \omega) -2$; (3) $G$ has exactly one vertex; and (4) $g(\Gamma, \omega) \leq 2$ and $(G, \omega)$ is stable. 
Caporaso has informed us about her very recent and unpublished  work with M. Melo proving one direction of the conjecture, i.e., 
$r_{(G, \omega)}^{\alg, k}(D) \leq r_{(\Gamma, \omega)}(D)$. 

To make the relation between (C) and \eqref{eqn:Caporaso:conj} precise, we consider a variant of the condition (C), which is concerned with the existence of a lifting as a divisor over $\KK$ (not just as a divisor over $\overline{\KK}$) of a divisor $D$ on $G$ (not just on $\Gamma_\QQ$). Let the notation be as in Question~\ref{q}. Let $\rho_*: \Div(X) \to \Div(G)$ be the specialization map (see \eqref{eqn:specialization:map:fin}). 
\begin{enumerate}
\item[(F)]
For any $D \in \Div(G)$, 
there exists a divisor $\widetilde{D} \in \Div(X)$ 
such that $D = \rho_*(\widetilde{D})$ and 
$r_{(\Gamma, \omega)}(D) = r_{X}(\widetilde{D})$. 
\end{enumerate}

The following proposition, which is due to Caporaso, 
shows that the condition (F) leads to 
the other direction in her conjecture. 

\begin{Proposition}
\label{prop:Caporaso}
Let $\KK, R$ and $k$ be as in Theorem~\ref{thm:main}.  Let $(G, \omega)$ be a vertex-weighted graph, and let $\Gamma$ be the metric graph associated to $G$. Let $\Xscr$ be a regular, generically smooth, semi-stable $R$-curve with generic fiber $X$ and reduction graph $(G, \omega)$. Assume that $\Xscr$ satisfies the condition \textup{(F)}.  Then, 
for any divisor $D \in \Div(G)$, 
we have 
\[
  r_{(G, \omega)}^{\alg, k}(D) \geq r_{(\Gamma, \omega)}(D).  
\]
\end{Proposition}

For a hyperelliptic vertex-weighted graph $(G, \omega)$, we can show the following (see Theorem~\ref{thm:main:F} for a stronger result, which considers a variant of the condition (C')). 

\begin{Theorem}
\label{thm:finite}
Let $\KK, R$ and $k$ be as in Theorem~\ref{thm:main}. Let $(G, \omega)$ be a hyperelliptic graph such that for every vertex $v$ of $G$, there are at most $(2 \omega(v) +2)$ positive-type bridges emanating from $v$. Then, 
there exists a regular, generically smooth, semi-stable $R$-curve $\Xscr$ with reduction graph $(G, \omega)$ 
which satisfies the condition \textup{(F)}.  
\end{Theorem}

Thus we obtain the following corollary. 

\begin{Corollary}
\label{cor:finite}
Let $k$ be an algebraically closed field with ${\rm char}(k) \neq 2$. Let $(G, \omega)$ be a hyperelliptic graph such that for every vertex $v$ of $G$, there are at most $(2 \omega(v) +2)$ positive-type bridges emanating from $v$. Then, 
for any $D \in \Div(G)$, we have $r_{(G, \omega)}^{\alg, k}(D) \geq r_{(\Gamma, \omega)}(D)$.  
\end{Corollary}

\subsection{Remarks}
A number of remarks are in order. 

\begin{Remark}
In this paper, we consider vertex-weighted graphs (i.e., not only vertex-weightless finite graphs), for vertex-weighted graphs appear naturally in tropical geometry and Berkovich spaces. (Indeed, a vertex-weighted metric graph is seen as a Berkovich skeleton of an algebraic variety over $\KK$. For the interplay between Berkovich spaces and tropical varieties over $\KK$, 
see, for example,~\cite{ABBR, BPR, Gu, Pe}.) 
\end{Remark}

\begin{Remark}
Theorem~\ref{thm:main} treats vertex-weighted hyperelliptic 
graphs of genus at least $2$. 
We also show that, for any vertex-weighted graph of  genus $0$ or $1$, 
there exists a regular, generically smooth, semi-stable $R$-curve $\Xscr$ with reduction graph $(G, \omega)$ that satisfies the condition \textup{(C)} and \textup{(C')} (see Proposition~\ref{prop:tree}). 
\end{Remark}

\begin{Remark}
\label{rmk:1:5}
The condition (C') is in general {\em not} equivalent to the following condition: 
\begin{enumerate}
\item[(C'')]
For any $D =\sum_{i=1}^k n_{i} [v_i]  \in  \Div(\Gamma_\QQ)$, 
there exists $P_i \in X(\overline\KK)$ with $\tau(P_i) = v_i$ for each $1 \leq i \leq k$ such that $r_{(\Gamma, \omega)}(D) = r_{X}(\sum_{i=1}^k n_{i} P_i)$. 
\end{enumerate}
See Example~\ref{eg:3}, where we give a hyperelliptic graph $G$ and a model $\Xscr$ that satisfy the conditions (C) and (C'), but does not satisfy the condition (C"). This example is interesting in two senses. First, for the divisor $D$ in 
Example~\ref{eg:3}, by the condition (C), there exists $\widetilde{D} \in \Div(X_{\overline\KK})$ with  $\tau_*(\widetilde{D}) = D$ and $r_{(\Gamma, \omega)}(D) = r_{X}(\widetilde{D})$. This example shows, however, that  $\widetilde{D}$ is not simply of the form 
$\sum_{i=1}^k n_{i} P_i$ with $\tau(P_i) = v_i$. Secondly, by the condition (C'), if we replace $D$ by a divisor $E =  \sum_{j=1}^\ell m_{j} [w_j]$ with $E \sim D$, then we can indeed lift $E$ in $X$ as a simple form $\widetilde{E} = \sum_{j=1}^\ell m_{j} Q_j$ with $\tau(Q_j) = w_j$ preserving the ranks $r_{(\Gamma, \omega)}(E) = r_{X}(\widetilde{E})$. 
\end{Remark}

\begin{Remark}
\label{rmk:ABBR}
In a very recent paper \cite{ABBR}, Amini, Baker, Brugall\'e and Rabinoff studied lifting of harmonic morphisms of metrized complexes, among others, to morphisms of algebraic curves (see also Theorem~\ref{thm:lifting} below).
In \cite[\S10.11]{ABBR}, they discussed lifting divisors of given rank, giving several examples for which various specialization lemmas do not attain the equality.  Question~\ref{q} will be interesting from this perspective, and 
Theorem~\ref{thm:main} gives a clean picture for the case of hyperelliptic graphs. We also remark that Cools, Draisma, Payne and Robeva considered a certain graph $G_\circ$ of $g$ loops to give a tropical proof of the Brill--Noether theorem and that their conjecture \cite[Conjecture~1.5]{CDPR} concerns lifting of divisors that preserves the ranks between $G_\circ$ and a regular, generically smooth, strongly semi-stable, totally degenerate $R$-curve with reduction graph $G_\circ$. 
\end{Remark}

\subsection{Strategy of the proof and other results}
\label{subsec:intro:3}
We now explain our strategy to prove Theorem~\ref{thm:main}. 
Our starting point is the following theorem. 

\begin{Theorem}[cf. {\cite[Theorem~4.8]{CaGonal}} and {\cite[Theorem~1.10]{ABBR}}]
\label{thm:lifting}
Let $\KK, R$ and $k$ be as in Theorem~\ref{thm:main}, and let $(G, \omega)$ be a vertex-weighted  hyperelliptic graph. Then the condition \textup{(i)} in Theorem~\ref{thm:main} is equivalent to the existence of  a regular, generically smooth, semi-stable $R$-curve $\Xscr$ with reduction graph $(G, \omega)$ such 
that the generic fiber $\Xscr_\KK$ is hyperelliptic. 
\end{Theorem}

Caporaso \cite[Theorem~4.8]{CaGonal} proved that 
the condition \textup{(i)}  in Theorem~\ref{thm:main} is equivalent to the existence of a hyperelliptic semi-stable curve $X_0$ over $k$. Based on \cite[Theorem~4.8]{CaGonal}, we will give a proof of  Theorem~\ref{thm:lifting} using equivariant deformation. 
We remark that there is another approach to Theorem~\ref{thm:lifting}. Amini, Baker, Brugall\'e and Rabinoff \cite[Theorem~1.10]{ABBR} recently showed a skeleton-theoretic version of Theorem~\ref{thm:lifting} as a corollary of their deep studies of  canonical gluing and star analytic spaces over an algebraically closed field with a non-Archimedean valuation (during the preparation of this paper).  With an argument of ``descent'' to the case of a discrete valuation field, 
it may be possible that one derives Theorem~\ref{thm:lifting} from \cite[Theorem~1.10]{ABBR}. 

Theorem~\ref{thm:lifting} shows that (ii)  implies (i) in Theorem~\ref{thm:main}. Since (C') implies (C) (see Lemma~\ref{lemma:prime}), 
the condition (iii) implies (ii) in Theorem~\ref{thm:main}.  
The main part of the proof of Theorem~\ref{thm:main} is 
to show that (i) implies (iii). 

For a metric graph $\Gamma$ 
and $v_0\in \Gamma$, 
a divisor $D \in  \Div(\Gamma)$ is said to be {\em $v_0$-reduced} if 
$D$ is effective away from $v_0$ and satisfies several nice properties (see Definition~\ref{def:reduced}). This notion was introduced by Baker and Norine \cite{BN}, and is a powerful tool in computing the ranks of divisors. 
With the notion of moderators (see \cite[Theorem~3.3]{BN}, \cite[Section~7]{MZ}, \cite[Corollary~2.3]{HKN}), we have the following 
properties of reduced divisors. 

\begin{Theorem}
\label{thm:main:3}
Let $\Gamma$ be a compact connected metric graph of genus $g \geq 2$. 
We fix a point $v_0 \in \Gamma$.  
Let $D \in \Div(\Gamma)$ be a $v_0$-reduced divisor on $\Gamma$, 
and let  $D(v_0)$ denote the coefficient of $D$ at $v_0$. Then, if $\deg(D) -  D(v_0) \leq g-1$, 
then there exists $w\in \Gamma \setminus\{v_0\}$ 
such that $D + [w]$ is a $v_0$-reduced divisor. 
\end{Theorem}

Let $\Gamma$ be a hyperelliptic metric graph. 
We fix $v_0 \in \Gamma$ satisfying \eqref{eqn:v0}. 
We set,  for an effective divisor $D \in \Div(\Gamma)$,  
$
  p_{\Gamma}(D) = \max\{
  r \in \ZZ_{\geq 0} \mid |D - 2 r [v_0] | \neq \emptyset
  \}$. 
We similarly define   $p_{(\Gamma, \omega)}(D)$ on a hyperelliptic 
vertex-weighted graph $(\Gamma, \omega)$ (See Sect.~\ref{subsec:p:Gamma:D}). 
Using Theorem~\ref{thm:main:3}, we 
compute $r_{(\Gamma, \omega)}(D)$ in terms of 
$p_{(\Gamma, \omega)}(D)$, which is a key ingredient  
of the proof of  Theorem~\ref{thm:main}. 

\begin{Theorem} 
\label{thm:main:2}
Let $(G, \omega)$ be a hyperelliptic vertex-weighted graph 
of genus $g$, 
and $\Gamma$ the metric graph associated to $G$. 
Then, for any effective divisor $D$ on $\Gamma$, we have 
\[
 r_{(\Gamma, \omega)}(D) = \begin{cases}
  p_{(\Gamma, \omega)}(D) & \textup{(if $\deg(D) - p_{(\Gamma, \omega)}(D) \leq g$), } \\
  \deg(D) -g & \textup{(if $\deg(D) - p_{(\Gamma, \omega)}(D) \geq g+1$).} 
  \end{cases} 
\]  
\end{Theorem}

There is a corresponding formula in the classical setting 
of ranks of divisors on hyperelliptic curves (see Proposition~\ref{prop:Yamaki}). 
We deduce  (iii)  from (i) in Theorem~\ref{thm:main}, combining Theorem~\ref{thm:lifting}, 
Theorem~\ref{thm:main:2} and Proposition~\ref{prop:Yamaki}. 

The organization of this paper is as follows. 
In Sect.~\ref{sec:prelim}, we briefly recall the theory of divisors 
on metric graphs. In Sect.~\ref{sec:hyp:graph}, 
we consider hyperelliptic graphs.
In Sect.~\ref{sec:hyp:semistable:curves}, 
we consider hyperelliptic semi-stable curves 
and prove Theorem~\ref{thm:lifting} using equivariant deformation. 
In Sect.~\ref{sec:reduced}, we prove Theorem~\ref{thm:main:3}. 
In Sect.~\ref{sec:rank:graph}, we study ranks of divisors on 
a hyperelliptic graph, and prove Theorem~\ref{thm:main:2}.  
In Sect.~\ref{sec:proof}, we prove Theorem~\ref{thm:main} and 
Proposition~\ref{prop:relation}. We also consider Question~\ref{q} for vertex-weighted graphs of genus $0$ or $1$. In Sect.~\ref{sec:Caporaso}, 
we consider variants of the condition (C) and (C'), and show Proposition~\ref{prop:Caporaso}, Theorem~\ref{thm:finite} and Corollary~\ref{cor:finite}.    
In the appendix, we put together some results on the 
deformation theory which are needed in Sect.~\ref{sec:hyp:semistable:curves}. 

\smallskip
{\sl Acknowledgments.}\quad 
The authors express their deep gratitude to Professor Lucia Caporaso for invaluable comments on the previous version of this paper, sharing her insight regarding relationship between  Question~\ref{q} and her conjecture, and letting them know the paper \cite{CaGonal}. The authors express their deep gratitude to Professors Omid Amini and Matthew Baker for invaluable comments on 
the previous version of this paper. The authors express their deep gratitude 
to the referees for carefully reading the paper, giving many invaluable comments and simplifying the proofs of Theorem~\ref{thm:main:3}, Theorem~\ref{thm:main:2} and Proposition~\ref{prop:Yamaki}. 

\setcounter{equation}{0}
\section{Preliminaries}
\label{sec:prelim}
In this section, we briefly recall the theory of divisors on a compact 
metric graph, Baker's Specialization Lemma, and the notion of reduced divisors on a metric graph, which we use later. 
We also recall some properties of a vertex-weighted graph and a contraction of metric graphs. 

\subsection{Theory of divisors on a metric graph}  
\label{subsec:RR}
We briefly recall the theory of divisors on metric graphs. 
We refer the reader to  \cite{BN, GK, HKN, MZ} 
for details and further references. 

Throughout this paper, a {\em finite graph} means an unweighted, finite connected multigraph. Notice that we allow the existence of loops. For a finite graph $G$, let $V(G)$ denotes the set of vertices, and $E(G)$ the set of edges. The genus of $G$ is defined to be $g(G) = |E(G)| - |V(G)| + 1$. 
For $v \in V(G)$, the {\em valence} $\val(v)$ of $V$ is the number of edges emanating from $v$. 
Recall from the introduction that $e \in E(G)$ 
is called a {\em bridge} if the deletion of $e$ makes $G$ disconnected. 
A vertex $v$ of $G$ is a {\em leaf end} if $\val(v) = 1$.  A {\em leaf edge} is an edge of $G$ that has a leaf end. In particular, a leaf edge is a bridge. 

An {\em edge-weighted graph} $(G, \ell)$ is the pair of a finite graph $G$ and 
a function (called a length function) $\ell: E(G) \to \RR_{>0}$.  In other words, 
an edge-weighted graph means a finite graph having each edge assigned a positive length.  
A {\em compact connected metric graph} $\Gamma$
is the underlining metric space of an edge-weighted graph $(G, \ell)$. 
We say that $(G, \ell)$ is a model of $\Gamma$. There are many possible models 
for $\Gamma$. However, if $\Gamma$ is not a circle, we can canonically construct 
a model $(G_\circ, \ell)$ of $\Gamma$ as follows (cf. \cite{Chan}). 
The set of vertices is given by 
$
  V(G_\circ) := \{v \in \Gamma \mid \val(v) \neq 2\},  
$
where the valence $\val(v)$ is the number of connected components of $U_v \setminus\{v\}$ with $U_v$ being any small neighborhood of 
$v$ in $\Gamma$. The set of edges $E(G_\circ)$ corresponds to the set of connected components  
of $\Gamma\setminus V(G_\circ)$. Since each connected component of 
$\Gamma\setminus V(G_\circ)$ is an open interval, its length determines 
the length function $\ell$. The model $(G_\circ, \ell)$ is called 
the {\em canonical model} of $\Gamma$. 

Let $\Gamma$ be a compact connected metric graph. 
By a cut-vertex of $\Gamma$, we mean a point $v$ of $\Gamma$ such that 
$\Gamma\setminus\{v\}$ is disconnected. 
By an edge of $\Gamma$, 
we mean an edge of the underlining graph $G_\circ$ of the canonical 
model $(G_\circ, \ell)$. Similarly, by a bridge (reps. 
a leaf edge) of $\Gamma$, we mean a bridge (reps. a leaf edge) of $G_\circ$. Let $e$ be an edge of $\Gamma$ that is not a loop. 
We regard $e$ as a closed subset of $\Gamma$, i.e., 
including the endpoints $v_1, v_2$ of $e$. We set 
$\overset{\scriptsize{\circ}}{e} = e \setminus\{v_1, v_2\}$. 
  
The genus $g(\Gamma)$ of a compact connected metric graph $\Gamma$ 
is defined to be its first Betti number, 
which equals $g(G)$ of any model $(G, \ell)$  of $\Gamma$. 
An element of the free abelian group $\Div(\Gamma)$ generated 
by points of $\Gamma$ is called a {\em divisor} on $\Gamma$. 
For $D = \sum_{v\in \Gamma} n_v [v]\in \Div(\Gamma)$,  
its {\em degree} is defined by 
$
  \deg(D) = \sum_{v\in \Gamma} n_v
$. 
We write the coefficient $n_v$ at $[v]$ for $D(v)$. 
A divisor 
$D = \sum_{v\in  \Gamma} n_v [v]\in \Div( \Gamma)$ is said to be {\em effective} 
if $D(v) \geq 0$ for any $v \in \Gamma$.
If $D$ is effective, we write $D \geq 0$.

A {\em rational function} on $\Gamma$ is a piecewise linear function 
on $\Gamma$
with integer slopes.
We denote by $\Rat(\Gamma)$ the set of rational functions on $\Gamma$. For 
$f \in \Rat(\Gamma)$ and a point $v$ in $\Gamma$, the sum of the outgoing  
slopes of $f$ at $v$ is denoted by $\ord_v(f)$. This sum is $0$ except for all but finitely 
many points of $\Gamma$, and thus 
\[
  \zero(f) := \sum_{v \in \Gamma} \ord_v(f) [v] 
\]
is a divisor on $\Gamma$. 
The set of {\em principal divisors} on $\Gamma$ is defined to be  
$ \Prin(\Gamma) := \{\zero(f) \mid f \in \Rat(\Gamma)\}$. 
Then $\Prin(\Gamma)$ is a subgroup of $\Div(\Gamma)$. 
Two divisors $D, E \in \Div(\Gamma)$ are said to be {\em linearly equivalent}, 
and we write $D \sim E$, if $D- E \in \Prin(\Gamma)$.  
For $D \in \Div(\Gamma)$, the complete linear system $|D|$ is defined by 
\[
  |D| = \{E \in \Div(\Gamma) \mid E \geq 0, \quad E \sim D\}. 
\]

Let $G$ be a finite graph. We say that $\Gamma$ is the 
{\em metric graph associated to} $G$  
if $\Gamma$ is the underlining metric space of $(G, \mathbf{1})$, where $\mathbf{1}$ 
denotes the length 
function which assigns to each edge of $G$ length~$1$. If this is the case, let $\Gamma_\QQ$ 
denote the set of points on $\Gamma$ whose distance from every vertex of $G$ is rational, and let $\Div(\Gamma_\QQ)$ denote the free abelian group  generated by the elements of $\Gamma_\QQ$. 

\begin{Definition}[Rank of a divisor, cf.~\cite{BN}]
Let $\Gamma$ be a compact connected metric graph. 
Let $D \in \Div(\Gamma)$. If $|D| = \emptyset$, then we set 
$r_\Gamma(D) := -1$. If $|D| \neq \emptyset$, we set 
\[
 r_\Gamma(D) := \max\left\{s \in \ZZ \;\left\vert\; 
 \begin{gathered}
 \text{For any effective divisor $E$ with $\deg(E) = s$,} \\
 \text{ we have $|D-E| \neq \emptyset$}
 \end{gathered}
 \right. \right\}. 
\]
\end{Definition}

We compare divisors on a compact connected metric graph 
$\Gamma$ and those on the metric graph obtained by contracting 
a bridge of $\Gamma$. Let $\Gamma$ be a compact connected metric graph. 
Suppose that $\Gamma$ has a bridge $e$, and let $\Gamma_1$ be the graph 
obtained by contracting $e$. 
Let 
$
  \varpi_1: \Gamma \to \Gamma_1
$ be the retraction map. 

\begin{Lemma}[{\cite[Lemma~3.11]{Chan}}]
\label{lemma:linear:equiv}
Let $\Gamma, \Gamma_1$ and $\varpi_1$ be as above. 
Let $D \in \Div(\Gamma)$ and $D_1 \in \Div(\Gamma_1)$. 
\begin{enumerate}
\item
We have $D \in \Prin(\Gamma)$ if and only if $\varpi_{1*}(D) \in \Prin(\Gamma_1)$.  
\item
We have $r_\Gamma(D) = r_{\Gamma_1}(\varpi_{1*}(D))$. 
\item
Suppose that the contracted bridge $e$ is a leaf edge, so that we have 
the natural embedding $\jmath_1: \Gamma_1 \hookrightarrow \Gamma$. 
Then we have $D \sim \jmath_{1*}(\varpi_{1*}(D))$ on $\Div(\Gamma)$. 
\item
Under the assumption of \textup{(3)}, 
we have $r_\Gamma(\jmath_{1*}(D_1)) = r_{\Gamma_1}(D_1)$. 
\end{enumerate}
\end{Lemma}

\Proof
(1) See \cite[Lemma~3.11]{Chan}. 
(2) This follows from (1) by the argument in 
\cite[Corollaries~5.10, 5.11]{BN2}. 
(3) Since $\varpi_{1*} (D - \jmath_{1*}(\varpi_{1*}(D))) = 0$, 
the assertion follows from (1). 
(4) Since $\varpi_{1*}(\jmath_{1*}(D_1)) 
= D_1$, the assertion follows from (2).  
\QED

\subsection{Reduced divisors on a metric graph}
\label{subsec:reduced}
We briefly recall the notion of {\em reduced divisors} 
on a graph, which is a powerful tool in computing the ranks of divisors. 
Reduced divisors were introduced in \cite{BN} 
to prove the Riemann--Roch formula on a finite graph.  

Let $\Gamma$ be a compact connected metric graph. For any closed subset $A$ of $\Gamma$ and $v \in \Gamma$, the {\em out-degree} of $v$ from $A$, denoted by $\mathrm{outdeg}^\Gamma_A(v)$, is defined to be the maximum number of internally disjoint segments of $\Gamma\setminus A$ with an open end $v$. Note that if $v \in A \setminus \partial A$, then $\mathrm{outdeg}^\Gamma_A(v) = 0$. For $D \in \Div(\Gamma)$, a point $v \in \partial A$ is  {\em saturated} 
for $D$ with respect to $A$ if $D(v) \geq \mathrm{outdeg}^\Gamma_A(v)$, 
and {\em non-saturated} otherwise. 

\begin{Definition}[$v_0$-reduced divisor]
\label{def:reduced}
We fix a point $v_0\in\Gamma$.  
A divisor $D \in \Div(\Gamma)$ is called a {\em $v_0$-reduced divisor} if 
$D$ is non-negative on $\Gamma\setminus\{v_0\}$, and every compact  
subset $A$ of $\Gamma\setminus\{v_0\}$ contains a non-saturated point 
$v \in \partial A$ for $D$ with respect to $A$. 
\end{Definition}
 
We remark that we may require that a compact subset $A$ of 
$\Gamma\setminus\{v_0\}$ be connected in the above definition. 

We put together useful
properties of $v_0$-reduced divisors 
in the following theorem. 

\begin{Theorem}[\cite{B, BN, HKN}]
\label{thm:HKN}
Let $D \in \Div(\Gamma)$ and $v_0 \in \Gamma$. 
\begin{enumerate}
\item
There exists a unique $v_0$-reduced divisor 
$D_{v_0}$ that is linearly equivalent to 
$D$. 
\item
The divisor $D$ is linearly equivalent to an effective divisor if and only if $D_{v_0}$ 
is effective. 
\item
Suppose that $\Gamma$ is the metric graph associated to a finite graph $G$ and that $v_0 \in \Gamma_\QQ$. Then, 
if $D \in \Div(\Gamma_\QQ)$, then $D_{v_0} \in \Div(\Gamma_\QQ)$. 
\end{enumerate}
\end{Theorem} 

For a given divisor $D \in \Div(\Gamma)$, Luo \cite{Luo} gives a criterion that 
$D$ is a $v_0$-reduced divisor based on Dhar's algorithm. 
Here we give a slightly modified version of \cite[Algorithm~2.5]{Luo}. 

\begin{Theorem}[cf. \cite{Luo}]
\label{thm:DL}
Let $v_0 \in \Gamma$. 
Let $D$ be an effective divisor on $\Gamma$ 
such that $D(v_0) = 0$.  
Then $D$ is $v_0$-reduced if and only if there exists a sequence 
\begin{equation}
\label{eqn:S}
  \mathbf{a} = (a_1, a_2, \ldots, a_k)
\end{equation}
with the following properties\textup{:} 
\begin{enumerate}
\item[(i)]
The points $a_1, a_2, \ldots, a_k$ are mutually distinct points of $\Gamma \setminus \{v_0\}$. 
\item[(ii)]
We have 
$\Supp(D) \subseteq \{a_1, a_2, \ldots, a_k\}$. 
\item[(iii)]
For $1 \leq i \leq k$, let $\Ucal_i$ be the connected component of 
$\Gamma\setminus \{a_i, a_{i+1}, \ldots, a_k\}$ that contains $v_0$, and put 
$A_i := \Gamma\setminus \Ucal_i$. Then 
$a_i \in \partial A_i$ and $a_i$ is a non-saturated point for $D$ with respect to 
$A_i$. 
\end{enumerate}
\end{Theorem}

\Proof
Because Theorem~\ref{thm:DL} is slightly different from \cite{Luo}, 
we give a brief proof. Suppose that $D$ is $v_0$-reduced. 
We construct a sequence $\mathbf{a}$ inductively. If 
$(a_1, \ldots, a_{i-1})$ is chosen, we put $S_{i-1} := \Supp(D) \setminus 
\{a_1, \ldots, a_{i-1}\}$. (For the first stage, we let $S_0 := \Supp(D)$.)  
Let $\Vcal$ be the connected component of 
$\Gamma\setminus S_{i-1}$ which contains $v_0$, and put 
$B := \Gamma\setminus \Vcal$. 
Since $D$ is $v_0$-reduced, there exists a non-saturated point 
$b \in \partial B$ for $D$ with respect to $B$. Then $b \in S_{i-1}$. We 
define $a_i := b$. Then $\mathbf{a} = (a_1, a_2, \ldots, a_k)$ satisfies (i)(ii) and (iii). 
We remark that in this construction 
we have $\Supp(D) = \{a_1, a_2, \ldots, a_k\}$, which is stronger than~(ii).  

On the other hand, suppose that there exists a sequence $\mathbf{a}$ 
satisfying (i), (ii) and (iii). Let $S$ be any subset of $\Supp(D)$. 
Let $\Ucal$ be the connected component of 
$\Gamma\setminus S$ which contains $v_0$, and put 
$A:= \Gamma\setminus \Ucal$. We take an $i$ with 
$S \subseteq \{a_{i}, a_{i+1}, \ldots, a_k\}$ and 
$S \not\subseteq \{a_{i+1}, a_{i+2}, \ldots, a_k\}$.  
Then $a_i \in \partial A$, and we have  
$D(a_i) < \outdeg^\Gamma_{A_i}(a_i) \leq \outdeg^\Gamma_{A}(a_i)$. Thus $a_i$ 
is a non-saturated point. Then \cite[Lemma~2.4]{Luo} tells us  that  
$D$ is $v_0$-reduced. 
\QED

\subsection{Specialization lemma}
\label{subsec:sp}
In this subsection, following \cite{B}, 
we briefly recall the relationship between linear systems 
on curves and those on graphs, and Baker's Specialization Lemma. 

Let  $\KK$ be a complete discrete valuation field with ring of integers $R$ 
and algebraically closed residue field $k$. 
Let  $X$ be a geometrically irreducible smooth projective curve over $\KK$. 
We assume that $X$ has a {\em semi-stable} model over $R$, i.e., 
there exists a regular $R$-curve $\Xscr$ whose generic fiber is isomorphic to $X$ 
and whose special fiber $\Xscr_0$ is a reduced scheme with at most nodes (i.e., ordinary double points) as singularities. 

The {\em dual graph} $G$ associated to 
$\Xscr_0$ is defined as follows.  
Let $X_1, \ldots, X_r$ be the irreducible components of $\Xscr_0$. Then 
$G$ has vertices $v_1, \ldots, v_r$ which correspond to $X_1, \ldots, X_r$, respectively. Two vertices $v_i, v_j$ ($i \neq j$) of $G$ are connected by 
$a_{ij}$ edges if $\# X_i \cap X_j = a_{ij}$.  A vertex $v_i$ has $b_i$ loops if 
$\#{\rm Sing}(X_i) = b_i$.  
We 
call the dual graph of $\Xscr_0$ the
{\em reduction graph} of the $R$-curve $\Xscr$.

Let $\Gamma$ be the metric graph associated to $G$, where each edge of $G$ is 
assigned length $1$. Let $P \in X(\KK)$. By the valuative criterion of properness, 
$P$ gives the section $\Delta_P$ over $R$, which meets an irreducible component of the 
special fiber in the smooth locus. Let $v \in G$ be the vertex corresponding to this component. 
We denote by $\tau: X(\KK) \to \Gamma$ the map which assigns $P$ to $v$. 
Suppose that $\KK^\prime$ is a finite extension field of $\KK$ with ring of integers $R^\prime$.  Let $e(\KK^\prime/\KK)$ denote the ramification index of 
$\KK^\prime/\KK$. Let $\Xscr^\prime$ be the minimal resolution of 
$\Xscr \times_{\Spec(R)} \Spec(R^\prime)$. 
Then the generic fiber of $\Xscr^\prime$  
is $X \times_{\Spec(\KK)} \Spec(\KK^\prime)$. Let $G^\prime$ be the dual graph of 
the special fiber of $\Xscr^\prime$. Let $\Gamma^\prime$ be a metric graph  
whose underlining graph is $G^\prime$, where each edge of $G^\prime$ is assigned length 
$1/e(\KK^\prime/\KK)$.  
Then $\Gamma^\prime$ is naturally isometric to $\Gamma$. 
We can extend $\tau$ to a map 
(again denoted by $\tau$ by slight abuse of notation)
\begin{equation}
\label{eqn:specialization:map:0}
  \tau : X ( \overline{\KK} ) \to \Gamma,
\end{equation}
which is called the {\em specialization map} (cf. \cite{CR}). Let
\begin{equation}
\label{eqn:specialization:map}
  \tau_* : \Div ( X_{\overline{\KK}})
  \to \Div ( \Gamma )
\end{equation}
be the induced group homomorphism.

\begin{Proposition}[\cite{B}]
\label{prop:Baker:properties}
\begin{parts}
\Part{(1)}
One has $\Image (\tau) =  \Gamma_\QQ$ and 
$\Image (\tau_*) = \Div(\Gamma_\QQ)$. 
\Part{(2)}
The map $\tau_*$ respects the linear equivalence. 
\Part{(3)}
For any $\widetilde{D} \in\Div(X_{\overline{\KK}})$, 
$\deg \tau_*(\widetilde{D}) = \deg \widetilde{D}$. 
\end{parts} 
\end{Proposition}

\Proof
For (1), see \cite[Remark~2.3]{B}. 
For (2), we refer to \cite[Lemma~2.1]{B}. 
The statement (3) is obvious from the definition of $\tau$. 
We note that, in \cite{B}, each component of the special fiber $\Xscr_0$ is assumed to be smooth, but the arguments in \cite{B} also hold when a component of 
$\Xscr_0$ has a node. 
\QED

We state Baker's Specialization Lemma  \cite{B}. Again, the arguments in \cite{B} hold when a component of $\Xscr_0$ has a node. (This is because the rank of a divisor is measured 
by $r_\Gamma$, not by $r_G$.)

\begin{Theorem}[Baker's Specialization Lemma \cite{B}] 
\label{thm:specialization:lemma}
For any $\widetilde{D} \in\Div(X_{\overline{\KK}})$, one has 
$
 r_\Gamma(\tau_*(\widetilde{D})) \geq r_X(\widetilde{D})
$. 
\end{Theorem}

\subsection{Vertex-weighted graph}
\label{subsec:vertex:weighted:graph}
In this subsection, following \cite{AC}, 
we briefly recall some properties of vertex-weighted graphs.  

A {\em vertex-weighted} graph $(G, \omega)$ is the pair of a finite graph $G$ and a function (called a vertex-weight function) $\omega: V(G) \to \ZZ_{\geq 0}$. The genus of $(G, \omega)$ is defined to be  
$g(G, \omega) =  g(G) + \sum_{v \in V(G)} \omega(v)$. 
For each vertex $v \in V(G)$, we add $\omega(v)$ loops to $G$ at the vertex $v$ to make a new finite graph $G^{\omega}$. The graph $G^{\omega}$ is called 
the {\em virtual weightless finite graph} associated to a vertex-weighted graph $(G, \omega)$. The attached loops are called {\em virtual loops}. 

Let $(G, \omega)$ be a vertex-weighted graph, and 
$e$ a bridge of $G$. Let $G_1, G_2$ denote the connected components of $G\setminus\{e\}$, which are equipped with the vertex-weight functions $\omega_1, \omega_2$ given by the restriction of $\omega$. 
We say that $e$ is a {\em positive-type} bridge if 
each of $(G_1, \omega_1)$ and $(G_2, \omega_2)$ has genus at least $1$. 

A {vertex-weighted metric} graph $(\Gamma, \omega)$ is the pair of a compact connected metric graph $\Gamma$ and a function 
$\omega: \Gamma \to \ZZ_{\geq 0}$ such that $\omega(v) = 0$ except for all but finitely many points $v$ in $\Gamma$. The genus of $(\Gamma, \omega)$ is defined to be $g(\Gamma, \omega) 
= g(\Gamma) + \sum_{v \in \Gamma} \omega(v)$. 
For each point $v \in \Gamma$ with $\omega(v) >0$, we add $\omega(v)$ length-one-loops to the point $v$ to make a new metric graph $\Gamma^{\omega}$. 
We call $\Gamma^{\omega}$ the {\em virtual weightless metric graph} associated to $(\Gamma, \omega)$. We note that, in \cite{AC}, Amini and Caporaso also define the virtual weightless metric graph 
$\Gamma^{\omega}_\epsilon$, where each attached loop is assigned length $\epsilon > 0$. 
In this paper, we only use the case of $\epsilon = 1$ (i.e., $\Gamma^{\omega} = \Gamma^{\omega}_1$). 

To a vertex-weighted graph $(G, \omega)$, one can naturally associate a 
{\em vertex-weighted metric graph} $(\Gamma, \omega)$. Indeed, we define $\Gamma$ to be the metric graph associated to $G$, where each edge of $G$ is assigned length $1$. We extend $\omega: V(G) \to \ZZ_{\geq 0}$ to $\omega: \Gamma \to \ZZ_{\geq 0}$ by assigning $\omega(v) = 0$ for any $v \in \Gamma \setminus V(G)$. Then $\Gamma^{\omega}$ is the metric graph associated to $G^{\omega}$ (i.e., each edge of $G^{\omega}$ is assigned length $1$), and we have 
$g(G^{\omega}) = g(G, \omega) = g(\Gamma^{\omega}) = g(\Gamma, \omega)$. 
 
Let $(\Gamma, \omega)$ be a vertex-weighted metric graph. 
We have the natural embeddings $\jmath: \Gamma \to \Gamma^{\omega}$ and 
$\jmath: \Gamma_\QQ \to \Gamma^{\omega}_\QQ$. 
Let $D \in \Div(\Gamma)$. Via $\jmath$, 
we have $\jmath_*(D) \in \Div(\Gamma^{\omega})$. 
The rank $r_{(\Gamma, \omega)}(D)$ of $D$ for $(\Gamma, \omega)$ is defined 
by 
\begin{equation}
 \label{eqn:r:G:w} 
  r_{(\Gamma, \omega)}(D) := r_{\Gamma^{\omega}}(\jmath_*(D)). 
\end{equation}

\begin{Remark}
Vertex-weighted graphs are generalization of finite graphs. Indeed, 
let $G$ be a finite graph with associated metric graph $\Gamma$. 
Let $\mathbf{0}: V(G) \to \ZZ_{\geq 0}$ be the zero function. Then 
$(G, \mathbf{0})$ is a vertex-weighted graph, and we have 
$r_{(\Gamma, \mathbf{0})}(D) = r_{\Gamma}(D)$ for any $D \in \Div(\Gamma)$. 
We will often identify a finite graph $G$ with the vertex-weighted graph 
$(G, \mathbf{0})$ equipped with the zero function $\mathbf{0}$. 
\end{Remark}

Vertex-weighted graphs naturally appear as the reduction graphs of $R$-curves, 
as we now explain. 
Let  $\KK$ be a complete discrete valuation field with ring of integers $R$ 
and algebraically closed residue field $k$ as in \S\ref{subsec:sp}. Let  $X$ be a geometrically irreducible smooth projective 
curve over $\KK$, and $\Xscr$ a semi-stable model of $X$ over $R$. 
Let $\Xscr_0$ be the special fiber of $\Xscr$. 
Recall from \S\ref{subsec:sp} that we have the dual graph 
$G$ of $\Xscr_0$. Let $v$ be a vertex of $G$, and let $C_v$ be the corresponding 
irreducible component of  $\Xscr_0$. We define $\omega(v)$ to be the geometric genus of $C_v$. 
Then $\omega: V(G) \to \ZZ_{\geq 0}$ is a vertex-weight function,  and we obtain a 
vertex-weighted graph  $(G, \omega)$. We call 
$(G, \omega)$ the (vertex-weighted) {\em reduction graph} of 
$\Xscr$. Compared with $G$, the vertex-weighted graph $(G, \omega)$ captures more information 
of $\Xscr$, encoding the genera of irreducible components of the special fiber. 

We remark that Amini and Caporaso \cite{AC} obtained the Riemann--Roch formula and the specialization lemma for 
vertex-weighted graphs. 

\medskip
In the rest of this subsection, we show some properties of divisors on vertex-weighted metric graphs. Let $(\Gamma, \omega)$ be a vertex-weighted metric graph. Let $\Gamma^{\omega}$ 
be the virtual weightless metric graph associated to $(\Gamma, \omega)$. 
Let $\jmath: \Gamma \to \Gamma^{\omega}$ be the natural embedding. Let 
 $\jmath_*: \Div(\Gamma) \to \Div(\Gamma^{\omega})$ be the induced 
 injective map. 

\begin{Lemma} 
\label{lem:compatibility:reducedness}
We keep the notation above. Let $D \in \Div ( \Gamma )$. 
\begin{enumerate}
\item
If $E\in \Div ( \Gamma )$ satisfies $D \sim E$ on $\Gamma$, then 
$\jmath_*(D) \sim \jmath_*(E)$ on $\Gamma^\omega$. 
\item
Fix a point $v_{0} \in \Gamma$. Then $D$ is a $v_{0}$-reduced divisor on 
$\Gamma$ if and only if $\jmath_*(D)$ is a $v_{0}$-reduced divisor on 
$\Gamma^\omega$.  
\item
$r_\Gamma(D) \geq 0$ if and only if $r_{( \Gamma , \omega )} (D) \geq 0$. 
\item
Let $e$ be a leaf edge of $\Gamma$ with leaf end $v$ such that 
$\omega(v) = 0$. Let $\Gamma_1$ be the metric graph obtained 
by contracting $e$ in $\Gamma$, 
and $\omega_1$ the restriction of $\omega$ to $\Gamma_1$. 
Let $\varpi_1 : \Gamma \to \Gamma_1$ be the retraction map. 
Then $r_{( \Gamma , \omega )} (D) = r_{( \Gamma_1 , \omega_1)} (\varpi_{1\ast}(D))$. 
\end{enumerate}
\end{Lemma} 

\Proof
(1) Let $f$ be a rational function on $\Gamma$ such that 
$D - E = \zero(f)$. 
For a virtual loop $C \subset \Gamma^\omega$ that is added at a vertex 
$v \in \Gamma$ 
with positive weight, 
we set $\widetilde{f}(w) = f(v)$ for any 
$w \in C$. Then we obtain a  rational function $\widetilde{f}$ on 
$\Gamma^\omega$. 
Since $\jmath_*(D) -  \jmath_*(E) =  \zero(\widetilde{f})$, 
we have $\jmath_*(D) \sim \jmath_*(E)$ on $\Gamma^\omega$. 

(2) 
By induction on the number of loops added to $\Gamma$, we may assume that 
$\Gamma^\prime$ is the one-point sum of $\Gamma$ and a loop $\ell$. We put 
$v := \Gamma \cap \ell$ and $\overset{\scriptsize{\circ}}{\ell} := \ell\setminus\{v\}$. 

First we show the ``only if'' part. Suppose that $A^\prime$ is a closed subset of $\Gamma^\prime$ with $v_0 \not\in A^\prime$. If $\partial A^\prime 
\cap \overset{\scriptsize{\circ}}{\ell}$ is non-empty, then any point 
$a^\prime \in \partial A^\prime \cap \overset{\scriptsize{\circ}}{\ell}$ 
is non-saturated for $\jmath_*(D)$ with respect to $A^\prime$. 
If $\partial A^\prime \cap 
\overset{\scriptsize{\circ}}{\ell} = \emptyset$, then we set $A := A^\prime \setminus 
\overset{\scriptsize{\circ}}{\ell}$. We regard $A$ as a closed subset of $\Gamma$. 
Since $D$ is $v_0$-reduced, we have a non-saturated 
point $a \in \partial A$ for $D$ with respect to $A$. 
Then $a$ is in $\partial A^\prime$ and 
is non-saturated for $\jmath_*(D)$ with respect to $A^\prime$.  
Thus $\jmath_*(D)$ is $v_0$-reduced on 
$\Gamma^\prime$. 

Next we show the ``if"' part. Suppose that $A$ is a closed subset of $\Gamma$ with $v_0 \not\in A$. If $v \in A$, then we put $A^\prime := A \cup \overset{\scriptsize{\circ}}{\ell}$. Then $A^\prime$ is a closed subset of $\Gamma^\prime$ with 
$v_0 \not\in A^\prime$. Since $\jmath_*(D)$ is $v_0$-reduced, there exists a non-saturated point $a^\prime \in \partial A^\prime$ for $\jmath_*(D)$ with respect to $A^\prime$. Since 
$a^\prime \not\in \overset{\scriptsize{\circ}}{\ell}$, we find that $a^\prime$ is in 
$ \partial A \subset \Gamma$ and is non-saturated for $D$ with respect to $A$. If $v \in \Gamma\setminus A$, then we regard $A$ as a closed subset of $\Gamma^\prime$. 
Since $\jmath_*(D)$ is $v_0$-reduced, there exists a non-saturated point 
$a \in \partial A$ in $\Gamma^\prime$ that is non-saturated for $\jmath_*(D)$ with respect to $A$. 
We find that $a \in \partial A$ in $\Gamma$ and that $a$ is non-saturated 
for $D$ with respect to $A$. Thus $D$ is $v_0$-reduced on $\Gamma$. 

(3) The ``only if'' part is obvious. Indeed, if there exists an 
effective divisor $D^\prime$ on $\Gamma$ with 
$D \sim D^\prime$, then, by (1), $\jmath_*(D^\prime)$ is 
an effective divisor on $\Gamma^\omega$ with $\jmath_*(D) \sim 
\jmath_*(D^\prime)$. Hence $r_{(\Gamma, \omega)}(D) := 
r_{\Gamma^\omega}(\jmath_*(D)) \geq 0$. 
We show the ``if'' part. Let $v_0$ be a point on $\Gamma$, and 
let $E$ be the $v_{0}$-reduced divisor linearly equivalent to $D$ on $\Gamma$. 
By (2), $\jmath_*(E)$ is a $v_{0}$-reduced divisor on $\Gamma^\omega$, 
and by (1), $\jmath_*(E) \sim \jmath_*( D)$ on $\Gamma^\omega$.
Since $r_{( \Gamma , \omega )} (D) \geq 0$, Theorem~\ref{thm:HKN} tells us that $\jmath_*(E)$ is effective, and thus $E$ is also effective. 

(4) The retraction map $\varpi_1$ extends to the retraction map 
$\varpi_1^\omega: \Gamma^\omega \to \Gamma_1^\omega$, where $e \;(\subset 
\Gamma \subset \Gamma^\omega)$ is contracted. 
Let $\jmath_1: \Gamma_1 \hookrightarrow \Gamma_1^\omega$ be the natural embedding. Then Lemma~\ref{lemma:linear:equiv} implies that 
\[
r_{( \Gamma , \omega )} (D) 
= r_{\Gamma^\omega}(\jmath_*(D)) 
= r_{\Gamma_1^\omega}\left(\varpi^\omega_{1*}(\jmath_*(D))\right) 
= r_{\Gamma_1^{\omega_1}}\left(\jmath_{1*}(\varpi_{1*}(D))\right)
= r_{( \Gamma_1 , \omega_1)} (\varpi_{1\ast}(D)), 
\]
which completes the proof. 
\QED

\setcounter{equation}{0}
\section{Hyperelliptic graphs} 
\label{sec:hyp:graph}
In this section, we put together some properties of hyperelliptic metric graphs and 
hyperelliptic vertex-weighted graphs.  
We also define a quantity $p_\Gamma(D)$ (resp. $p_{(\Gamma, \omega)}(D)$) for a divisor $D$ on a hyperelliptic metric graph $\Gamma$ (resp. a hyperelliptic vertex-weighted metric graph $(\Gamma, \omega)$), which will play an important role in this paper. 

\subsection{Hyperelliptic metric graphs}
We recall some properties of hyperelliptic metric graphs. 
We refer the reader to \cite{BN2} and \cite{Chan} for details. 

We recall the definition of hyperelliptic metric graphs. 

\begin{Definition}[Hyperelliptic metric graph, cf. {\cite[\S~5.1]{BN2}} 
and {\cite[Definition~2.3]{Chan}}]
\label{def:hyp:metric:graph}
A compact connected metric graph $\Gamma$ is said to 
be \emph{hyperelliptic}  
if the genus of $\Gamma$ is at least $2$ and 
there exists a divisor on $\Gamma$ of degree $2$ and rank $1$. 
\end{Definition}

\begin{Definition}[Hyperelliptic finite graph, cf. {\cite[\S~5.1]{BN2}} 
and {\cite[Definition~2.3]{Chan}}]
\label{def:hyp:graph}
Let $G$ be a finite graph, and let $\Gamma$ be the metric graph associated to $G$. 
A graph $G$ is said to be \emph{hyperelliptic} if 
$\Gamma$ is hyperelliptic. 
\end{Definition}

Originally, in \cite{BN2}, Baker and Norine define the notion of hyperelliptic graphs 
for {\em loopless} finite graphs $G$ by the existence of a divisor of degree $2$ and rank $1$. This condition is equivalent to the metric graph $\Gamma$ associated to $G$ being hyperelliptic.  However, for a finite graph $G$ with a {\em loop}, this equivalence does not hold.
In this paper, we adopt the above definition of hyperelliptic finite graphs, 
for we consider finite graphs with loops in general. 

Let $\langle \iota \rangle$ be the group of order $2$ with generator $\iota$. 
We say that $\langle \iota \rangle$ acts non-trivially on $\Gamma$ if there exists an injective group homomorphism $\langle \iota \rangle \to \Isom(\Gamma)$, where $\Isom(\Gamma)$ is the group of isometries of $\Gamma$. Let $\Gamma/\langle \iota \rangle$ denotes the metric graph defined as the topological quotient with quotient metric. 
(Notice that our $\Gamma/\langle \iota \rangle$ is a little different from the one given in \cite[\S2.2]{Chan}, where certain leaf edges are removed from $\Gamma/\langle \iota \rangle$ for the compatibility with the {loopless} quotient graph $G/\langle\iota\rangle$ defined in  \cite[\S5.2]{BN2}.)

\begin{Definition}[Hyperelliptic involution]
\label{def:hyp:invol}
Let $\Gamma$ be a compact connected metric graph 
of genus at least $2$. A {\em hyperelliptic} involution of $\Gamma$ is 
an $\langle \iota \rangle$-action on $\Gamma$
such that $\Gamma/\langle \iota \rangle$ is a tree.
\end{Definition}

First we study the action of involution on bridges. 

\begin{Lemma} 
\label{lemma:characterization:of:bridge}
Let $\Gamma$ be a compact connected metric graph 
of genus at least $2$ without points of valence $1$.
Assume that $\Gamma$ has a hyperelliptic involution $\iota$. 
Let $e$ be an edge of $\Gamma$ with endpoints $v_1$ and $v_2$. Assume that 
$e$ is not a loop. 
Then $e$ is a bridge if and only if $\iota(e) = e$ and 
$\iota(v_i) = v_i$ for $i = 1, 2$. 
\end{Lemma}

\Proof
Recall that an edge of $\Gamma$ means an edge of the canonical model of $\Gamma$, which is regarded as a closed subset of $\Gamma$ (i.e., including the endpoints). For a bridge $e$ of $\Gamma$ with endpoints $v_1$ and $v_2$, we set $\overset{\scriptsize{\circ}}{e}  = e \setminus\{v_1, v_2\}$ as before. 

We first show the ``if'' part. Let $e$ be an edge of $\Gamma$ such that 
$\iota(e) = e$ and $\iota(v_i) = v_i$ for $i = 1, 2$. 
Since $\langle \iota \rangle$-action on $e$ is trivial
and $\Gamma / \langle \iota \rangle$ is a tree,
the metric graph $\Gamma \setminus \overset{\scriptsize{\circ}}{e}$ 
is not connected.
Thus $e$ is a bridge.

Next we show the ``only if'' part. Let $e$ be a bridge with 
endpoints $v_1$ and $v_2$. 
Then one has $\Gamma \setminus \overset{\scriptsize{\circ}}{e} = \Gamma_{1} \amalg \Gamma_{2}$ (disjoint union), 
where $\Gamma_{1}$ and $\Gamma_{2}$ are the connected components such that $v_{1} \in \Gamma_{1}$ and $v_{2} \in \Gamma_{2}$. 
Since $\Gamma$ does not have points of valence $1$, 
each $\Gamma_{i}$ is not a point and has at most one point of valence $1$. 
In particular,  $\Gamma_i$ is not a tree. 

Let us show that $\iota (e) = e$.
To argue by contradiction, suppose that $\iota (e) \neq e$.
Then, without loss of generality, 
we may assume 
that $\iota (e) \subseteq \Gamma_{2}$. 
Then $\iota (e) \cap \Gamma_{1} = \emptyset$.
It follows that $e \cap \iota(\Gamma_1) = \emptyset$. 
Since $\iota(\Gamma_1)$ is connected and $e \cap \iota(\Gamma_1) = \emptyset$, 
we have either $\iota(\Gamma_1) \subseteq \Gamma_1$ or $\iota(\Gamma_1) \subseteq \Gamma_2$. The former does not occur. Indeed, if   
$\iota(\Gamma_1) \subseteq \Gamma_1$, then $\iota(\Gamma_1) = \Gamma_1$ (we apply $\iota$), which leads to $\emptyset = e \cap \iota(\Gamma_{1}) = e \cap \Gamma_{1} = \{ v_{1} \} \neq \emptyset$, a contradiction. Thus we have $\iota(\Gamma_1) \subseteq \Gamma_2$, so that 
$\iota ( \Gamma_{1}) \cap \Gamma_{1} = \emptyset$. 
Since $\Gamma / \langle \iota \rangle$
is a tree, $\Gamma_{1}$ is a tree. 
This is a contradiction. We conclude that $\iota (e) = e$. 

It remains to show that $\iota (v_{1}) = v_{1}$ and $\iota (v_{2}) = v_{2}$. 
It suffices to show $\iota (v_{1}) = v_{1}$, which amounts to  
$\iota ( \Gamma_{1}) = \Gamma_{1}$. 
If $\iota ( \Gamma_{1}) \neq \Gamma_{1}$, then 
the above argument implies that $\Gamma_{1}$ is a tree, 
which is a contradiction as before. 
This completes the proof.
\QED

The following theorem relates hyperelliptic metric graphs and hyperelliptic involutions. 

\begin{Theorem}[{\cite[Proposition~5.5~and~Theorem~5.12]{BN2}}, {\cite[Corollary~3.9~and~Theorem~3.13]{Chan}}]
\label{thm:HMG}
Let $\Gamma$ be a compact connected metric graph 
with genus at least $2$ without points of valence $1$. 
Then the following are equivalent\textup{:} 
\begin{enumerate}
\item[(i)]
$\Gamma$ is hyperelliptic\textup{;}
\item[(ii)]
$\Gamma$ has a hyperelliptic involution. 
\end{enumerate}
Further, a hyperelliptic involution is unique. 
\end{Theorem}

\Proof
By Lemma~\ref{lemma:linear:equiv} and Lemma~\ref{lemma:characterization:of:bridge}, 
we may assume that $\Gamma$ is bridgeless. For the bridgeless case, 
see \cite[Proposition~5.5~and~Theorem~5.12]{BN2} and \cite[Corollary~3.9~and~Theorem~3.13]{Chan}.  
\QED

\begin{Remark}
The uniqueness of hyperelliptic involution for hyperelliptic
graphs is shown in \cite[Corollary~3.9]{Chan}.
The proof there is based on 
\cite[Proposition~3.8]{Chan},
and the proof of \cite[Proposition~3.8]{Chan}
uses the Riemann--Roch formula on metric graphs.
(The idea of the proof is the same as that of \cite[Proposition~5.5]{BN2}.)
Since we would like to give a proof of the Riemann--Roch formula on a loopless hyperelliptic graph by applying Theorem~\ref{thm:main} and Proposition~\ref{prop:relation}, and since Theorem~\ref{thm:HMG} will be used in the proof of  Theorem~\ref{thm:main}, we remark here that one can give a proof of the uniqueness of hyperelliptic involution free from the Riemann--Roch formula. 

The idea is as follows (we leave the details to the interested readers). Suppose that $\iota, \iota^\prime$ are involutions on $\Gamma$. 
If $\Gamma$ has a bridge $e$, then any point on $e$ is fixed by $\iota$ and $\iota^\prime$ by Lemma~\ref{lemma:characterization:of:bridge}. Thus contracting $e$, we may assume that $\Gamma$ is bridgeless. 
Then one can find a point $v \in \Gamma$ such that 
 $\iota(v) =  \iota^\prime(v)$. 
Now let $x \in \Gamma$ be an arbitrary point. 
Since $\Gamma/\langle \iota \rangle$ and 
$\Gamma/\langle \iota^\prime \rangle$ are trees  
and since any two points in a tree are linearly equivalent to each other, 
we have $[\iota (v)] + [v] \sim [\iota (x)] + [x]$
and $[\iota^\prime (v)] + [v] \sim [\iota^\prime (x)] + [x]$. 
It follows from $[\iota (v)] + [v] = [\iota^\prime (v)] + [v]$ that 
$[\iota (x)] + [x] \sim [\iota^\prime (x)] + [x]$ and thus  
$[\iota (x)] \sim [\iota^\prime (x)]$. 
Since $\Gamma$ is bridgeless, we then have 
$\iota (x) = \iota^\prime (x)$. We obtain $\iota = \iota^\prime$. 
\end{Remark}

The following lemmas show the compatibility of the notion of being hyperelliptic 
under a contraction. 

\begin{Lemma}
\label{lemma:hyp:cont}
Let $\Gamma$ be a compact connected metric graph. 
Suppose that $\Gamma$ has a bridge, and let $\Gamma_1$ be the graph 
obtained by contracting a bridge. Then  $\Gamma$ is hyperelliptic if and only if 
$\Gamma_1$ is hyperelliptic. 
\end{Lemma}

\Proof
This follows from Lemma~\ref{lemma:linear:equiv} and the 
definition of a hyperelliptic metric graph. 
\QED

Let $\Gamma$ be a hyperelliptic metric graph. 
Let $\Gamma^\prime$ be 
the metric graph obtained by 
contracting all the leaf edges of $\Gamma$. 
By Lemma~\ref{lemma:hyp:cont}, $\Gamma^\prime$ is a hyperelliptic 
metric graph. By Theorem~\ref{thm:HMG}, $\Gamma^\prime$ has the 
hyperelliptic involution $\iota^\prime: \Gamma^\prime \to \Gamma^\prime$. 
We denote by 
  $\varpi: \Gamma \to \Gamma^\prime
$
the retraction map, which induces 
$\varpi_*: \Div(\Gamma) \to \Div(\Gamma^\prime)$.  
We have the natural 
embedding $\Gamma^\prime \hookrightarrow \Gamma$, 
and we regard $\Gamma^\prime$ as a subgraph of $\Gamma$. 

\begin{Lemma}
\label{lemma:iota:iota}
Let $\Gamma^\prime$ be as above, and 
let $v, w \in \Gamma^\prime$. Then 
$[v] + [\iota(v)] \sim [w] + [\iota(w)]$ as divisors on $\Gamma^\prime$. 
Further, $[v] + [\iota(v)] \sim [w] + [\iota(w)]$ as divisors on $\Gamma$. 
\end{Lemma}

\Proof
Let $\overline{\Gamma}$ be the metric graph contracting all the bridges 
of $\Gamma^\prime$ and let $\overline{\varpi}^\prime: \Gamma^\prime \to \overline{\Gamma}$ be the retraction map. 
By Lemma~\ref{lemma:hyp:cont}, $\overline{\Gamma}$ is a hyperelliptic 
metric graph. 
By Lemma~\ref{lemma:characterization:of:bridge}, the action $\iota^\prime$ on 
$\Gamma^\prime$ descends to an action $\overline{\iota}$ on 
$\overline{\Gamma}$, which gives the hyperelliptic involution of 
$\overline{\Gamma}$. 
Since $\overline{\varpi}^\prime(v) + \overline{\iota}(\overline{\varpi}^\prime(v)) \sim \overline{\varpi}^\prime(w) + \overline{\iota}(\overline{\varpi}^\prime(w))$ as divisors on $\overline{\Gamma}$ by \cite[Theorem~3.2 and its proof]{Chan} (see also \cite[Corollary~5.14]{BN2}), we have $[v] + [\iota(v)] \sim [w] + [\iota(w)]$ 
as divisors on $\Gamma^\prime$ by 
Lemma~\ref{lemma:linear:equiv}(3). 
The second assertion follows from 
Lemma~\ref{lemma:linear:equiv}(1). 
\QED

\subsection{Hyperelliptic vertex-weighted graphs}
\label{subsec:p:hyp:vertex}
We recall some properties of hyperelliptic vertex-weighted graphs studied by Caporaso \cite{CaGonal}. We also introduce hyperelliptic vertex-weighted metric graphs and see some of their properties. 
Since our focus on this paper is to prove Theorem~\ref{thm:main}, we restrict our attention to the necessary properties, which will be used later. 

\begin{Definition}[Hyperelliptic vertex-weighted metric graph]
\label{def:hyp:vw}
Let $(\Gamma, \omega)$ be a vertex-weighted metric graph. 
We say that $(\Gamma, \omega)$ is {\em hyperelliptic} 
if the genus of $(\Gamma, \omega)$ is at least $2$ and 
there exists a divisor $D$ on $\Gamma$ 
such that $\deg(D) = 2$ and $r_{(\Gamma, \omega)}(D) = 1$. 
\end{Definition}

\begin{Definition}[Hyperelliptic vertex-weighted graph, cf. {\cite{CaGonal}} and Definition~\ref{def:hyp:graph}]
\label{def:hyp:vw:2}
Let $(G, \omega)$ be a vertex-weighted graph, 
and $\Gamma$ the metric graph associated to $G$. 
We say that $(G, \omega)$ is {\em hyperelliptic} 
if $(\Gamma, \omega)$ is hyperelliptic. 
\end{Definition}

Let $(G, \omega)$ be a vertex-weighted graph, 
and $\Gamma$ the metric graph associated to $G$. Let $\Gamma^{\omega}$ 
be the virtual weightless metric graph associated to $(\Gamma, \omega)$. 
Recall that we have the natural embedding 
$\jmath: \Gamma \to \Gamma^{\omega}$ and that we denote by 
$\jmath_*: \Div(\Gamma) \to \Div(\Gamma^{\omega})$ 
the induced injective map. 

The following proposition is a metric graph version of \cite[Lemma~4.1]{CaGonal}. 

\begin{Proposition}
\label{prop:hyp:weightless}
With the above notation, 
$(\Gamma,  \omega)$ is hyperelliptic if and only if
$\Gamma^{\omega}$ is hyperelliptic.
\end{Proposition}

\Proof
The ``only if'' part is obvious. Indeed, suppose that $(\Gamma , \omega )$ is hyperelliptic, 
and we take a divisor $D$ on $\Gamma$ with $\deg(D) = 2$ 
and $r_{(\Gamma , \omega )} (D) = 1$. Since $r_{(\Gamma , \omega )} (D) = 1$ 
means by definition $r_{\Gamma^\omega} (\jmath_* (D)) = 1$, we see that 
$\jmath_* (D) \in \Div ( \Gamma^{\omega} )$ is a divisor
with $\deg\jmath_* (D) = 2$ and $r_{\Gamma^\omega} (\jmath_* (D)) = 1$. 
Thus $\Gamma^{\omega}$ is hyperelliptic.

We show the ``if'' part. 
Suppose that $\Gamma^{\omega}$ is hyperelliptic.
If $\omega$ is trivial, then there is nothing to prove, so 
that we assume that there exists a point
$v_{1} \in \Gamma$ 
with $\omega (v_{1}) > 0$.
We put $D := 2 [ v_{1} ] \in \Div(\Gamma)$. We are going to show that 
$r_{(\Gamma, {\omega})} (D) = 1$.

Let $\overline{\Gamma^{\omega}}$
be the metric graph obtained from
$\Gamma^{\omega}$
by contracting all the bridges, and 
let $\varpi^{\omega} : \Gamma^{\omega} \to \overline{\Gamma^{\omega}}$ be the retraction map. 
By Lemma~\ref{lemma:hyp:cont} and Theorem~\ref{thm:HMG}, 
$\overline{\Gamma^{\omega}}$ is a hyperelliptic metric graph,
and let $\iota^{\omega}$ be the hyperelliptic 
involution of $\overline{\Gamma^{\omega}}$. 
By Lemma~\ref{lemma:iota:iota}, 
the divisor $D' := [\varpi^{\omega}(v_{1})] + [\iota^{\omega} (\varpi^{\omega}(v_{1}) )] \in \Div (\overline{\Gamma^{\omega}})$ has rank $1$. 
Since we have added loops at $v_1$, the vertex 
$v_{1}$ is a cut-vertex of $\Gamma^{\omega}$. 
Then $\varpi^\omega(v_{1})$ is a cut-vertex of $\overline{\Gamma^{\omega}}$.
We then have
$\iota^{\omega} ( \varpi^{\omega} (v_{1}) ) = \varpi^{\omega} (v_{1})$
by \cite[Lemma~3.9]{Chan}, 
so that 
$\varpi^{\omega}_* (\jmath_*(D)) = 2 [ \varpi^{\omega} (v_{1}) ]= D'$.
It follows that 
$r_{\overline{\Gamma^{\omega}}} (\varpi^\omega_* (\jmath_*(D))) = 1$,
and thus 
$r_{\Gamma^{\omega}} (\jmath_*(D)) = 1$ 
by Lemma~\ref{lemma:linear:equiv}.
We obtain
$r_{(\Gamma, {\omega})} (D) =r_{\Gamma^{\omega}} (\jmath_*(D)) = 1$.  
\QED

The next proposition is a metric graph version of \cite[Lemma~4.4]{CaGonal}, 
and gives a vertex-weighted version of Theorem~\ref{thm:HMG}. 

\begin{Proposition}
\label{prop:HMG:2}
Let $(G, \omega)$ be a vertex-weighted graph of genus at least $2$. 
Assume that any leaf end $v$ of $G$ satisfies $\omega(v) > 0$.  Let 
$\Gamma$ be the metric graph associated to $G$, and 
$\Gamma^\omega$ the virtual weightless metric graph of $(\Gamma, \omega)$. 
Then the following are equivalent\textup{:}
\begin{enumerate} 
\item[(i)]
$(\Gamma, \omega)$ is hyperelliptic\textup{;} 
\item[(ii)]
$\Gamma^\omega$ has a unique hyperelliptic involution. 
\end{enumerate} 
Further, the hyperelliptic involution preserves $\Gamma$, where  
 $\Gamma$ is seen as a subgraph of $\Gamma^\omega$ via the natural embedding $\Gamma \hookrightarrow \Gamma^\omega$. 
\end{Proposition}

\Proof
By the assumption on $(G, \omega)$, $\Gamma^\omega$ has no points of valence $1$. Thus the condition (ii) is equivalent to $\Gamma^\omega$ being hyperelliptic, 
which is equivalent to the condition (i) (see 
Theorem~\ref{thm:HMG} and 
Proposition~\ref{prop:hyp:weightless}). 

Let $\iota^\omega$ denote the hyperelliptic involution of $\Gamma^\omega$. 
Let $C$ be a virtual loop which is added at a vertex $v \in V(G)$ with 
$\omega(v) >0$.  
To show that $\iota^\omega(\Gamma) = \Gamma$,  it suffices to show that 
$\iota^\omega(C) = C$. Since $v$ is a cut-vertex of $\Gamma^\omega$ and 
any cut-vertex is $\iota^\omega$-fixed by \cite[Lemma~3.10]{Chan}, 
we have $\iota^\omega(v) = v$. Then $\iota^\omega(C)$ is a loop containing $v$. 
If  $\iota^\omega(C) \neq C$, then $\Gamma^\omega/\langle\iota^\omega\rangle$ 
has a loop corresponding to $C$, which is impossible. Thus  $\iota^\omega(C) = C$ and $\iota^\omega(\Gamma) = \Gamma$. 
\QED

\begin{Definition}[Hyperelliptic involution on a hyperelliptic vertex-weighted graph]
\label{def:HMG:2}
Let $(G, \omega)$ be a hyperelliptic vertex-weighted graph such that 
any leaf end $v$ of $G$ satisfies $\omega(v) > 0$, and let $\Gamma$ be the metric graph associated to $G$. Let $\iota: \Gamma \to \Gamma$ be the involution defined by the restriction of the hyperelliptic involution of $\Gamma^\omega$ 
to $\Gamma$ (cf. Proposition~\ref{prop:HMG:2}). We call 
$\iota$ the hyperelliptic involution of $(\Gamma, \omega)$. 
\end{Definition}

Since $\Gamma/\langle \iota \rangle$ is a subtree of  
$\Gamma^\omega/\langle \iota^\omega \rangle$, 
the above definition agrees with 
Definition~\ref{def:hyp:invol}.

\subsection{Quantities $p_\Gamma(D)$ and $p_{(\Gamma, \omega)}(D)$}
\label{subsec:p:Gamma:D}
We introduce a quantity $p_\Gamma(D)$ for 
a divisor $D$ on a hyperelliptic metric graph $\Gamma$. 
We also introduce $p_{(\Gamma, \omega)}(D)$ for a divisor $D$ on 
hyperelliptic vertex-weighted metric graph $(\Gamma, \omega)$. 
The quantities 
$p_\Gamma(D)$ and $p_{(\Gamma, \omega)}(D)$ 
will play important roles in this paper. 

Let $\Gamma$ be a hyperelliptic metric graph. 
Let $\Gamma^\prime$ be 
the metric graph obtained by 
contracting all the leaf edges of $\Gamma$. 
We denote by 
  $\varpi: \Gamma \to \Gamma^\prime
$
the retraction map, which induces 
$\varpi_*: \Div(\Gamma) \to \Div(\Gamma^\prime)$. 

Since $\Gamma^\prime$ is hyperelliptic by Lemma~\ref{lemma:hyp:cont}, 
$\Gamma^\prime$ has a unique hyperelliptic involution $\iota^\prime$ 
by Theorem~\ref{thm:HMG}. 
We fix a point $v_0 \in \Gamma^\prime$ with 
\begin{equation}
\label{eqn:v0}
\iota^\prime(v_0) = v_0
\end{equation}
We note that such $v_0$ always exists (see Lemma~\ref{lemma:cont:p} below). 
We regard $v_0$ as an element of $\Gamma$ via the natural 
embedding $\Gamma^\prime \hookrightarrow \Gamma$. 
For an effective divisor $D$ on $\Gamma$, we set 
\begin{equation}
\label{eqn:p}
  p_\Gamma(D) = \max\{
  r \in \ZZ_{\geq 0} \mid \left|D - 2r [v_0]\right| \neq \emptyset
  \}. 
\end{equation}

We put together several results that will be used later. 

\begin{Lemma}
\label{lemma:cont:p}
Let $\Gamma, \Gamma^\prime$ and $\varpi$ be as above. 
\begin{enumerate}
\item
There exists $v_0 \in \Gamma^\prime$ with 
$\iota^\prime(v_0) = v_0$. 
\item
The quantity $p_\Gamma(D)$ defined in \eqref{eqn:p} 
is independent of the choice of 
$v_0 \in \Gamma^\prime$ with $\iota^\prime(v_0) = v_0$. 
\item
Let $D$ be an effective divisor $D$ on $\Gamma$, 
and let $D_{v_0}$ be the $v_0$-reduced divisor linearly equivalent to $D$. 
Then $p_\Gamma(D) = \left\lfloor \frac{D_{v_0}(v_0)}{2} \right\rfloor$. 
\item
For any effective divisor $D$ on $\Gamma$, we have 
$p_\Gamma(D) = p_{\Gamma^\prime}(\varpi_*(D))$. 
\end{enumerate}
\end{Lemma}

\Proof
(1)
Recall that $\langle\iota^\prime\rangle$ acts non-trivially on $\Gamma^\prime$ and that 
$T^\prime := \Gamma^\prime/\langle\iota^\prime\rangle$ is a tree. 
Let $\pi : \Gamma^\prime \to T^\prime$ be the quotient map.
Take a leaf end $\pi (v_0) \in T^\prime$.
If $\pi^{-1} (\pi (v_{0}))$ consists of two points,
then these two points should be leaf ends of $\Gamma^\prime$,
but that contradicts the assumption on $\Gamma^\prime$.
Thus $\pi^{-1} (\pi (v_{0})) = \{ v_0 \}$,
which shows that $\iota^\prime (v_0) = v_0$. 

(2) For $w \in \Gamma^\prime$, Lemma~\ref{lemma:iota:iota} tells us that 
$2 [v_0] \sim [w] + [\iota^\prime(w)]$ in $\Div(\Gamma)$. Thus 
\begin{equation}
\label{eqn:p:g:1}
  p_\Gamma(D) = \max\{r \in \ZZ_{\geq 0} \mid \left|
  D - r \left([w] + [\iota(w)] \right)
  \right| \neq \emptyset \}. 
\end{equation}
Suppose that $\widetilde{v_0} \in \Gamma^\prime$ is another point 
with $\iota^\prime(\widetilde{v_0}) = \widetilde{v_0}$. Then, setting 
$w = \widetilde{v_0}$ in \eqref{eqn:p:g:1}, 
we obtain the assertion. 

(3)
We set $s = \left\lfloor \frac{D_{v_0}(v_0)}{2} \right\rfloor$. Then 
$D_{v_0} - 2s [v_0]$ is a $v_0$-reduced effective divisor, so that 
$p_{\Gamma}(D) \geq s$. On the other hand, $D_{v_0} - 2(s+1) [v_0]$ 
is a $v_0$-reduced divisor with negative coefficient at $v_0$. Hence 
$|D_{v_0} - 2(s+1) [v_0]| = \emptyset$, so that $p_{\Gamma}(D) < s+1$. 
We conclude $p_{\Gamma}(D) = s$. 

(4) 
We note that $D - 2r [v_0] \sim \varpi_*(D) - 2r [v_0]$ in $\Div(\Gamma)$ 
by Lemma~\ref{lemma:linear:equiv}(2), from which the assertion follows. 
\QED

Now let $(\Gamma, \omega)$ be a hyperelliptic vertex-weighted metric graph. Let 
$\Gamma^{\omega}$ be the virtual weightless metric graph 
of $(\Gamma, \omega)$.  By Proposition~\ref{prop:hyp:weightless}, $\Gamma^{\omega}$ is a hyperelliptic metric graph. Let $\jmath: \Gamma \hookrightarrow \Gamma^{\omega}$ be the natural 
embedding.  For an effective divisor $D \in \Div(\Gamma)$, we set 
\begin{equation}
\label{eqn:p:G:w}
p_{(\Gamma, \omega)}(D) 
:= p_{\Gamma^{\omega}}(\jmath_*(D)).   
\end{equation}

\setcounter{equation}{0}
\section{Hyperelliptic semi-stable curves} 
\label{sec:hyp:semistable:curves}
In this section, we study hyperelliptic semi-stable curves, and show Theorem~\ref{thm:lifting} via the equivariant deformation based on \cite[Theorem~4.8]{CaGonal}. As we write in the introduction, there is another approach to Theorem~\ref{thm:lifting} due to Amini--Baker--Brugall\'e--Rabinoff \cite[Theorem~1.10]{ABBR}. 

\subsection{Hyperelliptic semi-stable curves}
\label{subsec:hyp:ss}
Let $\Omega$ be an algebraically closed field with ${\rm char}(\Omega) \neq 2$. 
Let $\OO$ be an $\Omega$-algebra.
We call $\OO$ a node 
if there is an isomorphism $\OO \cong \Omega [[ x,y ]] / (xy)$
as an $\Omega$-algebra.
Let $X_{0}$ be an algebraic scheme of dimension $1$ over $\Omega$
and let $c \in X_{0}$ be a closed point.
We call $c$ a node 
if the complete local ring $\widehat{\OO_{X_{0},c}}$ 
is a node in the above sense.
A \emph{semi-stable} curve is a 
connected reduced proper curve over $\Omega$  which 
has at most nodes as singularities. 
A \emph{stable curve} over $\Omega$ is a semi-stable curve
with ample dualizing sheaf.
Recall that $\langle \iota \rangle$ denotes the group of order~$2$. 

\begin{Definition}[Hyperelliptic curve]
\label{def:hyp:ss}
A semi-stable (resp. stable) curve $X_{0}$ over $\Omega$
with an $\langle \iota \rangle$-action
on $X_{0}$
is called a \emph{hyperelliptic} semi-stable (resp. stable)  curve 
if 
\begin{enumerate}
\item[(i)]
for any irreducible component $C$ of $X_{0}$ with
$\iota (C) = C$,
the $\langle \iota \rangle$-action restricted to $C$ is nontrivial (i.e., not the identity), 
and 
\item[(ii)]
$X_{0} / \langle \iota \rangle$ is a semi-stable curve of arithmetic 
genus $0$.
\end{enumerate}
\end{Definition}

\begin{Definition}[Hyperelliptic $S$-curve]
\label{def:hyp:S}
\begin{enumerate}
\item
Let $\mathscr{X} \to S$ be a proper and flat morphism over a scheme $S$. We say that $\mathscr{X}$ is a semi-stable $S$-curve (resp. a stable $S$-curve) if, for any geometric point $\overline{s}$ of $S$, the geometric fiber $\mathscr{X}_{\overline{s}}$ is a semi-stable curve (resp. a stable curve). 
\item
A  semi-stable (resp. stable) $S$-curve $\mathscr{X}$ 
equipped with an $\langle \iota \rangle$-action on $\mathscr{X}/S$ is called
a \emph{hyperelliptic}  semi-stable (resp. stable) $S$-curve  
if any geometric fiber of $\mathscr{X}_{\overline{s}}$ equipped with the restriction of the  
$\langle \iota \rangle$-action
is a hyperelliptic semi-stable curve.
\end{enumerate}
\end{Definition}

As in the introduction, let $\KK$ be a complete discrete valuation field with 
ring of integers $R$ and algebraically closed residue field $k$ such that 
${\rm char}(k) \neq 2$. 

\begin{Proposition} 
\label{properness}
Let $\mathscr{X}$ be a semi-stable $R$-curve
whose generic fiber is a smooth hyperelliptic curve $X$.
Assume that there exists an $\langle \iota \rangle$-action on 
$\mathscr{X}/\Spec(R)$ such that the restriction of $\iota$ to the generic fiber 
is the hyperelliptic involution on $X$. 
Then $\mathscr{X}$ equipped with the $\langle \iota \rangle$-action 
is a hyperelliptic semi-stable $R$-curve.
\end{Proposition}

\Proof
Let $X_{0}$ denote the the special fiber of $\mathscr{X} \to \Spec(R)$. 
Let $C$ be an irreducible component of $X_{0}$ such that 
with $\iota (C) = C$. We show that the $\langle \iota \rangle$-action 
on $C$ is nontrivial. 
Let $q : \mathscr{X} \to \mathscr{Y}$ be the quotient by $\iota$. 
Then, $q_{\ast} \OO_{\mathscr{X}}$ is a coherent 
$\OO_{\mathscr{Y}}$-module of rank $2$.
Let $\eta$ be the generic point of $C$. 
Then we have 
\[
\dim q^{-1} ( q ( \eta ) ) =
\dim_{\kappa ( q ( \eta ) )} q_{\ast} ( \OO_{\mathscr{X}} ) 
\otimes \kappa ( q (\eta) ) \geq 2, 
\]
where $\kappa ( q ( \eta ) )$ is the residue field at $q ( \eta )$.

On the other hand, since ${\rm char}(k) \neq 2$,  
the order $2$ of the action is invertible in $R$. 
Hence the restriction of $q$ to the special fiber
coincides with the quotient $X_{0} \to X_{0} / \langle \iota \rangle$. 
Since $\eta \in C$ and $\dim q^{-1} ( q ( \eta ) ) \geq 2$,
the $\langle \iota \rangle$-action on $C$ is not trivial.

It follows from \cite[Proposition~1.6]{Sa} that $\mathscr{Y} \to \Spec (R)$ 
is semi-stable. 
Since $\mathscr{Y} \to \Spec (R)$ is flat and since the arithmetic genus 
of the generic fiber of $\mathscr{Y} \to \Spec (R)$ 
is $0$, the arithmetic genus is of the special fiber 
$X_0/\langle \iota \rangle $ is also $0$. We obtain that 
$X_0/ \langle \iota \rangle $ is a semi-stable curve of genus $0$. 
\QED

\subsection{Equivariant specialization}
In this subsection, we prove Theorem~\ref{thm:lifting}. 
Let $\KK, R$ and $k$ be as in Theorem~\ref{thm:main}. 

Let $( G , \omega)$ be a vertex-weighted graph, and 
let $\Gamma$ be the metric graph associated to $G$. 
Let $(G_{\omega \circ}, \ell)$ be the model of $\Gamma$ 
with the set of vertices  
\[
  V(G_{\omega \circ}) 
  = 
  \{
  v \in V(G) 
  \mid
  \text{$w (v) > 0$ or $\val(v) \neq 2$} 
  \}. 
\] 
We define the vertex-weight function $\omega: V(G_{\omega \circ}) \to \ZZ_{\geq 0}$ by the restriction to vertex-weight function $\omega: V(G) \to \ZZ_{\geq 0}$ to $V(G_{\omega \circ})$. 
We call $(G_{\omega \circ}, \ell, \omega)$ 
the \emph{vertex-weighted canonical model} of $(\Gamma , \omega)$, and 
call $(G_{\omega \circ}, \omega)$ the underlining vertex-weighted graph of 
the canonical model of $(\Gamma , \omega)$. 

The following characterization is proved by Caporaso \cite{CaGonal}. 

\begin{Theorem}[{\cite[Theorem~4.8]{CaGonal}}]
\label{thm:model:of:configuration}
Let $(G , \omega)$ be a hyperelliptic vertex-weighted graph of genus $g$. 
Assume that any leaf end of $v$ of $G$ satisfies  $\omega (v) > 0$. 
Let $\Gamma$ be the metric graph associated to $G$, and 
$(G_{\omega \circ} , \omega )$ the underlining vertex-weighted graph of 
the canonical model of $(\Gamma , \omega)$. Then the following are equivalent. 
\begin{enumerate}
\item
For any $v \in V (G_{\omega \circ})$, 
there are at most
$(2 \omega (v) + 2)$
positive-type bridges emanating from $v$.
\item
There exists a hyperelliptic stable curve
$X_{0}$ of genus $g$ such that
\begin{enumerate}
\item[(i)]
the (vertex-weighted) dual graph of $X_{0}$ is $(G_{\omega \circ} , \omega )$, and
\item[(ii)] 
the $\langle \iota \rangle$-action on $X_{0}$ is compatible with
the hyperelliptic involution on $(\Gamma , \omega )$
in the following sense\textup{:}
For any $v \in V ( G_{\omega \circ} )$,
we have $\iota (C_v) = C_{\iota (v)}$,
where $C_v$ denotes the irreducible component of $X_0$ 
corresponding to $v$\textup{;} 
For any $e \in E ( G_{\omega \circ} )$,
we have $\iota ( p_{e} ) = p_{\iota (e)}$,
where $p_{e}$ is the node of $X_0$ corresponding to $e$.
\end{enumerate}
\end{enumerate}
\end{Theorem}

Based on Theorem~\ref{thm:model:of:configuration}, we use 
the equivariant deformation to show the existence of a regular model ~$\Xscr$. 

\begin{Theorem}
\label{thm:sec:4}
Let $(G, \omega)$ be a hyperelliptic vertex-weighted graph of genus $g(G, \omega) \geq 2$ such that, for every vertex $v$ of $G$, there are at most $(2 \omega(v) + 2)$ positive-type bridges emanating from $v$. Assume that any leaf end $v$ of $G$ satisfies $\omega(v) > 0$. 
Let $\Gamma$ be the metric graph associated to $G$. 
Then there exists a regular, generically smooth, semi-stable $R$-curve $\Xscr$ with reduction graph 
$(G, \omega)$ such that the generic fiber $X$ of $\Xscr$ is hyperelliptic. 
Further, for the specialization map $\tau: X(\overline{\KK}) \to \Gamma_\QQ$, we have $\tau \circ \iota_X = \iota \circ \tau$, where $\iota_X$ is the hyperelliptic involution of $X$, and $\iota$ is the hyperelliptic involution of~$\Gamma$. 
\end{Theorem}

\Proof
Let $(G_{\omega \circ}, \ell, \omega)$ be the vertex-weighted canonical model of $( \Gamma , \omega )$.
We take a hyperelliptic stable curve $X_{0}$ as in Theorem~\ref{thm:model:of:configuration}.
Let $p_{1} , \ldots , p_{r}$ be the 
$\langle\iota\rangle$-fixed nodes of $X_{0}$
and let $p_{r+1}, \ldots , p_{r+s}$ be the
nodes such that
$p_{r+1}, \ldots , p_{s} , \iota (p_{r+1}) , \ldots , \iota (p_{r+s} )$
are the distinct non-$\langle\iota\rangle$-fixed nodes.

For $1 \leq i \leq r+s$, let $\Def_{p_{i}}$ denote the deformation functor 
for the node $\widehat{\OO_{X_0, p_i}}$ (see \S\ref{subsec:A:2} for details). 
Let $\Phi_{\iota}^{gl}:
\Def_{(X_{0} , \iota)} \to 
\prod_{i = 1}^{r} \Def_{p_{i}} \times
\prod_{i = r + 1}^{r+s} \Def_{p_{i}}$ be the $\langle\iota\rangle$-equivariant 
global-local morphism, which assigns, to any $\langle\iota\rangle$-equivariant deformation of $X_0$, the deformation of the node at $p_i$ 
for $1 \leq i \leq r+s$ (see \S\ref{subsec:A:3} for details). 

Let $\pi$ be a uniformizer of $R$. For a functor $F$, we set 
$\widehat{F}(R) := \varprojlim_{n} F(R/\pi^n)$. 
For $1 \leq i  \leq r+s$, 
let $d_{i}$ be an element 
in $\widehat{\Def_{p_{i}}} ( R )$
that has a representative of form
\begin{align*}
\begin{CD}
\widehat{\OO_{X_{0} , p_{i}}} @<<< R [[ x , y ]] / ( xy - \pi^{\ell_{i}} ) \\
@AAA @AAA \\
k @<<< R
,
\end{CD}
\end{align*}
where $\ell_{i}$ is the length of the edge of $G_{\omega \circ}$ corresponding
to $p_{i}$.

We set 
$
d := (d_{i}) 
\in
\left(
\prod_{i = 1}^{r} \widehat{\Def_{p_{i}}} \times
\prod_{i = r + 1}^{r+s} \widehat{\Def_{p_{i}}} 
\right)
( R )
$. By Corollary~\ref{cor:surjectivitiy:EGLM},
we find an $\langle\iota\rangle$-equivariant diagram
\[
\begin{CD}
X_{0} @>>> \bar{\Xscr} \\
@VVV @VVV \\
\Spec(k) @>>> \Spf(R)
\end{CD}
\]
whose isomorphism class in 
$\widehat{\Def_{(X_{0} , \iota)}}  (R)$
is a lift of $d$ by $\widehat{\Phi_{\iota}^{gl}} (R)$.
This diagram of formal curves is algebraizable
(cf. Remark~\ref{algebraizable}),
and we write for the algebrization $\bar{\mathscr{X}} \to \Spec(R)$.
Let $\mathscr{X} \to \Spec(R)$ be the minimal resolution
of $\bar{\mathscr{X}} \to \Spec(R)$.
Then $\mathscr{X} \to \Spec R$ has 
the vertex-weighted reduction graph $(G, \omega)$. 

It remains to show that the specialization map 
$\tau : X ( \overline{\KK} ) \to \Gamma_{\QQ}$
is compatible with
the hyperelliptic involutions.
To see that,
let $\KK'$ be a finite extension of $\KK$
and $R'$ be the ring of integer of $\KK'$. 
Let $e(\KK^\prime/\KK)$ denote the ramification index of $\KK^\prime/\KK$. 
Let $\bar\Xscr' \to \Spec (R')$ be the base-change of $\bar{\mathscr{X}} \to \Spec (R)$
to $\Spec (R')$
and let $\Xscr'$ be the minimal resolution
of $\bar\Xscr'$.
Then the vertex-weighted dual graph of the special fiber 
$\bar\Xscr' \to \Spec (R')$
equals $(G_{\omega \circ} , \omega)$. 
The vertex-weighted dual graph $(G' , \omega')$ of the special fiber of ${\Xscr'}
\to \Spec (R')$,
where each edge is assigned length $1 / e(\KK^\prime/\KK)$,
is a model of $(\Gamma , \omega)$.
The $\langle\iota\rangle$-action on $\bar\Xscr'$ lifts to that on $\Xscr'$, which we denote by 
$\iota_{\Xscr^\prime}$. 
Let $v'$ be a vertex of $G'$ and let $C'_{v'}$ be the corresponding irreducible
components in the special fiber of $\Xscr' \to \Spec (R')$.
Let $e$  be an edge of $G_{\omega \circ}$ with $v' \in e$
and $p_e$ the corresponding node of $X_0$.
From the construction of the hyperelliptic involution
on $X_{0}$ in Theorem~\ref{thm:model:of:configuration}, 
we have
$\iota_X ( p_e ) = p_{\iota (e)}$ and 
$\iota_{\Xscr^\prime}  (C'_{v'} ) = C'_{\iota (v')}$.

Let $P \in X ( \overline{\KK} )$ be a point
and take a finite extension $\KK'$ such that $P \in X ( \KK' )$.
Then the corresponding section of $\Xscr^\prime \to \Spec (R')$
intersects with a unique irreducible component $C'_{v'}$ for some $v' \in V (G')$.
We have $\tau (P) = v'$ by definition.
Since the section corresponding to $\iota_X (P)$
intersects with $\iota_{\Xscr^\prime} ( C'_{v'} )$ and since $\iota_{\Xscr^\prime}  (C'_{v'} ) = C'_{\iota (v')}$
as noted above,
we obtain $\tau ( \iota_X ( P )) = \iota (v')$.
\QED

We are ready to prove Theorem~\ref{thm:lifting}. 

\begin{CorollaryNoNum}[$=$ Theorem~\ref{thm:lifting}]
Let $(G, \omega)$ be a hyperelliptic vertex-weighted graph such that every vertex $v$ of $G$ has at most $(2 \omega(v) + 2)$ positive-type bridges emanating from $v$. 
Then there exists a regular, generically smooth, semi-stable $R$-curve $\Xscr$ with reduction graph 
$(G, \omega)$ such that the generic fiber $X$ of $\Xscr$ is hyperelliptic. 
\end{CorollaryNoNum}

\Proof
Successively contracting the leaf edges
with a leaf end $v$ of $G$ such that $\omega(v) = 0$, 
we obtain a vertex-weighted hyperelliptic graph $( \overline{G} , 
\overline{\omega} )$.
Then we apply Theorem~\ref{thm:sec:4} 
to obtain a desired regular, generically smooth, semi-stable
$R$-curve for $( \overline{G} , 
\overline{\omega} )$.
Taking successive blowing-ups,
we obtain a desired $R$-curve for $(G, \omega)$.
\QED

\setcounter{equation}{0}
\section{Reduced divisors on a (hyperelliptic) graph}
\label{sec:reduced}
In this section, we prove Theorem~\ref{thm:main:3} using the notion of moderators (see \cite[Theorem~3.3]{BN}, \cite[Section~7]{MZ}, \cite[Corollary~2.3]{HKN}).  
The proof of Theorem~\ref{thm:main:3} is due to the referees and is significantly simplified from the original version.  

We begin by recalling the definition of moderators and some of their properties. 
Let $\Gamma$ be a compact connected metric graph of genus $g \geq 2$. 
Let $G$ be a model of $\Gamma$ without loops. We give an orientation on $G$, so that each edge $e$ of $G$ has 
head vertex $h(e)$ and tail vertex $t(e)$. 
An orientation on $G$ is said to be {\em cyclic} if there exist  
edges $e_1, \ldots, e_k$ of $G$ such that $h(e_i) = t(e_{i+1})$ for 
$i = 1, \ldots, k-1$ and $h(e_{k}) = t(e_1)$. An orientation on $G$ 
is {\em acyclic} if it is not cyclic. 

\begin{Definition}[{\cite[Definition~7.8]{MZ}}]
A divisor $K_+ \in \Div(\Gamma)$ is called a {\em moderator} if 
there exist a model $G$ of $\Gamma$ without loops and an 
acyclic orientation on $G$ such that 
\[
  K_+ = \sum_{v \in V(G)} (\val_+(v)-1)[v], 
\]
where $\val_+(v)$ denotes the number of outgoing edges from $v$ 
with respect to the orientation. 
\end{Definition}

\begin{Proposition}[{\cite[Theorem~3.3]{BN}, \cite[Section~7]{MZ}, \cite[Corollary~2.3]{HKN}}]
\label{prop:MZ}
Let $\Gamma$ be a compact connected metric graph of genus $g \geq 2$. 
\begin{enumerate}
\item
Any moderator $K_+$ on $\Gamma$ has degree $g-1$. 
\item
Let $D \in \Div(\Gamma)$ be a $v_0$-reduced divisor on $\Gamma$ 
with $D(v_0) < 0$. Then there exists a $v_0$-reduced moderator 
$K_+$ such that $D \leq K_+$ and $K_+(v_0) = -1$.  
\end{enumerate}
\end{Proposition}

\Proof
See \cite[Sect.~3.2]{BN}, 
\cite[Proposition~7.9]{MZ} and \cite[Sect.~2.1]{HKN} for (1). 

The assertion (2) is proved in \cite[Theorem~3.3]{BN} and \cite[Section~7]{MZ}, \cite[Corollary~2.3]{HKN}. Because the formulation is slightly different, we recall how to construct $K_+$. 

We set $D^\prime := D - D(v_0)[v_0]$. Then 
$D^\prime$ is an effective $v_0$-reduced divisor. We take a sequence $\mathbf{a} = (a_1, a_2, \ldots, a_k)$ with $\Supp(D^\prime) = \{a_1, a_2, \ldots, a_k\}$ as in 
(the proof of) Theorem~\ref{thm:DL}. 
We put $a_0 := v_0$. 
We give an ordering on $\{v_0\} \cup \Supp(D^\prime)$ by defining 
$a_0 < a_1 < a_2 < \cdots < a_k$. 

Let $G_\circ$ be the canonical model of $\Gamma$. We make a new finite graph 
$G_\circ^\prime$ by adding the middle points of all loops of $G_\circ$ (if exist),  so that $G_\circ^\prime$  is a loopless finite graph. 
Let $V(G_\circ^\prime)$ be the set of vertices of $G_\circ^\prime$. 
We set 
\begin{align*}
  V & := \{v \in V(G_\circ^\prime) \mid \val(v) \geq 2, \; v \neq v_0,\; 
v \not\in \Supp(D^\prime) \}, \\
  W & := \{v \in V(G_\circ^\prime) \mid \val(v) = 1, \; v \neq v_0\}. 
\end{align*}
Note that, since $D^\prime$ is $v_0$-reduced, we have $W \cap \Supp(D^\prime) = \emptyset$. 

We are going to give an ordering on 
$\{v_0\} \cup \Supp(D^\prime) 
\cup V \cup W$  (disjoint union). For $1 \leq i \leq k$, 
let $\Ucal_i$ be the connected component of $\Gamma\setminus\{a_i, a_{i+1}, \ldots, a_k\}$ that contains $v_0$. 
We write $\Ucal_1 \cap V = \{b_{11}, b_{12}, \ldots, b_{1 j_1}\}$. We give an ordering $b_{11} < b_{12} < \cdots < b_{1 j_1}$ so that $b_{1 \alpha}$ is contained in the connected component of $\Ucal_1 \setminus \{b_{1 \alpha+1}, \ldots, b_{1  j_1}\}$ that contains $v_0$ for any $\alpha = 1, \ldots, j_1-1$. Then we define 
$a_0 < b_{11} < b_{12} < \cdots < b_{1 j_1} < a_1$. 
Suppose now that an ordering $a_{i-2} < b_{i-1\,1} < \cdots < b_{i-1\, j_{i-1}} < a_{i-1}$ is defined. 
Inductively, we write $\Ucal_i \cap \left(V \setminus \{b_{11}, b_{12}, \ldots, 
b_{i-1\, j_{i-1}-1}, b_{i-1\, j_{i-1}}\}\right) = \{b_{i1}, b_{i2}, \ldots, b_{i\, j_i}\}$. 
We give an ordering $b_{i1} < b_{i2} < \cdots < b_{i\, j_i}$ so that 
$b_{i \alpha}$ is contained in the connected component of $\Ucal_i \setminus \{b_{i\, \alpha+1}, \ldots, b_{i\, j_i}\}$  that contains $v_0$ 
for any $\alpha = 1, \ldots, j_i-1$. 
Then we define $a_{i-1} < b_{i1} < b_{i2} < \cdots < b_{i\, j_i} < a_i$. 
At the stage $k+1$, we write $V \setminus \{b_{11}, b_{12}, \ldots, 
b_{k\, j_{k}-1}, b_{k\, j_{k}}\}
= \{b_{k+1\, 1}, b_{k+1\, 2}, \ldots, b_{k+1\, j_{k+1}}\}$, 
and we give an ordering $b_{k+1\,1} < b_{k+1\,2} < \cdots < b_{k+1\,\, j_{k+1}}$ 
so that $b_{k+1\, \alpha}$ is contained in the connected component of $\Gamma\setminus \{b_{k+1\, \alpha+1}, \ldots, b_{k+1\, j_{k+1}}\}$  that contains $v_0$ for any $\alpha = 1, \ldots, j_{k+1}-1$. 
Then we define $a_k < b_{k+1\, 1} < b_{k+1\, 2} < \cdots < b_{k+1\, j_{k+1}}$. 
Finally we write $W = \{c_1, \ldots, c_\ell\}$ and define 
$b_{k+1\, j_{k+1}} < c_1 < \cdots < c_\ell$. In conclusion, 
we have given an ordering on 
$\{v_0\} \cup \Supp(D^\prime) \cup V \cup W$. 

Let $G$ be the model 
of $\Gamma$ whose vertices are given by 
$\{v_0\} \cup \Supp(D^\prime) \cup V \cup W$. 
For each edge of $e$ of $G$, we define the head vertex 
$h(e)$ of $e$ and the tail vertex of $t(e)$ of $e$ so that $h(e)$ is smaller than $t(e)$ 
with respect to the above ordering on $V(G)$. This gives an acyclic 
orientation on $G$. Let $K_+ \in \Div(\Gamma)$ be the moderator with respect to this orientation. Then $K_+$ is $v_0$-reduced (cf. Theorem~\ref{thm:DL}). 
Further, by the 
construction, $K_+(v_0) = -1$ and  
$D(w) \leq K_+(w)$ for any $w \neq v_0 \in \Gamma$. 
By the assumption of $D$, we have $D(v_0) \leq -1 = K_+(v_0)$. We conclude that  $D \leq K_+$ on $\Gamma$. Thus $K_+$ has all the desired properties. 
\QED


\begin{TheoremNoNum}[$=$ Theorem~\ref{thm:main:3}]
Let $\Gamma$ be a compact connected metric graph of genus $g \geq 2$. 
We fix a point $v_0 \in\Gamma$. 
Let $D \in \Div(\Gamma)$ be a $v_0$-reduced divisor 
on $\Gamma$. Then, if $\deg(D) - D(v_0)\leq g-1$, 
then there exists $w\in \Gamma \setminus\{v_0\}$ 
such that $D + [w]$ is a $v_0$-reduced divisor. 
\end{TheoremNoNum}


\Proof
We set $D^{\prime\prime} := D - (D(v_0)+1)[v_0] \in \Div(\Gamma)$. 
Since $D^{\prime\prime}$ is $v_0$-reduced 
and $D^{\prime\prime}(v_0) = -1$,  
Proposition~\ref{prop:MZ} tells us that there exists a $v_0$-reduced 
moderator $K_+$ such that 
\[
  D^{\prime\prime} \leq K_+
\]   
and $K_+(v_0) = -1$. 
Since $\deg\left(D^{\prime\prime}\right) \leq g-2$ and $\deg(K_+) = g-1$, 
there exists $w \in \Gamma$ such that 
$D^{\prime\prime} + [w] \leq K_+$. 
Since $D^{\prime\prime}(v_0) = K_+(v_0) = -1$, 
We have $w \neq v_0$.  

Since $D^{\prime\prime} + [w] \leq K_+$ and $K_+$ is $v_0$-reduced, 
$D^{\prime\prime} + [w]$ is $v_0$-reduced. It follows that 
$D + [w] = D^{\prime\prime} + [w] + (D(v_0) +1)[v_0]$ is $v_0$-reduced, which completes the proof of 
Theorem~\ref{thm:main:3}. 
\QED

\medskip
We have the following corollaries of Theorem~\ref{thm:main:3}, 
which will be needed to prove Theorem~\ref{thm:main:2}. 

\begin{Corollary}
\label{cor:hyp:graph:2}
Let $\Gamma$ be a hyperelliptic metric graph of genus $g$. 
Let $v_0$ be an element of $\Gamma$ satisfying \eqref{eqn:v0}.  
Let $D$ a $v_0$-reduced divisor on $\Gamma$. 
Assume that $p_\Gamma(D) = 0$ and $\deg(D) \leq g-1$. Then there exists 
a divisor $E$ on $\Gamma$ such that 
\[
 D \leq E, \quad \deg(E) = g, \quad p_\Gamma(E) = 0. 
\]  
\end{Corollary}

\Proof
Since $p_\Gamma(D)=0$, we have $D(v_0) \leq 1$.

\smallskip
{\bf Case 1.}\quad
Assume that $D(v_0)=0$. Using 
Theorem~\ref{thm:main:3} repeatedly, 
there exist $w_{\deg(D)+1}, \ldots$, $w_g \in \Gamma\setminus\{v_0\}$ such that 
$E := D + [w_{\deg(D)+1}] + \cdots +[w_g]$ is $v_0$-reduced. 
If $p_\Gamma(E) \geq 1$, then $E- 2[v_0]$ is linearly equivalent to an effective divisor. 
However, since $E -  2[v_0]$ is $v_0$-reduced and the coefficient at $v_0$ is $-2$, 
this is impossible. Thus we get $p_\Gamma(E) = 0$. 

\smallskip
{\bf Case 2.}\quad
Assume that $D(v_0)=1$. We put $D^\prime = D - [v_0]$. Then 
$D^\prime$ is $v_0$-reduced, and using 
Theorem~\ref{thm:main:3} repeatedly, 
there exist $w_{\deg(D^\prime)+1}, \ldots$, $w_{g-1} \in \Gamma\setminus\{v_0\}$ such that 
$E^\prime := D^\prime + [w_{\deg(D^\prime)+1}] + \cdots +[w_{g-1}]$ is $v_0$-reduced. 
Put $E = E^\prime +[v_0]$. Then $E$ is $v_0$-reduced, 
$D \leq E$ and $\deg(E) = g$.  Further, we obtain 
$p_{\Gamma} (E) = 0$ by the same argument as in Case 1. 
\QED


\begin{Corollary}
\label{cor:g}
Let $\Gamma$ be a hyperelliptic metric graph of genus $g$. 
Let $D$ be an effective divisor on $\Gamma$. 
Assume that $p_\Gamma(D) = 0$ and $\deg(D) = g$. 
Then $r_\Gamma(D) = 0$. 
\end{Corollary}

\Proof
Recall that we have fixed a point $v_0$ on $\Gamma$ satisfying \eqref{eqn:v0}. 
Let $D_0$ be the $v_0$-reduced divisor on $\Gamma$ which is linearly equivalent to $D$.  
Then $D_0$ is effective. Since $p_\Gamma(D)=0$, we have $D_0(v_0) \leq 1$. 
We may and do replace $D_0$ with $D$. 

\smallskip
{\bf Case 1.}\quad
Assume that $D(v_0)=0$. In this case, $D-[v_0]$ is also $v_0$-reduced and is not effective, 
so that $D-[v_0]$ is not linearly equivalent to an effective divisor. 
Thus $r_\Gamma(D-[v_0]) = -1$. Hence $r_\Gamma(D) \leq 0$. Since $D$ is effective, 
we have $r_\Gamma(D) = 0$. 

\smallskip
{\bf Case 2.}\quad
Assume that $D(v_0)=1$. 
We set $D^\prime = D - [v_0]$. Then $D^\prime$ is effective 
and $v_0$-reduced. Since $p_\Gamma(D) = 0$, 
we have $p_\Gamma(D^\prime) = 0$. 

By Theorem~\ref{thm:main:3}, there exists $w \in \Gamma\setminus\{v_0\}$ 
such that $D^\prime+[w]$ is $v_0$-reduced. 
To argue by contradiction, we assume that 
$r_\Gamma(D) \neq 0$. 
%
Since $r_\Gamma(D) \geq 1$, $D- [\iota(w)]$ is linearly equivalent to 
an effective divisor $D^{\prime\prime}$. By Theorem~\ref{thm:HKN}(2), we may 
assume that $D^{\prime\prime}$ is $v_0$-reduced. Then 
$D^{\prime\prime} + [v_0]$ is $v_0$-reduced. 
We have 
\begin{align*}
D^{\prime\prime} + [v_0] & 
\sim D- [\iota(w)] + [v_0] \sim (D^\prime +  [v_0] ) - [\iota(w)] + [v_0] \\
& \sim  D^\prime + 2 [v_0]- [\iota(w)] \sim 
D^\prime + ([w] + [\iota(w)]) - [\iota(w)] \sim 
D^\prime+ [w]. 
\end{align*}
Since $w$ is taken so that 
$D^\prime+ [w]$ is $v_0$-reduced, the uniqueness of 
$v_0$-reduced divisors (Theorem~\ref{thm:HKN}(1)) implies that 
$D^{\prime\prime} + [v_0] = D^\prime+ [w]$ in $\Div(\Gamma)$. 
However, the coefficient of $D^{\prime\prime} + [v_0]$ at $v_0$ is at least $1$, 
while that of $D^\prime+ [w]$ is $0$. This is a contradiction, and we obtain 
$r_\Gamma(D) = 0$. 
\QED

\setcounter{equation}{0}
\section{Rank of divisors on a hyperelliptic graph}
\label{sec:rank:graph}
In this section, we prove Theorem~\ref{thm:main:2}. We first note that 
Riemann's inequality on graphs, which is a weaker form of the Riemann--Roch theorem on graphs, is deduced from Baker's Specialization Lemma and Riemann's inequality on curves.  

\begin{Proposition}
\label{prop:Riemann}
Let $G$ be a finite graph of genus $g$ and $\Gamma$ the metric graph associated to $G$. Let $D$ be a divisor on $\Gamma$. Then we have 
$r_\Gamma(D) \geq \deg(D) -g$. 
\end{Proposition}

\Proof
By \cite[Proposition~3.2]{GK}, there exists a divisor $D^\prime \in \Div(\Gamma_\QQ)$ such that $\deg(D^\prime)  = \deg(D)$ and $r_\Gamma(D) = r_\Gamma(D^\prime)$. Replacing $D$ with $D^\prime$, we may assume that $D \in \Div(\Gamma_\QQ)$. Let $\KK$ be a complete discrete valuation field with ring of integers $R$ and algebraically closed residue field, and let $\Xscr$ be 
 a regular, generically smooth, semi-stable $R$-curve such that the reduction graph 
 equals $\Gamma$. Let $X$ be the generic fiber of $\Xscr$ and 
$\tau_*: \Div(X_{\overline{\KK}}) \to  \Div(\Gamma)$ the specialization map. 
We take $\widetilde{D} \in \Div(X_{\overline{\KK}})$ with 
$\tau_*(\widetilde{D}) = D$. Then Baker's Specialization Lemma (Theorem~\ref{thm:specialization:lemma}) and Riemann's inequality on $X$ give 
\[
r_\Gamma(D) \geq r_X(\widetilde{D}) \geq  
\deg(\widetilde{D}) - g = \deg(D) - g. 
\] 
This completes the proof. 
\QED

We prove Theorem~\ref{thm:main:2}, using  
Corollary~\ref{cor:hyp:graph:2}, Corollary~\ref{cor:g} and 
Proposition~\ref{prop:Riemann}. 
Recall that $r_{(\Gamma, \omega)}(D)$ and $p_{(\Gamma, \omega)}(D)$ are respectively 
defined in \eqref{eqn:r:G:w} and \eqref{eqn:p:G:w}. 

\begin{TheoremNoNum}[$=$ Theorem~\ref{thm:main:2}]
Let $(G, \omega)$ be a hyperelliptic vertex-weighted graph, and 
$\Gamma$ the metric graph associated to $G$. Set $g = g(\Gamma, \omega)$. 
Let $D$ be an effective divisor on $\Gamma$.
Then 
\[
 r_{(\Gamma, \omega)}(D) = \begin{cases}
  p_{(\Gamma, \omega)}(D) & \textup{(if $\deg(D) - p_{(\Gamma, \omega)}(D) \leq g$), } \\
  \deg(D) -g & \textup{(if $\deg(D) - p_{(\Gamma, \omega)}(D) \geq g+1$).} 
  \end{cases} 
\]  
\end{TheoremNoNum}

\Proof
{\bf Step 1.}\quad 
Let $G^{\omega}$ be the virtual weightless graph associated to $(G, \omega)$, and 
let $\Gamma^{\omega}$ be the virtual weightless metric graph associated to $(G, \omega)$. 
Note that $\Gamma^{\omega}$ is the metric graph associated to $G^{\omega}$. 
By Proposition~\ref{prop:hyp:weightless}, $\Gamma^{\omega}$ is a hyperelliptic graph. 
Let $\jmath: \Gamma \hookrightarrow \Gamma^{\omega}$ be the natural embedding. 
Since $g(\Gamma, \omega) = g(\Gamma^{\omega})$, 
$r_{(\Gamma, \omega)}(D) = r_{\Gamma^{\omega}}(\jmath_*(D))$ 
and $p_{(\Gamma, \omega)}(D) = p_{\Gamma^{\omega}}(\jmath_*(D))$
by definition, it suffices to prove the theorem for 
the weightless graphs, i.e., for $G^{\omega}$ and $\Gamma^{\omega}$. 

\medskip
{\bf Step 2.}\quad 
By Step 1, we replace $G^{\omega}$ by $G$, and $\Gamma^{\omega}$ by $\Gamma$. 
Let $\overline{\Gamma}$ be the metric graph obtained by 
contracting all the 
leaf edges of $\Gamma$, and $\varpi: \Gamma \to \overline{\Gamma}$ 
the retraction map. 
Since $r_\Gamma(D) = r_{\overline{\Gamma}}(\varpi_*(D))$  
by Lemma~\ref{lemma:linear:equiv}(2) and 
$p_\Gamma(D) = p_{\overline{\Gamma}}(\varpi_*(D))$ by 
Lemma~\ref{lemma:cont:p}(3) for any divisor $D$ on $\Gamma$, 
we may and do assume that $\Gamma$ has no points of valence $1$. 
Let $\iota$ be the  hyperelliptic involution of $\Gamma$ (cf. Theorem~\ref{thm:HMG}). We fix $v_0 \in \Gamma$ with  $\iota(v_0) = v_0$ (cf. Lemma~\ref{lemma:cont:p}). 

Let $D$ be an effective divisor on $\Gamma$. 
Let $D_0$ be the $v_0$-reduced divisor linearly equivalent to $D$. 
We set 
$r = \left\lfloor \frac{D_0(v_0)}{2} \right\rfloor$ and 
$s = \deg(D) - 2r$. 
Then $D_0$ is written as 
\[
  D_0 = 2r [v_0] + [w_1] + \cdots + [w_s]
\]
for some $w_1, \ldots, w_s \in \Gamma$. By Lemma~\ref{lemma:cont:p}(3), 
we have $p_\Gamma(D) = r$. 

If $\iota(w_i)  = w_j$ for some $i \neq j$, then $[w_i] + [w_j] \sim 2 [v_0]$ by 
Lemma~\ref{lemma:iota:iota},  
and $D_0 \sim 2(r+1) [v_0] + \sum_{k=1, k \neq i, j}^s [w_k]$.  
This contradicts $p_\Gamma(D) = r$. 
Thus $\iota(w_i) \neq w_j$ for any $i \neq j$. 
Also,  $p_\Gamma([w_1] + \cdots + [w_s]) = 0$ by 
Lemma~\ref{lemma:cont:p}(3).  

\medskip
{\bf Case 1.}\quad Assume that $\deg(D) - p_\Gamma(D) \leq g$. 
Note that $s \leq r+ s = \deg(D) - r \leq g$. 
Since $[w_1] + \cdots + [w_s]$ is $v_0$-reduced and 
$p_\Gamma([w_1] + \cdots + [w_s]) = 0$, 
Corollary~\ref{cor:hyp:graph:2} tells us that there exist $w_{s+1}, \ldots, 
w_g \in \Gamma$ with $p_\Gamma([w_1] + \cdots + [w_s]+[w_{s+1}]+ \cdots + [w_g]) = 0$. 
By Corollary~\ref{cor:g}, we have 
$r_\Gamma([w_1] + \cdots + [w_s]+[w_{s+1}]+ \cdots + [w_g]) = 0$. Thus  
$r_\Gamma([w_1] + \cdots + [w_s]+ [w_{s+1}] + \cdots + [w_{s+r}]) = 0$. 
Since $2 [v_0] \sim [v] + [\iota(v)]$ for any $v \in \Gamma$ by 
Lemma~\ref{lemma:iota:iota}, we have 
\begin{align*}
  D & \sim 2r [v_0] + [w_1] + \cdots + [w_s]  \\
  & \sim [w_1] + \cdots + [w_s] + [w_{s+1}] + \cdots + [w_{s+r}] 
  + [\iota(w_{s+1})] + \cdots + [\iota(w_{s+r})]. 
\end{align*}
Since $r_\Gamma(E) \leq r_\Gamma(E -[v]) + 1$ for any divisor $E$ and $v \in \Gamma$, 
we have 
\begin{equation}
\label{eqn:thm:1:8}
  r_\Gamma(D) \leq r_\Gamma([w_1] + \cdots + [w_s] + [w_{s+1}] + \cdots + [w_{s+r}]) + r  = r. 
\end{equation}
On the other hand,
for any $u_1 , \ldots , u_r \in \Gamma$,
we have
\begin{align*}
D - ( [u_1] + \cdots + [u_r]  ) 
& \sim
2r [v_0] - ( [u_1] + \cdots + [u_r]  ) + [w_1] + \cdots + [w_s] \\
& \sim
[\iota ( u_1 )] + \cdots + [\iota (u_r)] + [w_1] + \cdots + [w_s]
\end{align*}
by Lemma~\ref{lemma:iota:iota}.
This shows $r_\Gamma ( D ) \geq r$.
Thus we conclude that $r_\Gamma(D) = r$, which is the desired estimate when 
$\deg(D) - p_\Gamma(D) \leq g$. 

\medskip
{\bf Case 2.}\quad Assume that $\deg(D) - p_\Gamma(D) \geq g+1$. 

\smallskip
{\bf Subcase 2-1.}\quad Assume that $s \leq g$. 
Since $[w_1] + \cdots + [w_s]$ is $v_0$-reduced and 
$p_\Gamma([w_1] + \cdots + [w_s]) = 0$, 
Corollary~\ref{cor:hyp:graph:2} tells us that there exist $w_{s+1}, \ldots, 
w_g \in \Gamma$ with $p_\Gamma([w_1] + \cdots + [w_s]+[w_{s+1}]+ \cdots + [w_g]) = 0$. 
By Corollary~\ref{cor:g}, we have 
$r_\Gamma([w_1] + \cdots + [w_s]+[w_{s+1}]+ \cdots + [w_g]) = 0$. 
Recalling that $2 [v_0] \sim [v] + [\iota(v)]$ for any $v \in \Gamma$ by 
Lemma~\ref{lemma:iota:iota}, we have 
\begin{align*}
  D & \sim 2r [v_0] + [w_1] + \cdots + [w_s]  \\
  & \sim 2 (r+s -g) [v_0] + [w_1] + \cdots + [w_s] + [w_{s+1}] + \cdots + [w_{g}] 
  + [\iota(w_{s+1})] + \cdots + [\iota(w_{g})].
\end{align*}
As in \eqref{eqn:thm:1:8}, since $r+s = \deg(D) - p_\Gamma(D) \geq g+1$, 
we have 
\begin{align*}
  r_\Gamma(D) 
  & \leq r_\Gamma([w_1] + \cdots + [w_s] + [w_{s+1}] + \cdots + [w_{g}]) + 2r + s - g    \\
  & = 2r + s - g = \deg(D) - g. 
\end{align*}
Since the other direction $r_\Gamma(D) \geq \deg(D) - g$ is Riemann's inequality (Proposition~\ref{prop:Riemann}), 
we conclude that $r_\Gamma(D) = \deg(D) - g$. 

\smallskip
{\bf Subcase 2-2.}\quad Assume that $s \geq g+1$. 
Since $p_\Gamma([w_1] + \cdots + [w_s]) = 0$, we have 
$p_\Gamma([w_1] + \cdots + [w_g]) = 0$. 
By Corollary~\ref{cor:g}, we have 
$r_\Gamma([w_1] + \cdots + [w_g]) = 0$. 
As in \eqref{eqn:thm:1:8}, 
we have 
\[
  r_\Gamma(D) 
  \leq r_\Gamma([w_1] + \cdots + [w_g]) + 2r + s - g    
  = 2r + s - g = \deg(D) - g. 
\]
As in Subcase 2-1,  we have 
the other direction $r_\Gamma(D) \geq \deg(D) - g$ by Riemann's inequality. 
Thus $r_\Gamma(D) = \deg(D) - g$, which 
completes the proof of Theorem~\ref{thm:main:2}. 
\QED

\setcounter{equation}{0}
\section{Proofs of Theorem~\ref{thm:main} and 
Proposition~\ref{prop:relation}}
\label{sec:proof}
In this section, we prove Theorem~\ref{thm:main} and Proposition~\ref{prop:relation}  and give several examples. We also consider Question~\ref{q} for 
a vertex-weighted graph of genus $0$ or $1$. 

We begin by proving Theorem~\ref{thm:main}. 

\begin{Lemma}
\label{lemma:2:to:1}
The condition \textup{(ii)} implies the condition \textup{(i)}  
in Theorem~\ref{thm:main}. 
\end{Lemma}

\Proof
Let $(G, \omega)$ be a hyperelliptic vertex-weighted graph and $\Gamma$ 
the metric graph associated to $G$. 
By definition, there exists a divisor $D \in \Div(\Gamma)$ 
such that $\deg(D) = 2$ and $r_{(\Gamma, \omega)}(D) = 1$. In view of 
\cite[Proposition~3.1]{GK}, 
$D$ is taken in $\Div(\Gamma_\QQ)$. Assuming (ii), we take a regular, generically smooth, semi-stable $R$-curve $\Xscr$ with reduction graph $(G, \omega)$ and $\widetilde{D} \in \Div(X_{\overline{\KK}})$ such that $D = \tau_*(\widetilde{D})$ and $r_{(\Gamma, \omega)}(D) = r_X(\widetilde{D})$. (Here $X$ is the generic fiber of $\Xscr$ and $\tau$ is the specialization map.) It follows that $X$ is a hyperelliptic curve. 
Then Theorem~\ref{thm:lifting} tells us that $(G, \omega)$
satisfies the condition (i). 
\QED

We show that the condition (C') implies the condition (C) in the introduction. 

\begin{Lemma}
\label{lemma:prime}
Let $(G, \omega)$ be a vertex-weighted graph, and $\Gamma$ the metric graph associated to $G$. Assume that there exists a regular, generically smooth, semi-stable $R$-curve $\Xscr$ with reduction graph $(G, \omega)$ satisfying the condition \textup{(C')}. Then $\Xscr$ satisfies the condition \textup{(C)}. 
\end{Lemma}   

\Proof
Let $D \in \Div(\Gamma_\QQ)$. From the condition (C'), we infer that there exist  divisors $E \in \Div(\Gamma_\QQ)$ and $\widetilde{E} 
\in \Div(X_{\overline{\KK}})$ such that $D \sim E$, 
$\tau_*(\widetilde{E} ) = E$ and $r_{(\Gamma, \omega)}(E) = r_X(\widetilde{E})$. 
By \cite[Corollary~A.9]{B} for metric graphs, the restriction of the specialization map 
$\rest{\tau_*}{\Prin(X_{\overline{\KK}})}: \Prin(X_{\overline{\KK}}) \to \Prin(\Gamma_\QQ)$ is surjective, 
where  $\Prin(\Gamma_\QQ) := \Div(\Gamma_\QQ) \cap \Prin(\Gamma)$. Since $D-E \in \Prin(\Gamma_\QQ)$,  there exists a principal divisor $\widetilde{N}$ such that $\tau_*(\widetilde{N}) = D-E$. We set 
$\widetilde{D} = \widetilde{E} + \widetilde{N} \in \Div(X_{\overline{\KK}})$. 
Then $\widetilde{D}$ satisfies $D = \tau_*(\widetilde{D})$ and 
$r_{(\Gamma, \omega)}(D) = r_X(\widetilde{D})$. 
\QED

By Lemma~\ref{lemma:prime}, (iii) implies (ii) in Theorem~\ref{thm:main}. 
Thus it suffices to show that (i) implies (iii) in Theorem~\ref{thm:main}, 
which amounts to the following. 

\begin{Theorem}
\label{thm:1:3}
Let $(G, \omega)$  be a hyperelliptic vertex-weighted graph such that, for every vertex $v$ of $G$, there are at most $(2\omega(v)+2)$ positive-type bridges emanating from $v$. Let $\KK$ be a complete discrete valuation field with ring of integers $R$ and algebraically closed residue field $k$ with ${\rm char}(k) \neq 2$. Then there exists a regular, generically smooth, semi-stable $R$-curve $\Xscr$ with generic fiber $X$ and reduction graph $(G, \omega)$ which satisfies the following condition\textup{:} Let $\Gamma$ be the metric graph associated to $G$\textup{;} For any $D \in \Div(\Gamma_\QQ)$, 
there exist a divisor $E = \sum_{i=1}^k n_i [v_i] 
\in \Div(\Gamma_\QQ)$ that is linearly equivalent to $D$ 
and a divisor $\widetilde{E} = \sum_{i=1}^k n_{i} P_i \in \Div(X_{\overline\KK})$
such that $\tau(P_i) = v_i$ for any $1 \leq i \leq k$ and 
$r_{(\Gamma, \omega)}(E) = r_{X}(\widetilde{E})$. 
\end{Theorem}

Before proving Theorem~\ref{thm:1:3}, we give a formula for the ranks of divisors on hyperelliptic curves which corresponds to Theorem~\ref{thm:main:2}.  

\begin{Proposition}
\label{prop:Yamaki}
Let $F$ be a field and $\overline{F}$ an algebraic closure of $F$. 
Let $X$ be a connected smooth hyperelliptic curve of genus $g \geq 2$ defined over $F$, and let $\iota_X$ be the hyperelliptic involution of $X$.  Let $D$ be an effective divisor on $X_{\overline{F}}$. 
We express $D$ as 
\[
  D = P_1 + \cdots + P_r + \iota_X(P_1) + 
  \cdots + \iota_X(P_r) + Q_1 + \cdots + Q_s,  
\]
where $P_1, \ldots, P_r, Q_1, \ldots, Q_s \in X(\overline{F})$ and 
$\iota_X(Q_i) \neq Q_j$ for any $i \neq j$ with $1 \leq i, j \leq s$. 
Then we have 
\[
  r_X(D) = \begin{cases}
  r & \textup{(if $\deg(D) - r \leq g$)},  \\
  \deg(D) -g & \textup{(if $\deg(D) - r \geq g+1$)}.  
  \end{cases}
\]
\end{Proposition}

\Proof
We may and do assume that $F = \overline{F}$.
Let $K_X$ be a canonical divisor of $X$, and let $f: X \to \PP^{g-1}$ be the canonical map defined by the complete linear system $|K_X|$. 
We set $C = f(X)$,
and
let $H \in \Div(C)$ be a hyperplane section. 
Then the pull-back $f^* : |H| \to |K_{X}|$ is an isomorphism
between linear systems.
Since $X$ is hyperelliptic, we have 
$\deg (H) = g-1$. 

We put $E := f(P_1) + \cdots + f(P_r) + f(Q_1) + \cdots + f(Q_s) \in \Div(C)$. 
Then 
$\deg (H - E)=
g -1 - \deg (D) + r$.
We remark that the restriction of 
the pull-back map $f^*$ 
gives the 
isomorphism 
$\rest{f^*}{ |H - E|} : |H - E| \overset{\sim}{\to} |K_{X} - D|$.
Indeed,
since $f: X \to C$ is the quotient map of the hyperelliptic involution $\iota_X$ and since $\iota_X(Q_i) \neq Q_j$ for any $i \neq j$ with $1 \leq i, j \leq s$, 
we have, 
for any $H' \in |H|$, 
$f^* (H') \geq D$ if and only if $H' \geq E$. 

\medskip
{\bf Case 1.}\quad
Suppose that $\deg(D) - r \leq g -1$. 
Then
$\deg (H - E)
\geq 0$.
Since 
$C \cong \PP^{1}$,
it follows that
\[
\dim (|H - E|) = \deg (H-E)
= g -1 - \deg (D) + r
.
\]
Via the above identification $|H - E| \cong |K_{X} - D|$,
we obtain 
$
\dim(|K_X-D|) = g - 1 -\deg(D) + r$. 
Then the Riemann--Roch theorem tells us that
\[
  r_X(D) = 
\dim(|K_X-D|)
+ 1 - g + \deg(D) = r
,
\] 
which gives the desired equality for $\deg(D) - r \leq g-1$. 

\medskip
{\bf Case 2.}\quad
Suppose that $\deg(D) - r \geq g$. 
Then $\deg (H - E) < 0$, and hence 
$|K_{X} - D| \cong |H-E| = \emptyset$.
It follows from the Riemann--Roch theorem
that
$
  r_X(D) = \deg(D) -g
$. 
This gives the desired equality for $\deg(D) - r \geq g$. (We note that, 
if $\deg(D) - r = g$, then $ r_X(D) = \deg(D) -g = r$.)

This completes the proof. 
\QED

{\sl Proof of Theorem~\ref{thm:1:3}}.\quad
Let $g \geq 2$ denote the genus of $(G, \omega)$. 
If $e$ is a leaf edge with leaf end $v$ with $\omega(v) = 0$, then 
we contract $e$. Let $G^\prime$ be the graph obtained by 
successively contracting all such leaf edges. 
Then $G^\prime$ is a finite graph such that any leaf edge of $G^\prime$ (if exists) 
has an leaf end $v$ with $\omega(v) >0$. 
We note that $G^\prime$ is seen as a subgraph of $G$. 
Let $(G^\prime, \omega^\prime)$ be the vertex-weighted graph, 
where the vertex-weight function is given by the restriction of 
$\omega$ to $V(G^\prime)$. 

Let $\Gamma^\prime$ be the metric graph associated to  $G^\prime$. 
By Proposition~\ref{prop:HMG:2}, $\Gamma^\prime$ has the 
hyperelliptic involution $\iota^\prime: \Gamma^\prime \to 
\Gamma^\prime$ (see Definition~\ref{def:HMG:2}). 
We remark that $\Gamma^\prime$ is 
naturally seen as a subset of $\Gamma$. 

We take a regular, generically smooth, semi-stable $R$-curve  
$\Xscr^\prime$ as in Theorem~\ref{thm:sec:4}. 
In particular, the generic fiber $X$ of $\Xscr^\prime$ is a hyperelliptic curve, and 
the dual graph of the special fiber equals $(G^\prime, \omega^\prime)$. 
Further, we have 
$\tau^\prime\circ\iota_X = \iota^\prime\circ\tau^\prime$ 
for the specialization map $\tau^\prime: X(\overline\KK) \to \Gamma^\prime$ and the hyperelliptic involution $\iota_X: X \to X$. 
We take a Weierstrass point $P_0^\prime \in X(\overline\KK)$, i.e., a point 
satisfying $\iota_X(P_0^\prime) = P_0^\prime$, and put  
$v_0^\prime = \tau^\prime(P_0^\prime) \in \Gamma^\prime_\QQ$. 
Then we have $\iota^\prime(v_0^\prime) = v_0^\prime$. 

As we have seen in the proof of Theorem~\ref{thm:lifting} (Corollary of Theorem~\ref{thm:sec:4}), by successively 
blowing up at closed points on the special fiber, we obtain a 
regular, generically smooth, semi-stable $R$-curve $\Xscr$ 
such that the dual graph of the special fiber equals $(G, \omega)$. 
We are going to show that $\Xscr$ has the desired properties.  

Let $\tau: X(\overline\KK) \to \Gamma_\QQ$ be the specialization map defined by 
$\Xscr$. Let $\jmath: \Gamma^\prime \hookrightarrow \Gamma$ be the natural embedding and $\varpi: \Gamma \to \Gamma^\prime$ the natural retraction. 
Then we have $\tau^\prime = \varpi\circ \tau$. 

\smallskip
{\bf Case 1.}\quad 
Suppose that $r_{(\Gamma, \omega)}(D) = -1$. 
We put $E := D$, and write $E = \sum_{i=1}^k n_i [v_i] 
\in \Div(\Gamma_\QQ)$. Take any $P_i \in X(\overline{\KK})$ with 
$\tau(P_i) = v_i$ for $1 \leq i \leq k$ (cf. Proposition~\ref{prop:Baker:properties}(1)), and we set 
$\widetilde{E} = \sum_{i=1}^k n_{i} P_i \in \Div(X_{\overline\KK})$. 
We need to show that $r_{X}(\widetilde{E}) = -1$. To argue by contradiction, 
suppose that $r_{X}(\widetilde{E}) \geq 0$. 
Then there exists an effective divisor $\widetilde{F} \in \Div(X_{\overline\KK})$ 
with $\widetilde{E} \sim \widetilde{F}$. Then $\tau_*(\widetilde{F})$ is an effective divisor on $\Gamma$ and, by Proposition~\ref{prop:Baker:properties}, 
$D = \tau_*(\widetilde{E})  \sim \tau_*(\widetilde{F})$. This contradicts
our assumption that $r_{(\Gamma, \omega)}(E) = -1$ by 
Lemma~\ref{lem:compatibility:reducedness}. 
We obtain the assertion when $r_{(\Gamma, \omega)}(D) = -1$. 

\smallskip
{\bf Case 2.}\quad 
Suppose that $r_{(\Gamma, \omega)}(D) \geq 0$. By Lemma~\ref{lem:compatibility:reducedness}, we have $r_\Gamma(D) \geq 0$. 
We set $D^\prime = \varpi_*(D) 
\in \Div(\Gamma^\prime_\QQ)$. 
Let $E^\prime \in \Div(\Gamma^\prime_\QQ)$ be the $v_0^\prime$-reduced divisor that is linearly equivalent to $D^\prime$ on $\Gamma^\prime$.  By 
Lemma~\ref{lem:compatibility:reducedness} and Theorem~\ref{thm:HKN}, 
$E^\prime$ is an effective divisor. 

We set 
$r = \left\lfloor \frac{E^\prime(v_0^\prime)}{2} \right\rfloor$ and 
$s = \deg(E^\prime) - 2r$, then  
$E^\prime$ is written as 
\[
  E^\prime = 2r [v_0^\prime] + [w_1^\prime] + \cdots + [w_s^\prime]
\]
for some $w_1^\prime, \ldots, w_s^\prime \in \Gamma^\prime_\QQ$ such that 
$\iota^\prime(w_i^\prime) \neq w_j^\prime$ for $i \neq j$. 

We claim that $r = p_{(\Gamma^\prime, \omega^\prime)}(E^\prime)$. Indeed, 
let $\Gamma^{\prime\, {\omega^\prime}}$ be the virtual weightless metric graph associated to $(\Gamma^\prime, \omega^\prime)$ with  
hyperelliptic involution $\iota^{\prime\, {\omega^\prime}}$, 
and let 
$\jmath^{\prime\, {\omega^\prime}}: \Gamma^\prime \hookrightarrow \Gamma^{\prime\, {\omega^\prime}}$ be the natural embedding. 
By Lemma~\ref{lem:compatibility:reducedness}(2), 
$\jmath^{\prime\, {\omega^\prime}}_*(E^\prime) 
= 2r [v_0^\prime] + [w_1^\prime] + \cdots + [w_s^\prime]$ is a $v_0$-reduced 
divisor on $\Gamma^{\prime\, {\omega^\prime}}$, 
and $\iota^{\prime\, {\omega^\prime}}(w_i^\prime) \neq w_j^\prime$ for $i \neq j$ (cf. Definition~\ref{def:HMG:2}). By Lemma~\ref{lemma:cont:p}(3), 
we have $r = p_{\Gamma^{\prime\, {\omega^\prime}}}\left(\jmath^{\prime\, {\omega^\prime}}_*(E^\prime)\right)$. By definition, the right-hand side equals 
$p_{(\Gamma^\prime, \omega^\prime)}(E^\prime)$, and thus 
$r = p_{(\Gamma^\prime, \omega^\prime)}(E^\prime)$. 

By Proposition~\ref{prop:Baker:properties}(1), 
we take 
$Q_1, \ldots, Q_s \in  X(\overline{\KK})$ such that 
$\tau^\prime(Q_i) = w_i^\prime$ for $i = 1, \ldots, s$. 
Since $\tau^\prime\circ\iota_X = \iota^\prime\circ\tau^\prime$, 
we have $\iota_X(Q_i) \neq Q_j$ for $i \neq j$. 
We set $\widetilde{E} =  2r  P_0 + Q_1 + \cdots + Q_s 
\in \Div(X_{\overline\KK})$. Finally, we set 
$E = \tau_*(\widetilde{E}) = 
2r [\tau(P_0)] + [\tau(Q_1)] + \cdots + [\tau(Q_s)]
\in \Div(\Gamma_\QQ)$. 

We show that $E$ and $\widetilde{E}$ have desired properties. Indeed, since 
$\varpi_*(E) = \varpi_*(\tau_*(\widetilde{E})) = \tau^\prime_*(\widetilde{E}) 
= E^\prime \sim D^\prime = \varpi_*(D)$ on $\Gamma^\prime$, we have 
$E \sim D$ on $\Gamma$ by Lemma~\ref{lemma:linear:equiv}. 
By Theorem~\ref{thm:main:2} and Proposition~\ref{prop:Yamaki}, we then have 
\[
  r_{(\Gamma^\prime, \omega^\prime)}(E^\prime) = r_X(\widetilde{E}) = \begin{cases}
  r & \textup{(if $\deg(D) - r \leq g$)},  \\
  \deg(D) -g & \textup{(if $\deg(D) - r \geq g+1$)}.  
  \end{cases}
\]
By Lemma~\ref{lem:compatibility:reducedness}, 
we have $r_{(\Gamma, \omega)}(D) = 
r_{(\Gamma^\prime, \omega^\prime)}(D^\prime) 
=  r_{(\Gamma^\prime, \omega^\prime)}(E^\prime)$. 
Thus we obtain the assertion. 
\QED

Next we consider a vertex-weighted graph of genus $0$ or $1$. 

\begin{Proposition}
\label{prop:tree}
Let $\KK$ be a complete discrete valuation field with ring of integers $R$ and algebraically closed residue field $k$ with ${\rm char}(k) \neq 2$. 
Let $(G, \omega)$ be a vertex-weighted graph of genus $0$ or $1$, and $\Gamma$ 
the metric graph associated to $G$. Then there exists a regular, generically smooth, semi-stable $R$-curve $\Xscr$ with generic fiber $X$ and reduction graph $G$ which satisfies the condition \textup{(C')} in Theorem~\ref{thm:main}.  
\end{Proposition}

\Proof
{\bf Case 1.}\quad Suppose that $g(G, \omega) = 0$. This means that 
$\omega = \mathbf{0}$, and $G$ is a tree. 
There exists a regular, generically smooth, strongly semi-stable, 
totally degenerate $R$-curve $\Xscr$ with 
reduction graph $G$.  Let $X$ denote the generic fiber of $\Xscr$. Then 
$X_{\overline{\KK}} \cong \PP^1_{\overline{\KK}}$. 

Let $v_0$ be any vertex of $G$. 
Let $D$ be a divisor on $\Gamma_\QQ$.
Since $G$ is a tree, $D$ is linearly equivalent to 
$(\deg D) [v_0]$. It follows that $r_\Gamma(D) = \deg(D)$ if 
$\deg(D) \geq 0$ and that $r_\Gamma(D) = -1$ if $\deg(D) < 0$. 
Let $\widetilde{D}$ be any divisor on $X_{\overline{\KK}}$ 
such that $\tau_*(\widetilde{D}) = D$. Then $\deg(\widetilde{D}) = \deg(D)$ (cf. Proposition~\ref{prop:Baker:properties}(3)). 
Since $X_{\overline{\KK}} \cong \PP^{1}_{\overline{\KK}}$, 
we have $r_X(\widetilde{D}) = \deg(D)$ if 
$\deg(D) \geq 0$, and $r_X(\widetilde{D}) = -1$ if $\deg(D) < 0$. 
Thus we get $r_\Gamma(D) = r_X(\widetilde{D})$

\medskip
{\bf Case 2.}\quad Suppose that $g(G, \omega) = 1$. In this case, $\omega = \mathbf{0}$, or we have $\omega (v_1) = 1$ for some vertex $v_{1}$ of $G$ and $\omega (v) = 0$ for any other vertex $v$. 

{\bf Subcase 2-1.}\quad 
Suppose that $\omega = \mathbf{0}$. Then $g(\Gamma) = 1$.  
Let $D$ be a divisor on $\Gamma_\QQ$.
As in the Case 1 of the proof of Theorem~\ref{thm:1:3}, 
we may assume that $D$ is linearly equivalent to 
an effective divisor. Also, since the assertion is obvious if $D = 0$, 
we may assume that $\deg(D) \geq 1$. 

We note that if $\deg(D) \geq 2$, then $r_\Gamma(D) \geq 1$. Indeed, 
let $v$ be any point in $\Gamma$, and $D_v$ the $v$-reduced divisor 
that is linearly equivalent to $D$. Since $g(G) = 1$, the $v$-reduced divisor $D_v$ is of form $a [v] + b [w]$, where $a \in \ZZ$ and $b\in\{ 0,1 \}$. Since $\deg(D_v)  \geq 2$, it follows that $a \geq 1$ and thus $D_v - [v]$ is effective. Since 
$v$ is arbitrary, it follows that $r_\Gamma(D) \geq 1$. 

Repeating the above procedure, we obtain $r_\Gamma(D) \geq \deg(D)-1$. 
We claim that $r_\Gamma(D) = \deg(D)-1$. Indeed, if this is not the case, 
we will then have $\deg(D) [w_1] \sim \deg(D) [w_2]$ for any $w_1, w_2 \in \Gamma$, and thus $g(\Gamma) = 0$, which contradicts $g(\Gamma) = 1$. 

Let $\ell$ be the total length of the metric graph obtained by contracting all leaf edges of $\Gamma$. Notice that there exists an $R$-curve $\Xscr^\prime$ whose generic fiber $X$ is a smooth connected curve of genus $1$ and the special fiber is a geometrically irreducible rational curve with one node with multiplicity $\ell$. 
(For example, one takes 
$\Xscr^\prime = \Proj\left(R[x, y, z]/(y^2 z - x^3 - x z^2 -\pi^\ell z^3)\right)$, where $\pi$ is a uniformizer of $R$.) Then taking successive blow-ups on the special fiber, we have a regular, generically smooth, semi-stable $R$-curve $\Xscr$ such that the reduction  graph is $G = (G, \mathbf{0})$. 

Let $E$ be an effective divisor linearly equivalent to $D$. 
We write $E = \sum_{i=1}^k n_{v_i} [v_i]$ where $n_{v_i} \geq 0$ for all $i$. We take $\widetilde{E} = \sum_{i=1}^k n_{v_i} P_i$ such that $\tau(P_i) = v_i$ for $1 \leq i \leq k$. 
Since $\widetilde{E}$ is effective and $\deg(\widetilde{E}) > 0$, 
by the Riemann--Roch formula on $X$, we have $r_X(\widetilde{E}) =  \deg(\widetilde{E})-1$. Hence $r_\Gamma(E) = r_X(\widetilde{E})$. 
%

\medskip
{\bf Subcase 2-2.}\quad 
Suppose that there exists one vertex $v_1$ of $G$ with $\omega(v_1) = 1$ and $\omega(v) = 0$ for the other vertices. Let $\Gamma^\omega$ be the virtual weightless metric graph of $(G, \omega)$. Then $g(\Gamma^\omega) = 1$. 

As in the Case 1 of the proof of Theorem~\ref{thm:1:3}, 
we may assume that $D$ is linearly equivalent to 
an effective divisor. Also we may assume that $D \neq 0$, so that 
$\deg(D) \geq 1$. 
Let $E$ be an effective divisor linearly equivalent to $D$. 
Then the computation in the above subcase gives
$r_{(\Gamma, \omega)}(E) = r_{\Gamma^\omega}(E) = \deg(E) - 1$. 
Let $\Xscr^\prime$ be a regular $R$-curve  whose generic fiber $X$ and the special fiber are both smooth connected curves of genus $1$.
Then taking successive blow-ups on the special fiber, 
we have a regular, generically smooth, semi-stable $R$-curve $\Xscr$ of $X$ such that the reduction  graph is $(G, \omega)$. Then the argument in the above subcase shows that 
there exists $\widetilde{E} \in \Div(X_{\overline\KK})$ 
such that $\tau_*(\widetilde{E}) = E$  and 
$r_{(\Gamma, \omega)}(E) = r_X(\widetilde{E})$. 
\QED

Next we prove Proposition~\ref{prop:relation}.

\begin{PropositionNoNum}[$=$ Proposition~\ref{prop:relation}]
Let $G$ be a finite graph and $\Gamma$ the metric graph associated to $G$. Assume that there exist a complete discrete valuation field $\KK$ 
with ring of integers $R$, and a regular, generically smooth, strongly semi-stable, totally degenerate  $R$-curve $\Xscr$ with the reduction graph $G = (G, \mathbf{0})$ satisfying 
the condition \textup{(C)} in Question~\ref{q}. Then the Riemann--Roch formula on 
$\Gamma$ is deduced from the Riemann--Roch formula on~$X_{\overline\KK}$. 
\end{PropositionNoNum}

\Proof
We take any $D \in \Div(\Gamma_\QQ)$. 
By the condition (C), there exists $\widetilde{D} \in \Div(X_{\overline\KK})$
such that $r_\Gamma(D) = r_X(\widetilde{D})$ 
and $\tau_*(\widetilde{D}) = D$. 

By the Riemann--Roch formula on $X$, we have 
\[
  r_X(\widetilde{D}) - r_X(K_X-\widetilde{D}) 
  = 1 -g(X) + \deg(\widetilde{D}). 
\] 
Since $\Xscr$ is strongly semi-stable and totally degenerate, we have $g(X)= g(\Gamma)$. 
We have $\deg(\widetilde{D}) = \deg D$ (cf. Proposition~\ref{prop:Baker:properties}(3)). Further, by \cite[Lemma~4.19]{B}, we have $\tau(K_X) \sim K_\Gamma$. 
Then  
\[
  r_\Gamma(D) - r_X(K_X-\widetilde{D}) 
  = 1 -g(\Gamma) + \deg(D). 
\]
We put 
$\widetilde{\mathscr{D}} = \{
\widetilde{F} \in\Div(X_{\overline\KK}) \mid \tau_*(\widetilde{F}) \sim D
\}$. 
By the Riemann--Roch formula on $X$, 
we have 
\[
  \max_{\widetilde{F}\in\widetilde{\mathscr{D}}} \{ r_X(K_X-\widetilde{F}) \}
  = -1 +g(X) - \deg(\widetilde{D}) 
  + \max_{\widetilde{F}\in\widetilde{\mathscr{D}}} \{r_X(\widetilde{F}) \}. 
\] 
Since the right-hand side attains the maximum when $\widetilde{F}=\widetilde{D}$ by Baker's Specialization Lemma and our choice of 
$\widetilde{D}$, so does 
the left-hand side. By the condition (C)
and Baker's Specialization Lemma, the left-hand side equals  
$r_\Gamma(K_\Gamma-D)$. Hence we get 
$r_X(K_X-\widetilde{D}) = r_\Gamma(K_\Gamma- D)$, and thus 
\[
  r_\Gamma(D) - r_\Gamma(K_\Gamma- D) 
  = 1 -g(\Gamma) + \deg(D)
\]
The last equality is nothing but the Riemann--Roch formula on $\Gamma_\QQ$. 
Finally, by the approximation result by Gathmann--Kerber \cite[Proposition~1.3]{GK}, the Riemann--Roch formula on $\Gamma$ is deduced from that on $\Gamma_\QQ$. 
\QED

\begin{Remark}
\label{rmk:relation}
Let $G$ be a loopless hyperelliptic graph. 
Let $\overline{G}$ be the finite graph obtained 
by contracting all the bridges of $G$. Let $\Gamma$ and $\overline{\Gamma}$ be the metric graphs associated to $G$ and $\overline{G}$, respectively. By Theorem~\ref{thm:main} and Proposition~\ref{prop:relation}, the Riemann--Roch formula on $\overline{\Gamma}$ is deduced from the Riemann--Roch formula on a suitable hyperelliptic curve. Since the rank of divisors  is preserved under contracting bridges by \cite[Corollary~5.11]{B} and \cite[Lemma~3.11]{Chan} (cf. Lemma~\ref{lemma:linear:equiv}), the Riemann--Roch formula on $\Gamma$ is deduced. Since $r_G(D) = r_\Gamma(D)$ for $D \in \Div(G)$ by \cite{HKN}, the Riemann--Roch formula on $G$ is also deduced. 
\end{Remark}

We give some examples of ranks of divisors on 
metric graphs. 

\begin{Example}
\label{eg:2}
Let $G$ be the following graph of genus $g \geq 3$, 
where each vertex is given by a white circle or a black circle. Let $\Gamma$ be 
the metric graph associated to $G$. Let $D = [v_1] + [v_2]$. 
It is easy to see $r_\Gamma(D) = 1$. 
\[
\setlength\unitlength{0.08truecm}
  \begin{picture}(70,60)(0,0)
  \put(30, 10){\circle*{2}}
  \put(30, 50){\circle*{2}}
  \put(45, 30){\circle{2}} 
  \put(22.5, 30){\circle{2}}
  \put(37.5, 30){\circle{2}}
  \put(15, 30){\circle{2}} 
  \put(28.5, 30){$\ldots$} 
  \qbezier(30, 10)(60, 30)(30, 50)
  \qbezier(30, 10)(0, 30)(30, 50)
  \qbezier(30, 10)(45, 30)(30, 50)
  \qbezier(30, 10)(15, 30)(30, 50)
  \put(29, 5){$v_2$}
  \put(29, 53){$v_1$}
  \end{picture}
\]

We take a complete valuation field $\KK$ with ring of integers $R$ such that there exists a regular, generically smooth, strongly semi-stable, totally degenerate $R$-curve  $\Xscr$ such that the generic fiber $X$ is {\em non-hyperelliptic} and the dual graph of the special fiber equals $G$. 
There exists such $\Xscr$, see, e.g., \cite[Example~3.6]{B}. 

Let $\widetilde{D}$ be a divisor on $X_{\overline{\KK}}$ such that 
$\tau_*(\widetilde{D}) = D$. Then $\deg(\widetilde{D}) = 2$. Since 
$X$ is assumed to be non-hyperelliptic, we have $r_X(\widetilde{D}) \neq 1$. 
It follows that the condition (C) in Question~\ref{q} is not satisfied for this choice of 
$\Xscr$. (Indeed, we have to choose a model $\Xscr$ such that $X$ is hyperelliptic to make the condition (C) satisfied.)
\end{Example}

\begin{Example}
\label{eg:5}
Let $G$ be the following three petal graph of genus $3$, 
where each vertex is given by a white circle or a black circle. 
Let $\Gamma$ be the metric graph associated to $G$. Let $D = 2 [v_0]$. 
It is easy to see $r_\Gamma(D) = 1$.  Thus $\Gamma$ is a hyperelliptic graph. 
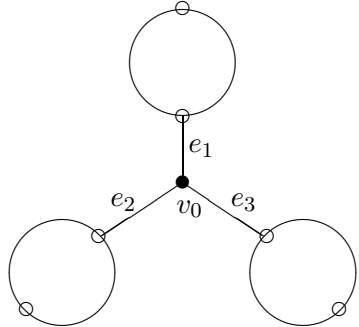
\[
  \setlength\unitlength{0.08truecm}
  \begin{picture}(70,70)(0,0)
  \put(30,55){\circle{20}}
  \put(10,20){\circle{20}}
  \put(50,20){\circle{20}}
  \put(30, 64){\circle{2}}
  \put(30, 46){\circle{2}}
  \put(4, 14){\circle{2}}
  \put(16, 26){\circle{2}}
  \put(56, 14){\circle{2}}
  \put(44, 26){\circle{2}} 
  \put(30, 35){\circle*{2}} 
  \put(30,35){\line(0,1){11}}
  \put(30,35){\line(3,-2){13.5}}
  \put(30,35){\line(-3,-2){13.5}}
  \put(29, 30){$v_0$}
  \put(31, 40){$e_1$}
  \put(18, 31){$e_2$}
  \put(38, 31){$e_3$}
  \end{picture}
\]

Let $\KK$ be a complete valuation field with ring of integers $R$ and algebraically residue field $k$ such that ${\rm char}(k) \neq 2$.  Let  $\Xscr$  be a regular, generically smooth, strongly semi-stable, totally degenerate $R$-curve with the reduction graph $G$. Let $X$ be the generic fiber of $\Xscr$. 

Since the vertex $v_0$ has three positive-type bridges $e_1, e_2, e_3$, the graph 
$G = (G, \mathbf{0})$ does not satisfy the condition (i) in Theorem~\ref{thm:main}. Then Theorem~\ref{thm:lifting}
tells us that $X$ is not hyperelliptic. 
The argument in Example~\ref{eg:2} (which agrees with 
Theorem~\ref{thm:main}) shows that 
there exists no divisor $\widetilde{D}$ on $X_{\overline{\KK}}$ with 
$r_X(\widetilde{D}) = 1$ such that 
$\tau_*(\widetilde{D}) = D$. 
\end{Example}

\begin{Example}
\label{eg:3}
This example shows that we need to replace $D$ with a divisor $E$ 
linearly equivalent to $D$ to satisfy the condition (C') in 
Theorem~\ref{thm:main} (see Remark~\ref{rmk:1:5}). 

Let $G$ be the following hyperelliptic graph of genus $4$, 
where each vertex is given by a white circle or a black circle. 
Let $\Gamma$ be the metric graph associated to $G$. 
The involution $\iota$ of $\Gamma$ is given by 
the reflection relative to the horizontal line through $w_2$. 

Let $D = 3 [v_1] + [v_2]$. 
We take a function $f$ on $\Gamma$ so that
$f(v_1) = 1, f(w) = 0$ for any $w \in V(G) \setminus\{v_1\}$ 
and $f$ is linear on each edge. 
Then $D + (f) = [v_2] + [w_1] + [w_2] + [w_3]$.  
Since $ [v_2] + [w_1]  \sim [w] + [\iota(w)]$ 
for any $w \in \Gamma$ by Lemma~\ref{lemma:iota:iota}, 
we have $r_\Gamma(D) \geq 1$. In fact, it is easy to see from 
Theorem~\ref{thm:main:2} that 
$r_\Gamma(D) = 1$. 
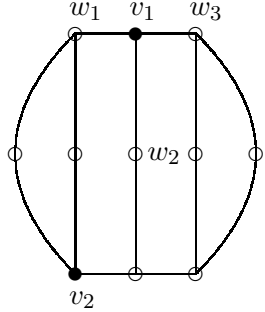
\[
  \setlength\unitlength{0.08truecm}
  \begin{picture}(50,60)(0,0)
  \put(10,50){\line(1,0){20}}
  \put(10,10){\line(1,0){20}}
  \put(10,10){\line(0,1){40}}
  \put(20,10){\line(0,1){40}}
  \put(30,10){\line(0,1){40}}
  \put(0, 30){\circle{2}} 
  \put(10, 30){\circle{2}}
  \put(20, 30){\circle{2}}
  \put(30, 30){\circle{2}} 
  \put(40, 30){\circle{2}} 
  \put(10, 50){\circle{2}}
  \put(20, 50){\circle*{2}}
  \put(30, 50){\circle{2}} 
\qbezier(0, 30)(0, 40)(10, 50)
\qbezier(0, 30)(0, 20)(10, 10)
\qbezier(40, 30)(40, 40)(30, 50)
\qbezier(40, 30)(40, 20)(30, 10)
  \put(10, 10){\circle*{2}}
  \put(20, 10){\circle{2}}
  \put(30, 10){\circle{2}} 
  \put(9, 5){$v_2$}
  \put(19, 53){$v_1$}
  \put(9, 53){$w_1$}
  \put(22, 29){$w_2$}
  \put(29, 53){$w_3$}
 \end{picture}
\]

The graph $G$ has no bridges. Let $\KK$ be a complete valuation field with ring integer $R$ and algebraically closed residue field $k$ such that ${\rm char}(k) \neq 2$.  
By Theorem~\ref{thm:lifting}, 
we have a regular, generically smooth, strongly semi-stable, totally degenerate $R$-curve $\Xscr$ with reduction graph $G = (G, \mathbf{0})$ 
such that the generic fiber $X$ is hyperelliptic. 
Let $\iota_X$ be the hyperelliptic involution on $X$. 
As we have shown, this model $\Xscr$ satisfies the condition (C') in 
the introduction. 

Let $P_1, P_2 \in X(\overline{\KK})$ be any points with 
$\tau(P_1) = v_1$ and $\tau(P_2) = v_2$.  Since 
$\tau\circ\iota_X = \iota\circ\tau$ and $\iota(v_1) \neq v_2$, 
we have $\iota_X(P_1) \neq P_2$. 
We set $\widetilde{D} = 3 P_1 + P_2$. By Proposition~\ref{prop:Yamaki}, 
we have $r_X(\widetilde{D}) = 0$. Hence 
$r_\Gamma(\tau_*(\widetilde{D})) \neq r_X(\widetilde{D})$. 
\end{Example}

\section{Rationality in lifting and a conjecture of Caporaso}
\label{sec:Caporaso}
In this section, we consider variants of the conditions (C) and (C') in the introduction, and discuss how they are related to the conjecture of Caporaso \cite[Conjecture~1]{Ca2}.  Finally, we show one direction of the conjecture 
for a hyperelliptic vertex-weighted graph satisfying the condition (i) 
in Theorem~\ref{thm:main}. 

\subsection{Terminology and properties of finite graphs}
In what follows, we consider divisors and linear equivalences on a finite graph $G$.  Let us first fix the notation and terminology. The {\em group of divisors} $\Div(G)$ on $G$ is defined to be the free $\ZZ$-module generated by the elements of $V(G)$. Then $\Div(G) = \bigoplus_{v \in V(G)} \ZZ[v]$ is naturally seen as a $\ZZ$-submodule of $\Div(\Gamma)$, where $\Gamma$ is the metric graph associated to $G$. 

A {\em rational function} on $G$ is a piecewise linear function on $\Gamma$, which is linear on edges and with integer value at each vertex. The set of rational functions on $G$ is denoted by $\Rat(G)$. Let $f \in \Rat(G)$. Then $f$ is naturally seen as an element of $\Rat(\Gamma)$, and $\zero(f) \in  \Div(\Gamma)$ is in fact an element of $\Div(G)$.  The set of {\em principal divisors} is defined by 
$\Prin(G) := \{ \zero(f) \mid f \in \Rat(G)\}$. Two divisors $D, E \in \Div(G)$ are said to be {\em linearly equivalent} in $\Div(G)$, and we write $D \sim_G E$, if $D - E \in  \Prin(G)$. Since $\Prin(G) = \Prin(\Gamma) \cap \Div(G)$,
we have, for $D, E \in \Div(G)$,  $D \sim_G E$ if and only if $D \sim E$. 

We will use the following lemma.  Recall that, by a hyperelliptic vertex-weighted graph $(G, \omega)$, we mean that $(\Gamma, \omega)$ is hyperelliptic, where  
$\Gamma$ the metric graph associated to $G$ (cf. Definition~\ref{def:hyp:vw:2}). 

\begin{Lemma}
\label{lemma:hyp:finite}
Let $(G, \omega)$ be a hyperelliptic vertex-weighted graph, and $\Gamma$ the metric graph associated to $G$. Then there exists a divisor $D \in \Div(G)$ with $\deg(D) = 2$ and $r_{(\Gamma, \omega)}(D) = 1$. 
\end{Lemma}

\Proof 
If $e$ is a leaf edge with a leaf end $v$ with $\omega(v) = 0$, then we contract $e$. Let $G^\prime$ be the finite graph that is obtained by contracting all such leaf edges, and give the vertex-weight function $\omega^\prime$ by the restriction of $\omega$ to $V(G^\prime)$.  

Let $\Gamma^\prime$ be the metric graph associated to  $G^\prime$. 
By Proposition~\ref{prop:HMG:2}, $\Gamma^\prime$ has the 
hyperelliptic involution $\iota^\prime: \Gamma^\prime \to 
\Gamma^\prime$ (see Definition~\ref{def:HMG:2}). We note that there exists a point $v \in \Gamma^\prime$ with $\omega(v) >0$ or $\val(v) \neq 2$. Then 
$v$ and $\iota^\prime(v)$ are both vertices of $G^\prime$. We set $D := [v] + [\iota^\prime(v)]$, which is seen 
as an element of $\Div(G)$. Then we have $\deg(D) = 2$ and $r_{(\Gamma, \omega)}(D) = 1$. 
\QED

\subsection{Conditions (F) and (F'), and a conjecture of Caporaso}
As before, let $\KK$ be a complete discrete valuation field with ring of integers $R$ and algebraically closed residue field $k$ such that ${\rm char}(k) \neq 2$. 
Let $(G, \omega)$ be a vertex-weighted graph, and let $\Gamma$ be the metric graph associated to $G$. Let $\Xscr$ be a regular, generically smooth, semi-stable $R$-curve with generic fiber $X$ and reduction graph $(G, \omega)$. For each vertex $v$ 
of $G$, let $C_v$ denote the irreducible component of the special fiber 
$\Xscr_0$ corresponding to $v$. 

Since $X$ is smooth (resp. $\Xscr$ is regular), the group of Cartier divisors 
on $X$ (resp. $\Xscr$) is the same as the group of Weil divisors. 
The Zariski closure of an effective divisor on $X$ in $\Xscr$ is a 
Cartier divisor. Extending by linearity, one can associate 
to any divisor on $X$ 
a Cartier divisor on $\Xscr$, which is also called the {\em Zariski closure} of the divisor. 

Let $\widetilde{D}$ be a divisor on $X$ and $\widetilde{\Dscr}$ the Zariski closure 
of $\widetilde{D}$. Let $\OO_{\Xscr}(\widetilde{\Dscr})$ be the locally-free sheaf on $\Xscr$ associated to $\widetilde{\Dscr}$. 
We set 
\[
  \rho_*(\widetilde{D}) :=  
  \sum_{v \in V(G)} \deg\left(
  \OO_{\Xscr}(\widetilde{\Dscr})\vert_{C_v}
  \right) [v]
  \in \Div(G). 
\]
We obtain the {\em specialization map} 
\begin{equation}
\label{eqn:specialization:map:fin}
  \rho_*: \Div(X) \to \Div(G). 
\end{equation}
We note that, if $\widetilde{D} \in \Div(X(\KK))$, i.e., 
$\widetilde{D} = \sum_{i=1}^k n_i P_i$ with $P_i \in X(\KK)$, 
then $\rho_*(\widetilde{D}) = \tau_*(\widetilde{D})$, where 
$\tau_*: \Div(X_{\overline\KK}) \to \Div(\Gamma)$ is the specialization map 
\eqref{eqn:specialization:map} induced by $\tau: X(\overline{\KK}) \to \Gamma$ 
in \eqref{eqn:specialization:map:0} (see \cite[\S2.3]{B}). 

Recall from the introduction that we consider the following condition (F), which is a variant of the condition (C). 
\begin{enumerate}
\item[(F)]
For any $D \in \Div(G)$, 
there exists a divisor $\widetilde{D} \in \Div(X)$ 
such that $D = \rho_*(\widetilde{D})$ and 
$r_{(\Gamma, \omega)}(D) = r_{X}(\widetilde{D})$. 
\end{enumerate}
We remark that the condition (F) is concerned with the existence of a lifting as a divisor over $\KK$ (not just as a divisor over $\overline{\KK}$) of a divisor $D$ on $G$ (not just on $\Gamma_\QQ$). We also consider the following condition (F'), which is a variant of the condition (C') in the introduction. 
\begin{enumerate}
\item[(F')]
For any $D \in \Div(G)$, 
there exist a divisor $E = \sum_{i=1}^k n_{i} [v_i] 
\in \Div(G)$ that is linearly equivalent to $D$ in $\Div(G)$, 
and $P_i \in X(\KK)$ for $1 \leq i \leq k$ such that 
$\tau(P_i) = v_i$ for any $1 \leq i \leq k$ and 
$r_{(\Gamma, \omega)}(E) = r_{X}\left(\sum_{i=1}^k n_{i} P_i\right)$. 
\end{enumerate}

Now we show Proposition~\ref{prop:Caporaso}, which is due to Caporaso. 

\begin{PropositionNoNum}[$=$ Proposition~\ref{prop:Caporaso}]
Let $\KK, R$ and $k$ be as above.  Let $(G, \omega)$ be a vertex-weighted graph, and let $\Gamma$ be the metric graph associated to $G$. Let $\Xscr$ be a regular, generically smooth, semi-stable $R$-curve with generic fiber $X$ and reduction graph $(G, \omega)$. Assume that $\Xscr$ satisfies the condition \textup{(F)}.  Then, 
for any divisor $D \in \Div(G)$, 
we have 
\[
  r_{(G, \omega)}^{\alg, k}(D) \geq r_{(\Gamma, \omega)}(D) .  
\]
\end{PropositionNoNum}

\Proof
Recall from the introduction that $r_{(G, \omega)}^{\alg, k}(D)$ is defined by   
\begin{align*}
  r_{(G, \omega)}^{\alg, k}(D)  & := \max_{X_0}\, r(X_0, D), \\
  r(X_0, D)  & := \min_{E}\, r^{\max}(X_0, E), \\
  r^{\max}(X_0, E) & := \max_{\Escr_0}\, \left( h^0(X_0, \Escr_0) -1 \right), 
\end{align*}
where $X_0$ runs over all connected reduced projective nodal curves defined over $k$ with dual graph $(G, \omega)$,  
$E$ runs over all divisors on $G$ that are linearly equivalent to $D$ in $\Div(G)$,  and $\Escr_0$ runs over all Cartier divisors on $X_0$ such that $\deg\left(\Escr_0\vert_{C_v}\right) = E(v)$ for any $v \in  V(G)$. 

Now we take $X_0$ as the special fiber of $\Xscr$. Let $E$ be any divisor on $G$ that is linearly equivalent to $D$ in $\Div(G)$. If the condition (F) is satisfied, then there exists $\widetilde{E} \in \Div(X)$ such that $\rho_*(\widetilde{E}) = E$ and  $r_{(\Gamma, \omega)}(E) 
= r_{X}(\widetilde{E})$. 

Let $\Escr$ be the Zariski closure of $\widetilde{E}$ in $\Xscr$, and we put 
$\Escr_0 := \Escr\vert_{X_0}$. By the definition of $\rho_*(\widetilde{E})$, 
we have $\deg\left(\Escr_0\vert_{C_v}\right) = E(v)$. On the other hand, 
the upper-semicontinuity of the cohomology implies 
that
\[
  h^0(X_0, \Escr_0) -1 \geq h^0(X, \widetilde{E}) - 1 
  = r_X(\widetilde{E}) =  r_{(\Gamma, \omega)}(E) = r_{(\Gamma, \omega)}(D).  
\]

Thus, letting $X_0$ be the special fiber of $\Xscr$, $E$ any divisor on $G$ that is linearly equivalent to $D$ in $\Div(G)$, and $\Escr_0$ the restriction of the Zariski closure of $\widetilde{E}$ to the special fiber, we obtain 
$r_{(G, \omega)}^{\alg, k}(D) \geq r_{(\Gamma, \omega)}(D)$. 
\QED

\subsection{Conditions (F) and (F') for hyperelliptic metric graphs}
We prove the following theorem, which is in a way refinement of Theorem~\ref{thm:main}. Theorem~\ref{thm:main:F} implies Theorem~\ref{thm:finite}. 

\begin{Theorem}
\label{thm:main:F}
Let $\KK$ be a complete discrete valuation field with ring of integers $R$ 
and algebraically closed residue field $k$ such that ${\rm char}(k) \neq 2$. 
Let $(G, \omega)$ be a hyperelliptic vertex-weighted graph. 
Then the following are equivalent. 
\begin{enumerate}
\item[(i)]
For every vertex $v$ of $G$, there are at most $(2\, \omega(v) + 2)$ positive-type bridges emanating from~$v$.  
\item[(ii)]
There exists a regular, generically smooth, semi-stable $R$-curve $\Xscr$ with reduction graph $(G, \omega)$ satisfying~\textup{(F)}. 
\item[(iii)]
There exists a regular, generically smooth, semi-stable $R$-curve $\Xscr$ with reduction graph $(G, \omega)$ satisfying~\textup{(F')}. 
\end{enumerate}
\end{Theorem}

\begin{Remark}
In the proof of Theorem~\ref{thm:main}, we see that the condition (i) in Theorem~\ref{thm:main:F}  is equivalent to the existence of a regular, generically smooth, semi-stable $R$-curve $\Xscr$ with generic fiber $X$ and reduction graph $(G, \omega)$ such that $X$ is hyperelliptic. Then any such $R$-curve $\Xscr$ satisfies the conditions (F) and (F') (and also (C) and(C')). 
\end{Remark}

\Proof
Let $g$ denote the genus of $(G, \omega)$. Let $\Gamma$ be the metric graph 
associated to $G$. 

\medskip
{\bf Step 1.}\quad
We show that (iii) implies (ii). 
By \cite[Corollary~A.9]{B}, the specialization map $\rho_*: \Prin(X) \to \Prin(G)$ is surjective. (In \cite{B}, a loopless finite graph is considered, and the general case is reduced to the case of a loopless finite graph.) Then arguing in exactly the same way as in Lemma~\ref{lemma:prime}, we find that (iii) implies (ii). 

\medskip
{\bf Step 2.}\quad
We show that (ii) implies (i). By Lemma~\ref{lemma:hyp:finite}, there exists a divisor $D \in \Div(G)$ such that $\deg(D) =2$ and $r_{(\Gamma, \omega)}(D) = 1$.  Then by the condition (F), there exists a divisor $\widetilde{D} \in \Div(X)$ with $\deg(\widetilde{D}) = 2$ and $r_{X}(\widetilde{D}) = 1$. Thus $X$ is a hyperelliptic curve, and by Theorem~\ref{thm:lifting},  the condition (i) holds.  

\medskip
{\bf Step 3.}\quad
We show that (i) implies (iii). This step is the main part of the proof of this theorem. 

We take a regular, generically smooth, semi-stable $R$-curve $\Xscr$ with reduction graph $(G, \omega)$ such that the generic fiber $X$ of $\Xscr$ is hyperelliptic as in the 
proof of Theorem~\ref{thm:1:3}. We are going to show that $\Xscr$ satisfies (F'). 

Let $\tau\vert_{X(\KK)} : X(\KK) \to V(G)$ be 
the restriction of the specialization map 
$\tau: X(\overline{\KK}) \to \Gamma$ to $X(\KK)$.  
Then $\tau\vert_{X(\KK)} : X(\KK) \to V(G)$ is surjective (see \cite[Remark~2.3]{B}). Note that $\tau(P) = \rho_*(P)$ for $P \in X(\KK)$, where $P \in X(\KK)$ is regarded as an element of $\Div(X(\KK)) \subset \Div(X)$ on the right-hand side. 

Let $D$ be any divisor on $G$. 

\smallskip
{\bf Case 1.}\quad
Suppose that $r_{(\Gamma, \omega)}(D) = -1$. We put $E := D$, and write 
$E = \sum_{i=1}^k n_i [v_i] \in \Div(G)$. By the surjectivity of $\tau\vert_{X(\KK)}$, 
we take $P_i \in X(\KK)$ such that $\tau(P_i) = v_i$ for $1 \leq i \leq k$. Then 
we have $r_{(\Gamma, \omega)}(D) = r_{X}\left(
\sum_{i=1}^k n_i P_i
\right)$ by a similar argument of the proof of Theorem~\ref{thm:1:3} (Case 1). 

\smallskip
{\bf Case 2.}\quad
Suppose that $r_{(\Gamma, \omega)}(D) \geq 0$. 
We follow the notation in the proof of Theorem~\ref{thm:1:3}. In particular, 
$(G^\prime, \omega^\prime)$ is the vertex-weighted graph obtained by contracting all the leaf edges of $G$ with leaf ends of weight zero, $\Gamma^\prime$ is the metric graph associated to $G^\prime$, and $\iota^\prime: \Gamma^\prime \to \Gamma^\prime$ is the hyperelliptic involution (cf. Definition~\ref{def:HMG:2}). 
Let $\varpi: \Gamma \to \Gamma^\prime$  be the retraction map, and 
$\jmath: \Gamma^\prime \hookrightarrow \Gamma$ be the natural embedding. 
By slight abuse of notation, we also write $\varpi: G \to G^\prime$ and 
$\jmath: G^\prime \to G$ for the induced maps on finite graphs. We regard 
$G^\prime$ as a subgraph of $G$. 

We take any $v \in V(G^\prime)$ such that $\iota^\prime(v) \in V(G^\prime)$ (cf. the proof of Lemma~\ref{lemma:hyp:finite}). By the surjectivity of $\tau\vert_{X(\KK)}$, we take $P \in X(\KK)$ with 
$\tau(P) = v$. We set $P^\prime := \iota_X(P) \in X(\KK)$ and 
$v^\prime := \tau(P^\prime) \in \Div(G)$. 
 Then we have $\varpi(v^\prime) = \iota^\prime(v)$, so that 
 $v^\prime \sim_G \iota^\prime(v)$.

We set
$r = p_{(\Gamma, \omega)}(D)$, and put 
\[
  F:= D - r \left([v] + [v^\prime]\right) \in \Div(G).
\]
Then $F \sim_G D - r \left([v] + [\iota^\prime(v)]\right)$. 

Let $\Gamma^\omega$ be the virtual weightless metric graph associated to 
$(\Gamma, \omega)$ and $\jmath^\omega: 
\Gamma \hookrightarrow \Gamma^\omega$ the natural embedding. 
Regarding $F$ as a divisor on $\Gamma$, we have 
\[
r_{(\Gamma, \omega)}(F) := 
r_{\Gamma^\omega}(\jmath^\omega_*(F)) = 
r_{\Gamma^\omega}\left(\jmath^\omega_*(D) - r \left([v] + [v^\prime]\right)\right) \geq 0
\]
by the definition of $p_{(\Gamma, \omega)}(D)$. 
By Lemma~\ref{lem:compatibility:reducedness}(3), we have $r_\Gamma(F) \geq 0$. By \cite[Lemma~2.3]{GK}, there exists an effective divisor on $G$ that is linearly equivalent to $F$. 
It follows that 
\[
  F \sim_G [u_1] + \cdots + [u_s]
\]
for some $u_1, \ldots, u_s \in V(G)$. By the surjectivity of $\tau\vert_{X(\KK)}$, 
we take $Q_j \in X(\KK)$ with $\tau(Q_j) = u_j$ for $j = 1, \ldots, s$. 
We find that $\iota_X(Q_i) \neq Q_j$ for $i \neq j$. Indeed, if $\iota_X(Q_i) 
= Q_j$, then $[u_i] + [u_j] \sim_G [\varpi(u_i)] + [\varpi(u_j)] \sim_G 
[v] + [\iota^\prime(v)]$. Then $|F - ([v] + [\iota^\prime(v)])| = |D - (r+1)([v] + [\iota^\prime(v)])| \neq \emptyset$,  which contradicts 
$r = p_{(\Gamma, \omega)}(D)$ (cf. \eqref{eqn:p:g:1}). 

We set $E := r \left([v] + [v^\prime]\right) + [u_1] + \cdots + [u_s] 
\in \Div(G)$ and $\widetilde{E} 
:= r (P + P^\prime) + Q_1 + \cdots + Q_s \in \Div(X(\KK))$. Then 
$\tau_*(\widetilde{E}) = E$. Further, 
$E$ is linearly equivalent to $D$, so that we have  
\[
 r_{(\Gamma, \omega)}(E) = \begin{cases}
  r & \textup{(if $\deg(D) - r \leq g$), } \\
  \deg(D) -g & \textup{(if $\deg(D) - r \geq g+1$)} 
  \end{cases} 
\]  
by Theorem~\ref{thm:main:3}. On the other hand, by Proposition~\ref{prop:Yamaki}, 
we have 
\[
  r_{X}(\widetilde{E}) = \begin{cases}
  r & \textup{(if $\deg(D) - r \leq g$), } \\
  \deg(D) -g & \textup{(if $\deg(D) - r \geq g+1$).} 
  \end{cases} 
\]
Hence we obtain $r_{(\Gamma, \omega)}(E) = r_{X}(\widetilde{E})$, and $\Xscr$ satisfies the condition (F'). 
\QED

\begin{CorollaryNoNum}[$=$ Corollary~\ref{cor:finite}]
Let $k$ be an algebraically closed field with ${\rm char}(k) \neq 2$. Let $(G, \omega)$ be a hyperelliptic graph such that for every vertex $v$ of $G$, there are at most $(2 \omega(v) +2)$ positive-type bridges emanating from $v$. Then, 
for any $D \in \Div(G)$, we have $r_{(G, \omega)}^{\alg, k}(D) \geq r_{(\Gamma, \omega)}(D)$.  
\end{CorollaryNoNum}

\Proof
We set $R := k[[t]]$ and $\KK := k(\!(t)\!)$, where $t$ is an indeterminate.  Then 
$\KK$ is a complete discrete valuation field with ring of integers $R$ and 
residue field $k$. It suffices to apply 
Proposition~\ref{prop:Caporaso} and Theorem~\ref{thm:main:F}.
\QED

\renewcommand{\theTheorem}{A.\arabic{Theorem}}
\renewcommand{\theClaim}{A.\arabic{Theorem}.\arabic{Claim}}
\renewcommand{\theequation}{A.\arabic{equation}}
\renewcommand{\thesubsection}{A.\arabic{subsection}}
\setcounter{Theorem}{0}
\setcounter{subsection}{0}
\setcounter{Claim}{0}
\setcounter{equation}{0}
\section*{Appendix. Deformation theory}
Let $\langle \iota\rangle$ denote the group of order $2$ with 
generator $\iota$. 
To prove Theorem~\ref{thm:lifting} in 
\S\ref{sec:hyp:graph}, we use the $\langle \iota\rangle$-equivariant 
deformation theory. Since we cannot find a suitable reference in the form we use in \S\ref{sec:hyp:graph} (i.e., over the ring of Witt vectors of a field $k$ of any characteristic $\neq 2$), we put together necessary results in this appendix. Note that one can find, among other things,  the $\langle \iota\rangle$-equivariant deformation theory over $k$ of characteristic $\neq 2$  
(i.e, not over the ring of Witt vectors) in Ekedahl \cite{E}. 
Unlike the previous sections, proofs 
of the results in this appendix are only sketched. 
Our basic references are \cite{DM, E, HaDef, Sch}. 

We fix the notation and terminology. 
Let $k$ be a field. We assume that ${\rm char}(k) \neq 2$. 
We put 
\[
\Lambda :=
\begin{cases}
k & \text{if ${\rm char}(k) = 0$,} \\
\text{the ring of Witt vectors over $k$}
&
\text{if ${\rm char}(k) > 0$.}
\end{cases}
\]
Let $\mathscr{A}$ be the category of Artin local $\Lambda$-algebras
with residue field $k$. 
Let $R$ be a complete local $\Lambda$-algebra with residue field $k$. 
Let $h_R: \mathscr{A} \to (Sets)$ 
be the functor given by $h_R(A) = \Hom(R, A)$ for $A \in \Ob(\mathscr{A})$. 
A functor $F : \mathscr{A} \to (Sets)$ is {\em pro-represented} by 
$R$ if $F$ is isomorphic to~$h_R$. 

Let $\mathscr{\widehat{A}}$ be the category of complete local $\Lambda$-algebras with residue field $k$. One can extend any functor $F : \mathscr{A} \to (Sets)$ to $\widehat{F} : \widehat{\mathscr{A}} \to (Sets)$ by defining $\widehat{F} (R) := \varprojlim F(R / \mathfrak{m}^{i})$, where $R \in \Ob(\widehat{\mathscr{A}})$ with maximal ideal $\mathfrak{m}$. If $F$ is pro-represented by $R$, then there is an isomorphism $\xi: h_R \to F$, and we can think of $\xi$ as an element of $\widehat{F}(R)$.  In this case, the pair $(R, \xi)$ is called the {\em universal family} of $F$. 

Let $F$ and $G$ be functors from  $\mathscr{A}$ to $(Sets)$. 
A morphism $G \to F$ is said to be {\em smooth} if for every surjective homomorphism $B \to A$ of local Artin $\Lambda$-algebras, the map $G(B) \to G(A) \times_{F(A)} F(B)$ is surjective. If $G \to F$ is smooth, then 
for every $A \in \Ob(\mathscr{A})$, the map $G(A) \to F(A)$ is surjective.  

It is useful to introduce a weaker notion of the pro-representability. Let $F: \mathscr{A} \to (Sets)$ be a functor. A pair $(R, \xi)$ with $R \in \widehat{\mathscr{A}}$ and $\xi \in \widehat{F}(R)$ is a {\em pro-representable hull} of $F$ if 
$h_R \to F$ is smooth and if the associated map $h_R(k[\epsilon]/(\epsilon^2)) \to F(k[\epsilon]/(\epsilon^2))$ is bijective. In this case, 
the pair $(R, \xi)$ is also called a {\em miniversal family} of $F$.  

\subsection{Equivariant deformation of curves}
\label{subsec:A:1}
In this subsection, we describe
the $\langle \iota \rangle$-equivariant deformation theory of  
curves.

Let $X_{0}$ be a stable curve of genus $g$ over $k$.
Let $A$ be an Artin local $\Lambda$-algebra with residue field $k$.
A \emph{deformation} of $X_{0}$ to $A$ is 
a stable curve 
$\mathcal{X} \to \Spec(A)$ with an identification
$\mathcal{X} \times_{\Spec(A)} \Spec(k) = X_{0}$.
Two deformations $\mathcal{X} \to \Spec(A)$ and 
$\mathcal{X}' \to \Spec(A)$ are said to be 
isomorphic if
there exists an 
isomorphism $\mathcal{X} \to \mathcal{X}'$ over $A$
which restricts to the identity on the special fiber $X_{0}$.

The \emph{deformation functor for $X_{0}$} is a functor 
\[
\Def_{X_{0}} : \mathscr{A} \to (Sets)
\]
that assigns to any $A \in \Ob ( \mathscr{A} )$ 
the set of isomorphism classes of deformations of $X_{0}$ to $A$. 

Suppose now that $X_{0}$ is a hyperelliptic stable curve of genus $g$ 
over $k$ (cf. Definition~\ref{def:hyp:ss}). 
For an Artin local $\Lambda$-algebra $A$ with residue field $k$, 
an \emph{$\langle \iota \rangle$-equivariant deformation
of $X_{0}$ to $A$} is 
the pair of a stable curve $\mathcal{X} \to \Spec(A)$ with an identification
$\mathcal{X} \times_{\Spec(A)} \Spec(k) = X_{0}$
and an $\langle \iota \rangle$-action on $\mathcal{X}$ whose restriction 
to the special fiber $X_{0}$ is the given $\langle \iota \rangle$-action.  
Two equivariant deformations $\mathcal{X} \to \Spec(A)$ 
and $\mathcal{X}' \to \Spec(A)$ of $X_{0}$ are  
said to be isomorphic if there is an 
$\langle \iota \rangle$-equivariant
isomorphism $\mathcal{X}' \to \mathcal{X}$ over $A$
whose restriction to the special fiber $X_{0}$ is the identity. 

The \emph{equivariant deformation functor} for $X_{0}$ is 
a functor 
\[
\Def_{(X_{0} , \iota)} : \mathscr{A} \to (Sets) 
\]
which assigns to $A \in \Ob(\mathscr{A})$  
the set of  isomorphism classes of equivariant deformations of
$X_{0}$ to~$A$. 

The deformation functor $\Def_{X_{0}}$ has a natural 
$\langle \iota \rangle$-action induced by the $\langle \iota \rangle$-action 
on $X_0$. We define $\Def_{X_{0}}^{\iota}$ to be the subfunctor of $\Def_{X_{0}}$ consisting of the $\langle\iota\rangle$-invariant elements of $\Def_{X_{0}}$.
We define a canonical morphism
$\Def_{(X_{0} , \iota)} \to \Def_{X_{0}}$ 
by forgetting the $\langle \iota \rangle$-action, which 
factors through $\Def_{X_{0}}^{\iota}$.

\begin{Lemma} \label{equivariant-invariant}
The canonical morphism $\Def_{(X_{0} , \iota)} \to \Def_{X_{0}}^{\iota}$
is an isomorphism.
\end{Lemma}

\Proof 
One can obtain the assertion by using \cite[Theorem~1.11]{DM}. 
\QED

\begin{Proposition} 
\label{universal-equivariant-deformation}
The functor $\Def_{(X_{0} , \iota)}$ is pro-represented by 
a formal power series over $\Lambda$.
\end{Proposition}

\Proof
The deformation functor $\Def_{X_{0}}$
is pro-represented by $\Spf \Lambda [[ t_{1} , \ldots , t_{3g-3} ]]$ 
by \cite[p.79]{DM}. 
Since $\Def_{(X_{0} , \iota)} = \Def_{X_{0}}^{\iota}$
by Lemma~\ref{equivariant-invariant},
$\Def_{(X_{0} , \iota)}$ can be pro-represented by the
formal subscheme of $\Spf \Lambda [[ t_{1} , \ldots , t_{3g-3} ]]$
consisting of the $\langle\iota\rangle$-invariants. 
Since the order $2$ of $\iota$ is invertible in $\Lambda$, one can take 
a suitable coordinate system such that  
the $\langle \iota \rangle$-action is expressed as
\[
\iota^{\ast} (t_{1}) = t_{1},
\ldots
,
\iota^{\ast} (t_{s}) = t_{s},
\iota^{\ast} (t_{s+1}) = - t_{s+1},
\ldots
,
\iota^{\ast} (t_{3g-3}) = - t_{3g-3}  
\]
for some $0 \leq s \leq 3g-3$. 
It follows that $\Def_{(X_{0} , \iota)}$ is a formal power series over $\Lambda$. 
\QED

\begin{Remark} \label{algebraizable}
Since the universal deformation
$\mathscr{C} \to
\Spf \Lambda [[ t_{1} , \ldots , t_{3g-3}]]$ is algebraizable (\cite[p.82]{DM}), 
the universal $\langle\iota\rangle$-equivariant deformation of $X_{0}$
is algebraizable. 
\end{Remark}

\subsection{Deformation of nodes with $\langle \iota \rangle$-actions}
\label{subsec:A:2}
In this subsection, 
we consider the deformation theory of nodes with 
$\langle \iota \rangle$-actions.

We begin by recalling the deformation theory of nodes.
Let $\OO \cong  k [[ x,y ]] / (xy)$ be a node over $k$. 
Let $A$ be an Artin local $\Lambda$-algebra with residue field $k$. 
A {\em deformation} of $\OO$ to $A$ is 
a co-cartesian diagram 
of local homomorphisms 
\begin{equation}
\label{eqn:cocartesian}
\begin{CD}
\OO @<<< B \\
@AAA @AAA \\
k @<<< A 
\end{CD}
\end{equation}
of $A$-algebras,
where $B$ is a flat local $A$-algebra. 
Two deformations $A \to B$ and $A \to B^\prime$ are said to be isomorphic 
if there exists an $A$-algebra isomorphism $B \to B^\prime$ which makes the 
co-cartesian diagrams for $B$ and $B^\prime$ commutative. 

Let $\mathscr{A}$ be the category of Artin local $\Lambda$-algebras
with residue field $k$ as in \S\ref{subsec:A:1}. 
The {\em deformation functor for $\OO$} is  
the functor
\[
\Def_{\OO} : \mathscr{A} \to (Sets) 
\]
that assigns to any $A \in \Ob(\mathscr{A})$ the
set of isomorphism classes of deformations
of $\OO$ to~$A$. 

The deformation functor $\Def_{\OO}$ has a pro-representable hull.
To be precise, 
by \cite[p.81]{DM}, 
\begin{align} \label{miniversal-deformation-node}
\begin{CD}
\OO = k [[ x , y  ]] / (xy) @<<<  \Lambda [[ x , y , t]] /( xy - t) \\
@AAA @AAA \\
k @<<< \Lambda [[ t ]] 
\end{CD}
\end{align}
is a pro-representable hull (i.e., a miniversal family) of $\Def_{\OO}$.

Suppose now that $\OO$ is equipped with an 
$\langle \iota \rangle$-action. 
Then we have an $\langle \iota \rangle$-action
$\iota_* : \Def_{\OO} \to \Def_{\OO}$ as follows. 
For $A \in \mathrm{Ob}( \mathscr{A} )$,
take any $\eta \in \Def_{\OO} (A)$ with a representative
\[
\begin{CD}
\OO @<{\alpha}<< B \\
@AAA @AAA \\
k @<<< A 
.
\end{CD}
\]
Then the diagram
\[
\begin{CD}
\OO @<{\iota \circ \alpha}<< B \\
@AAA @AAA \\
k @<<< A 
.
\end{CD}
\]
is also a deformation of $\OO$ to $A$. We define 
$\iota_{\ast} (\eta)$ is to be the isomorphism class of the above diagram.
We have $\iota_{*} ^{2} = \id$.

Typical examples of nodes with $\langle\iota\rangle$-actions 
arise from hyperelliptic stable curves.
Let $X_{0}$ be a hyperelliptic stable curve over $k$
with 
hyperelliptic involution $\iota_{X_0}$. Recall from the definition of a hyperelliptic stable curve (cf. Definition~\ref{def:hyp:ss}) that for any irreducible component $C$ of $X_{0}$ with $\iota (C) = C$, the $\langle \iota \rangle$-action restricted to $C$ is nontrivial.  Let $c$ be an $\iota_{X_0}$-fixed node. 
Then $\OO := \widehat{\OO_{X_{0} , c}}$ is a node 
equipped with the $\langle \iota \rangle$-action given by $\iota_{X_0}$. 
The following lemma concretely describes 
the $\langle \iota \rangle$-action on  $\OO$. 

\begin{Lemma} 
\label{lemma:action:on:node}
Let $\OO$ be a node equipped with the $\langle \iota \rangle$-action as above \textup{(}i.e., arising from a hyperelliptic stable curve\textup{)}. 
Then
there exists a $k$-algebra isomorphism $\OO \cong  k [[ x,y ]] / (xy)$
for which the $\langle \iota \rangle$-action on $k [[ x,y ]] / (xy)$ 
is given by either one of the following\textup{:}  
\addtocounter{Claim}{1}
\begin{align} \label{(1)}
& \iota (x) = y, \quad \iota (y) = x, \\
\addtocounter{Claim}{1}
\label{(2)}
& \iota (x) = -x, \quad \iota (y) = -y. 
\end{align}
\end{Lemma}

We remark that the above actions are ``admissible'' in the sense of 
Ekedahl \cite[Definition~1.2]{E}. 

In what follows, let $\OO$ be a node with an $\langle \iota \rangle$-action
as in Lemma~\ref{lemma:action:on:node},
and we identify  $\OO$ with $k [[ x,y ]] / (xy)$ via the above isomorphism. 

\begin{Lemma} \label{lem:trivial-action-on-node}
Let $\OO = k [[ x , y ]] / (xy)$ be the node over $k$ with the $\langle \iota \rangle$-action given by either \eqref{(1)} or  \eqref{(2)}.
Let $\iota_{*} : \Def_{\OO} \to \Def_{\OO}$ be the induced 
$\langle \iota \rangle$-action.
Then $\iota_{*} = \id$.
\end{Lemma}

\Proof
Let $A$ be an Artin local $\Lambda$-algebra with residue field $k$.
Take any element of $\Def_{\OO} (A)$ with a representative
\begin{align*} 
\begin{CD}
\OO = k [[ x,y]] / (xy) @<{\alpha}<< B \\
@AAA @AAA \\
k @<<< A .
\end{CD}
\end{align*}
Note that $\OO$ is equipped with the $\langle \iota \rangle$-action given by 
either \eqref{(1)} or  \eqref{(2)}.
To show that the $\langle \iota \rangle$-action on $\Def_{\OO} (A)$ is trivial, 
it is enough to define an $A$-involution on $\iota_B : B \to B$ such
that $\alpha \circ \iota_B = \iota \circ \alpha$.

We put an $\langle \iota \rangle$-action on $\Lambda [[ x , y , t]] /( xy - t)$
over $\Lambda [[ t ]]$ as follows.
If the $\langle\iota\rangle$-action on $\OO$ is given by \eqref{(1)},
then we
let $\iota: \Lambda [[ x , y , t]] /( xy - t) \to \Lambda [[ x , y , t]] /( xy - t)$ 
be the $\Lambda[[t]]$-algebra involution given by $\iota (x) = y$ and $\iota (y) = x$.  
If the $\langle\iota\rangle$-action on $\OO$ is given by \eqref{(2)},
then we
let $\iota: \Lambda [[ x , y , t]] /( xy - t) \to \Lambda [[ x , y , t]] /( xy - t)$ 
be the $\Lambda[[t]]$-algebra involution given by $\iota (x) = -x$ and $\iota (y) = -y$.  

Since (\ref{miniversal-deformation-node}) is a pro-representable hull
of $\Def_{\OO}$, we have the following commutative diagram
\begin{align*} 
\begin{CD}
\OO = k [[ x,y]] / (xy) @<{\alpha}<< B @<<< \Lambda [[ x , y , t]] /( xy - t) \\
@AAA @AAA @AAA \\
k @<<< A @<<< \Lambda [[ t ]],
\end{CD}
\end{align*}
where each square is co-cartesian.
Then the $\langle \iota \rangle$-action on $\Lambda [[ x , y , t]] /( xy - t)$ 
induces the $A$-algebra involution $\iota_B$ on
$B$ by co-cartesian product, which satisfies $\alpha \circ \iota_B = \iota \circ \alpha$.
Thus we obtain the assertion.
\QED

\subsection{Global-local morphism}
\label{subsec:A:3}
Let $X_{0}$ be a stable curve of genus $g$ over $k$,
and let $p_{1}, \ldots p_{t}$ be all the nodes of $X_{0}$.
We assume that any node is defined over $k$. To ease notation, 
we denote by $\Def_{p_{i}}$ the deformation functor $\Def_{\widehat{\OO_{X_0, p_i}}}$ for $\widehat{\OO_{X_0, p_i}}$. 

The {\em global-local} morphism is a morphism 
\[
\Phi^{gl} : \Def_{X_{0}} \to
\prod_{i=1}^{t} \Def_{p_{i}}
\]
that assigns to any deformation
$\mathcal{X} \to \Spec(A)$ of $X_{0}$ 
the deformation $A \to \widehat{\OO_{\mathcal{X} , p_{i}}}$ of
each node $\widehat{\OO_{X_{0} , p_{i}}}$ (cf. \cite[p.81]{DM}). 
The morphism $\Phi^{gl}$ is smooth by \cite[Prop.(1.5)]{DM}. 

We consider an $\langle \iota \rangle$-equivariant version of the
global-local morphism.
Assume that $X_{0}$ a hyperelliptic stable curve over $k$ with 
hyperelliptic involution $\iota = \iota_{X_0}$. 
Let $p_{1} , \ldots , p_{r}$ be the nodes of $X_{0}$ fixed by $\iota$,
and let $p_{r+1}, \ldots , p_{r+s}$ be nodes such that
$p_{r+1}, \ldots , p_{r+s}$, $\iota (p_{r+1}) , \ldots , \iota (p_{r+s} )$
are the distinct nodes that are not fixed by $\iota$. 
The {\em $\langle\iota\rangle$-equivariant global-local} morphism is a morphism 
\[
\Phi_{\iota}^{gl}
:
\Def_{(X_{0} , \iota)} \to
\prod_{i = 1}^{r} \Def_{p_{i}} \times
\prod_{i = r + 1}^{r+s} \Def_{p_{i}} 
\]
that assigns, to any $\langle\iota\rangle$-equivariant deformation
$\mathcal{X} \to \Spec(A)$ of $X_{0}$, 
the 
deformation 
$A \to \widehat{\OO_{\mathcal{X}, p_{i}}}$ 
of the node 
$\widehat{\OO_{X_0 , p_{i}}}$ 
for $1 \leq i \leq r+s$. 
Note that the target of $\Phi_{\iota}^{gl}$ is 
$\prod_{i = 1}^{r+s} \Def_{p_{i}}  = \prod_{i = 1}^{r} \Def_{p_{i}} \times
\prod_{i = r + 1}^{r+s} \Def_{p_{i}}$, and {\em not} $\prod_{i = 1}^{r} \Def_{p_{i}} \times \prod_{i = r + 1}^{r+s} \Def_{p_{i}} \times \prod_{i = r + 1}^{r+s} \Def_{\iota(p_{i})}$.  

The following proposition shows 
that the $\langle \iota \rangle$-equivariant global-local morphism $\Phi_{\iota}^{gl}$
is smooth, as in the case of the usual global-local morphism $\Phi^{gl}$.

\begin{Proposition} \label{smoothnessofglmorphism}
The morphism $\Phi_{\iota}^{gl}$ is smooth.
\end{Proposition}

\Proof
By Proposition~\ref{universal-equivariant-deformation}, $\Def_{(X_{0} , \iota)}$ is pro-represented by a formal power series over $\Lambda$. 
By \eqref{miniversal-deformation-node},
the pro-representable hull of 
$
\prod_{i = 1}^{r} \Def_{p_{i}} \times
\prod_{i = r + 1}^{r+s} \Def_{p_{i}} 
$ is a formal power series over $\Lambda$. 
Then by \cite[the~proof~of~Prop.(1.5)]{DM}, 
it suffices to show that $\Phi_{\iota}^{gl} ( k [ \epsilon ] / ( \epsilon^{2}))$ 
is surjective. 

To do that,
we regard $\Phi_{\iota}^{gl}$ as the restriction of $\Phi^{gl}$
to the subfunctors 
consisting of the $\langle\iota\rangle$-invariants
as we now explain. 
First, by Lemma~\ref{equivariant-invariant}, $\Def_{(X_0, \iota)}$ is regarded as 
the subfunctor consisting of the $\langle\iota\rangle$-invariants of 
$\Def_{X_0}$. 
Next, we focus on the targets of $\Phi_{\iota}^{gl}$ and $\Phi^{gl}$.
We consider the $\langle \iota \rangle$-action on
$
\prod_{i = 1}^{r} \Def_{p_{i}} \times 
\prod_{i = r + 1}^{r+s} \left( \Def_{p_{i}} \times \Def_{\iota (p_{i})} \right) 
$
given by $\eta \mapsto \iota_*(\eta)$ for $\eta \in \Def_{p_{i}}$ for $1 \leq i \leq r$ and $(\eta, \eta^\prime) \mapsto (\iota_*(\eta^\prime), \iota_*(\eta))$ for 
$(\eta, \eta^\prime) \in  \Def_{p_{i}} \times \Def_{\iota (p_{i})}$ for $r+1 \leq i \leq r+s$. 
Let 
\[
\Psi:
\prod_{i = 1}^{r} \Def_{p_{i}} \times
\prod_{i = r + 1}^{r+s} \Def_{p_{i}} 
\to
\prod_{i = 1}^{r} \Def_{p_{i}} \times
\prod_{i = r + 1}^{r+s} \left( \Def_{p_{i}} \times \Def_{\iota (p_{i})} \right) 
\]
be the morphism defined 
by the product of the identity morphisms
$\Def_{p_{i}} \to \Def_{p_{i}}$ 
for $1 \leq i \leq r$,
and the graph embeddings
$\Def_{p_{i}} \ni \eta \mapsto ( \eta , \iota_* (\eta)) 
\in \Def_{p_{i}} \times \Def_{\iota (p_{i})}$ of $\iota_{\ast}$
for $r+1 \leq i \leq r+s$.
For $1 \leq i \leq r$, the $\langle \iota \rangle$-action on 
$\Def_{p_i}$ is trivial by Lemma~\ref{lem:trivial-action-on-node}.
For $r+1 \leq i \leq r+s$, the morphism
$\Def_{p_{i}} \to \Def_{p_{i}} \times \Def_{\iota (p_{i})}$
is an isomorphism onto the subfunctor of
$\Def_{p_{i}} \times \Def_{\iota (p_{i})}$ consisting of the
$\langle\iota\rangle$-invariants. 
Thus 
$\prod_{i = 1}^{r} \Def_{p_{i}} \times
\prod_{i = r + 1}^{r+s} \Def_{p_{i}} $
is regarded via
$\Psi$
as the subfunctor 
of $\prod_{i = 1}^{r} \Def_{p_{i}} \times
\prod_{i = r + 1}^{r+s} \left( \Def_{p_{i}} \times \Def_{\iota (p_{i})} \right)$ 
consisting of the $\langle\iota\rangle$-invariants.

Through these identifications, 
$\Phi_{\iota}^{gl} ( k [ \epsilon ] / ( \epsilon^{2}))$ 
is regarded as the
restriction of $\Phi^{gl} (k [ \epsilon ] / ( \epsilon^{2}))$ 
to the $\langle\iota\rangle$-invariants. 
By \cite[Prop.(1.5)]{DM}, 
$\Phi^{gl}(k [ \epsilon ] / ( \epsilon^{2}))$ is surjective. 
Since $2$ is invertible in $k$,
the induced map between $\langle\iota\rangle$-invariants is also surjective,
so that 
$\Phi_{\iota}^{gl} ( k [ \epsilon ] / ( \epsilon^{2}))$ is surjective.
\QED

\begin{Corollary} \label{cor:surjectivitiy:EGLM}
For any $R \in \widehat{\mathscr{A}}$, 
$\widehat{\Phi^{gl}_\iota} ( R )$ is surjective.
\end{Corollary}

\Proof
The assertion follows from Proposition~\ref{smoothnessofglmorphism}
and \cite[Remark~2.4]{Sch}.
\QED

\end{document}